\documentclass[11pt,a4paper]{article}
\usepackage{amsmath,amssymb,amsthm,color,comment,pstricks,tikz}

\allowdisplaybreaks

\newtheorem{theorem}{Theorem}[section]
\newtheorem{thm}[theorem]{Theorem}

\newtheorem{cor}[theorem]{Corollary}
\newtheorem{lemma}[theorem]{Lemma}
\newtheorem{lem}[theorem]{Lemma}
\newtheorem{proposition}[theorem]{Proposition}
\newtheorem{propn}[theorem]{Proposition}
\newtheorem{prop}[theorem]{Proposition}
\newtheorem{definition}[theorem]{Definition}
\newtheorem{defn}[theorem]{Definition}

\newtheorem{remark}[theorem]{Remark}

\numberwithin{equation}{section}

\def\be{\begin{equation}}
\def\ee{\end{equation}}
\def\bes{\begin{equation*}}
\def\ees{\end{equation*}}

\def\dUS{{d^S_\sU}}

\def\dU{{d_\sU}}
\def\dE{{d_\infty}}
\def\BU{{B_\sU}}
\def\wt{\widetilde}
\def\df{{d_f}}
\def\dw{{d_w}}

\def\dSU{{d_\sU^\sS}}

 \def\sE {{\mathcal E}} \def\sF {{\mathcal F}}
  \def\sI {{\mathcal I}}
  \def\sL {{\mathcal L}}

\def\sS {{\mathcal S}} \def\sT {{\mathcal T}} \def\sU {{\mathcal U}}

 \def\bE {{\mathbf E}}

 \def\bN {{\mathbb N}} 
\def\bP {{\mathbf P}}  \def\bR {{\mathbb R}}

 \def\bZ {{\mathbb Z}}

\def\sm{\smallskip\noindent}

\def\ignore#1{}

\def\lam {\lambda} \def\Lam {\Lambda}
 \def\eps{\varepsilon}
\def\th{\theta} 

\def\Gam{\Gamma} \def\gam{\gamma}

\def\to {\rightarrow}

\def\pd {\partial}
\def\q{\quad} 
\def\dint{\int\kern-.6em\int}

\def \half {{\tfrac12}}
\def\fract{\tfrac}

\def\wt{\widetilde}

\def\be{\begin{equation}}
\def\ee{\end{equation}}
\def\bes{\begin{equation*}}
\def\ees{\end{equation*}}
\def\ba{\begin{align}}
\def\ea{\end{align}}
\def\xxea{\end{align}}
\def\bas{\begin{align*}}
\def\eas{\end{align*}}

\def\nn{\nonumber}

\def\proof{{ \sm {\em Proof. }}}
\def\qed{{\hfill $\square$ \bigskip}}

\definecolor{dgreen}{rgb}{0, 0.6, 0.1}
\definecolor{dblue}{rgb}{0, 0.0, 0.6}
\definecolor{vdblue}{rgb}{0,.08, 0.45}
\definecolor{dred}{rgb}{0.7, 0.0, 0.0}
\definecolor{vdblue}{rgb}{0,.08, 0.45}

\definecolor{purple}{rgb}{0.6, 0.0, 0.6}
\definecolor{mytext}{rgb}{0.1, 0.1, 0.1}

\def\Reff{R_{\rm eff}}

 \def\pP{{\mathbf P}}   \def\pE{{\mathbf E}}
 \def\bP{{\mathbf P}}   \def\bE{{\mathbf E}}
 \def\BE{{B_\infty}}

\let\OLDthebibliography\thebibliography
\renewcommand\thebibliography[1]{
  \OLDthebibliography{#1}
  \setlength{\parskip}{0pt}
  \setlength{\itemsep}{0pt plus 0.3ex}
}

\begin{document}

\font\titlefont=cmbx14 scaled\magstep1
\title{\titlefont \vspace{-30pt}Quenched and averaged tails of the heat kernel\\of the two-dimensional uniform spanning tree\footnote{Data sharing not applicable to this article as no datasets were generated or analysed during the current study.}}
\author{{M.~T.~Barlow}\thanks{Research partially supported by NSERC (Canada).},
D.~A.~Croydon\thanks{Research partially supported by
JSPS KAKENHI Grant Numbers
18H05832, 19K03540} and {T.~Kumagai}\thanks{Research partially supported by JSPS KAKENHI Grant
Number JP17H01093 and by the Alexander von Humboldt Foundation.}}

\maketitle


\begin{center}
\emph{This paper is dedicated to the memory of Harry Kesten,\\
a pioneer in the study of anomalous random walks in random media.}
\end{center}

\begin{abstract}
This article investigates the heat kernel of the two-dimensional uniform spanning tree. We improve previous work by demonstrating the occurrence of log-logarithmic fluctuations around the leading order polynomial behaviour for the on-diagonal part of the quenched heat kernel. In addition we give two-sided estimates  for the averaged heat kernel, and we show that the exponents that appear in the off-diagonal parts of the quenched and averaged versions of the heat kernel differ. Finally, we derive various scaling limits for the heat kernel, the implications of which include enabling us to sharpen the known asymptotics regarding the on-diagonal part of the averaged heat kernel and the expected distance travelled by the associated simple random walk.
\end{abstract}

\section{Introduction}

The focus of this article is the two-dimensional uniform spanning tree (UST), which is a random subgraph of $\mathbb{Z}^2$ that will henceforth be denoted by $\mathcal{U}$. Since the introduction of this object in \cite{Pem91}, considerable progress has been made in our understanding of the geometry of USTs (and, more generally, uniform spanning forests), see \cite{BLPS} for background. In this direction, a particularly useful viewpoint was provided by Wilson, who gave a construction of USTs via loop erased random walks (LERWs) \cite{Wil}. Indeed, the latter description was at the heart of Schramm's seminal work describing the subsequential scaling limits of two-dimensional LERW and $\mathcal{U}$ in terms of what is now called the Schramm-Loewner evolution (SLE) \cite{Schramm}, see also \cite{LSW}. In recent years, building on Lawler and Viklund's convergence result for the LERW in its natural parametrisation \cite{LV}, a more detailed picture of the scaling limit of $\mathcal{U}$ has been established \cite{BCK,HS}. And, closely related to this, properties of the simple random walk (SRW) on $\mathcal{U}$ have also been explored \cite{BCK, BM10, BM11}. The goal here is to provide further insight into the behaviour of the heat kernel (transition density) of the latter process.

Let us proceed to present some of the basic notation that will be used throughout the article. We will assume that the two-dimensional UST $\sU$ is built on a probability space $(\Omega, \sF, \pP)$; we write $\pE$ for the associated expectation. Note that, $\mathbf{P}$-a.s., $\sU$ is a one-ended tree with vertex set $\mathbb{Z}^2$ \cite{Pem91}. We write $\gam(x,y)$ for the unique self-avoiding path between $x,y \in \bZ^2$, and $\gam(x,\infty)$ for the unique infinite self-avoiding path started at $x$. By Wilson's algorithm (see \cite{Wil}, and the recollection of this at the start of Section \ref{sec:LERW}), $\gam(x,y)$ is equal in law to the loop erasure of a SRW started at $x$ and run until it hits $y$. We will denote by $\dU$ the intrinsic metric on the graph $\sU$, so that $\dU(x,y)$ is the length of the geodesic $\gam(x,y)$. We write $\mu_\sU$ for the measure on $\mathbb{Z}^2$ such that $\mu_\sU(\{x\})$ is given by the number of edges of $\sU$ that contain $x$; this is the invariant measure of the simple random walk. We denote balls in the intrinsic metric $\dU$ by $\BU(x,r) =\{ y \in \bZ^2: \dU(x,y) \le r\}$. We use $d_\infty$ to denote the $\ell_\infty$ metric on $\bZ^2$, and $B_\infty(x,r)$ to denote balls in the $d_\infty$-metric; these balls are of course boxes.

Many of the exponents that describe the behaviour of $\sU$ and the associated random walk can be expressed in terms of the {\em growth exponent of the two-dimensional LERW}, which is given by $\kappa=5/4$. More precisely, let $L_n$ be the loop erasure of a SRW in $\bZ^2$ run until its first exit from $[-n,n]^2$, $M_n$ be the number of steps in $L_n$, and $G(n) = \bE(M_n)$. By \cite{Law14}, we have
\be \label{e:growth}
  G(n) = \bE M_n \asymp n^\kappa,
 \ee
where $\asymp$ means `bounded above and below by constant multiples of'. (This improves earlier estimates in \cite{Ken, Mas09}, which establish that $\lim_{n\rightarrow\infty} \log G(n)/\log n = \kappa$.) The papers \cite{BM10, BM11} gave estimates for the heat kernel of $\sU$ in terms of the function $G$; these can now be written more simply using \eqref{e:growth}. When we cite results from \cite{BM10, BM11} we will give the simplified versions without further comment.

Next, we introduce the simple random walk on $\sU$, which is the discrete-time Markov process $X^\sU=((X^\sU_n)_{n\geq0},({P}_x^\sU)_{x\in \mathbb{Z}^2})$ that at each time step jumps from its current location to a uniformly chosen neighbour in the graph $\sU$. For $x\in \mathbb{Z}^2$, the law ${P}_x^\sU$ is called the {\em quenched} law of the simple random walk on $\sU$ started at $x$, and we write
\[p_n^\sU(x,y)=\frac{P_x^\sU\left(X^\sU_n=y\right)}{\mu_\sU\left(\{y\}\right)},\qquad \forall x,y\in\mathbb{Z}^2,\]
for the corresponding quenched heat kernel.

To understand the properties of random walk on a space such as $\sU$, a by now well-established approach is to first study volume growth and resistance growth (see, for example, \cite{BJKS, Kum, KM}). Regarding the volume growth, one would expect from \eqref{e:growth} and Wilson's algorithm that $\BU(x,r^\kappa)$ should be approximately the same as $B_\infty(x,r)$, and hence that $|\BU(x,R)|$ should be of order $R^{2/\kappa}$. This expectation was confirmed  by \cite[Theorem 1.2]{BM11}, which gives stretched exponential estimates for the upper and lower tails of $R^{-2/\kappa} |B_\sU(0,R)|$. We define  the `fractal dimension' of $\sU$ by
\begin{equation}\label{dfdef}
d_f = \frac2\kappa=\frac85.
\end{equation}
Using the estimates in \cite[Theorem 1.2]{BM11} an easy Borel-Cantelli argument gives that there exist deterministic constants $c_1,c_2\in(0,\infty)$ such that, $\mathbf{P}$-a.s.,
\[c_1r^{d_f}(\log\log r)^{-9}\leq\mu_\sU\left(B_\mathcal{U}(0, r)\right)\leq c_2r^{d_f}(\log\log r)^{3}\]
for large $r$. The first main result of this paper is that volume fluctuations of log-logarithmic magnitude really do occur.

{\thm\label{mainthm1} $\mathbf{P}$-a.s.,
\begin{equation}\label{bigvolas1}
\limsup_{r \to \infty} \frac{   \mu_\sU\left(B_\mathcal{U}(0, r) \right) } {r^{d_f}(\log\log r)^{1/5}} =\infty,
\end{equation}
and also
\begin{equation}\label{smallvolas1}
\liminf_{r \to \infty}
\frac{ (\log\log r)^{3/5} \mu_\sU\left(B_\mathcal{U}(0, r) \right) }{r^{d_f}} =0.
\end{equation}}
\medskip

Similar fluctuations have also been observed for Galton-Watson trees \cite[Proposition 2.8]{BK} (see also \cite[Lemma 5.1]{CrKu}). The proof here is more complicated as the correlations between different parts of the space are harder to control. The key ingredient is the argument of Section \ref{sec:control} below, in which we provide a general technique for estimating from below the probability of seeing a particular path configuration in the initial stages of the construction of the UST via Wilson's algorithm. This enables us to control the probability of seeing especially short or long paths in some region of $\sU$.

The volume  fluctuations of Theorem \ref{mainthm1} are associated with corresponding fluctuations in the on-diagonal part of the quenched heat kernel. From \cite[Theorem 4.5]{BM11}, we know there exist deterministic constants $c_1,c_2\in(0,\infty)$ and $\alpha_1,\alpha_2\in(0,\infty)$ such that, $\mathbf{P}$-a.s.,
\[c_1n^{-d_f/d_w}(\log\log n)^{-\alpha_1}\leq p^\mathcal{U}_{2n}(0,0)\leq c_2n^{-d_f/d_w}(\log\log n)^{\alpha_2}\]
for large $n$. Here
\begin{equation}\label{dwdef}
d_w = 1 + d_f = \frac{2+\kappa}{\kappa} =\frac{13}5
\end{equation}
is the so-called walk dimension; this represents the space-time scaling exponent with respect to the intrinsic metric. Applying Theorem \ref{mainthm1}, we are able to deduce that log-log fluctuations in the quenched heat kernel actually occur.

{\cor\label{cor1} There exists $\beta>0$ such that,  $\mathbf{P}$-a.s.,
\begin{align*}
  \liminf_{n\to \infty}(\log\log n)^{1/13}n^{d_f/d_w}p^\mathcal{U}_{2n}(0,0)&=0, \\
  \limsup_{n\to \infty}(\log\log n)^{-\beta} n^{d_f/d_w}p^\mathcal{U}_{2n}(0,0)&=\infty.
\end{align*}}
\medskip

These volume and heat kernel fluctuations arise from unlikely configurations of $\sU$ inside $B_\infty(0,r_k)$ at a (random) sequence of scales  $r_k \to \infty$. Another consequence of the occurrence of such exceptional configurations is the failure of the elliptic Harnack inequality in this setting. For a precise description of the particular form of the elliptic Harnack inequality that we consider, see Definition \ref{ehidef} below.

{\cor\label{ehicor} The large-scale elliptic Harnack inequality does not hold for the random walk on $\sU$.}
\medskip

We now consider the off-diagonal heat kernel. To avoid the issues of parity that arise from the fact $\sU$ is a bipartite graph, we introduce the following smoothed version of the heat kernel
\[ \tilde{p}_n^\sU(x,y):=\frac{p_n^\sU(x,y)+p_{n+1}^\sU(x,y)}{2}.\]
In \cite[Theorem 4.7]{BM11} it was shown that there exist deterministic constants $\alpha,C\in (0,\infty)$ such that, $\mathbf{P}$-a.s.:
\[\frac{n^{-\frac{d_f}{d_w}}}{A} \exp\left\{-A  \left(\frac{d_\sU(0,x)^{{d_w}}}{n}\right)^{\frac{1}{{d_w}-1}} \right\}\leq\tilde{p}_n^\sU(0,x)\leq {A}{ n^{-\frac{d_f}{d_w}}}\exp\left\{-\frac{1}{A}  \left(\frac{d_\sU(0,x)^{ {d_w}}}{n}\right)^{\frac{1}{{d_w}-1}} \right\}\]
holds whenever $n \ge d_\mathcal{U}(0,x)$ and $\max\{n^\dw,|x|\}$ is suitably large, where
\be \label{e:Adef}
A= A(n,x) :=C  \left(\log\left(\max\{n^\dw,|x|\} \right)\right)^\alpha.
\ee
The logarithmic correction factor $A$ represents the possible influence of exceptional environments on the heat kernel.

We are unlikely to see an exceptional configuration at any particular scale, so it is not surprising that for the averaged heat kernel the fluctuations of Corollary \ref{cor1} disappear: by \cite[Theorem 4.4(c)]{BM11}, we have that
\be \label{e:annondub}
c_1n^{-d_f/d_w} \leq \pE \, p^\mathcal{U}_{2n}(0,0)\leq c_2n^{-d_f/d_w},\qquad \forall n\geq 1.
\ee
As for the off-diagonal part of the averaged heat kernel, one might hope that one could replace the random distance $d_\sU(0,x)$ with its typical order with respect to the Euclidean metric, that is, $|x|^\kappa$, and that, as with  \eqref{e:annondub} one would be able to remove the errors associated with the term $A$ in \eqref{e:Adef}. We show that this is almost the case, however, in the annealed off-diagonal bounds the exponent $\frac{1}{d_w-1}$ needs to be replaced by a strictly smaller number.

\begin{thm}\label{mainthm3} There exist constants $c_1,c_2,c_3,c_4\in (0,\infty)$ and
$0< \theta_2 \le \theta_1 <1$ such that: for every $x=(x_1,x_2)\in\mathbb{Z}^2$ and $n\geq |x_1|+|x_2|$,
\[
n^{-\frac{d_f}{d_w}} \exp\left\{-c_2  \left(\frac{|x|^{\kappa {d_w}}}{n}\right)^{\frac{\theta_1}{{d_w}-1}} \right\}\leq
 \pE \tilde{p}_n^\sU(0,x) \le   c_3 n^{-\frac{d_f}{d_w}}\exp\left\{-c_4  \left(\frac{|x|^{\kappa {d_w}}}{n}\right)^{\frac{\theta_2}{{d_w}-1}} \right\}. \]
\end{thm}\medskip

Our argument indicates that we can take $\theta_1<1$ due to contribution to the averaged heat kernel from realisations of $\sU$ where the intrinsic distance from $0$ to $x$ is unusually short, and thus where the heat kernel $\tilde{p}^{\sU}_n(0,x)$ is unusually large. This phenomenon was not observed in the earlier study of random walk on a Galton-Watson tree of \cite{BK} (see Theorem 1.5 in particular), since the intrinsic metric of the trees was the only one involved there.

\begin{remark}
{\rm We have $\th_1 = \frac{ \dw-1 }{\kappa \dw -1}=\frac{32}{45}$, and we conjecture that this is also the correct value for $\th_2$. This would mean that the averaged heat kernel estimates of Theorem \ref{mainthm3} are of the usual sub-Gaussian form, but with respect to the extrinsic walk dimension $\kappa d_w$, rather than the intrinsic walk dimension that appears in the quenched bounds.
}\end{remark}

In the course of our proofs we obtain some new tail estimates on the length of the path $\gam(x,y)$ between points $x$ and $y$; by Wilson's algorithm this is also the length of a LERW run from $x$ to $y$.

\begin{thm} \label{T:dxy}
(a) There exist constants $c_i$ such that for $\lam \ge 1$, $x, y \in \bZ^2$,
\[ c_1 e^{ -c_2 \lam^4 } \le  \pP\big( \dU(x,y)  < \lam^{-1} d_\infty(x,y)^\kappa \big) \le c_3 e^{ -c_4 \lam^4 }.\]
(b) There exist constants $c,q$ such that  for $\lam \ge 1$, $x, y \in \bZ^2$,
\[ \pP\left(\dU(x,y) \ge \lam d_\infty(x,y)^\kappa\right) \le c ( \log \lam )^q\lam^{-(2-\kappa)/\kappa}. \]
\end{thm}\medskip

The upper bound in (a) is proved in Theorem \ref{T:LERW-lb}, (b) is proved at the end of Section \ref{sec:LERW}, and the lower bound in (a)  is proved at the end of Section \ref{sec:control}.

\smallskip
We now consider the scaling limit of the UST and its heat kernel. Schramm's original work encoded $\mathcal{U}$ in terms of a path ensemble (consisting of the shortest paths in $\mathcal{U}$ between pairs of vertices), which enabled basic topological properties of any possible scaling limit to be deduced. In \cite{BCK}, building on the work of \cite{BM10, BM11}, this scaling picture was extended to incorporate the intrinsic (i.e.\ shortest path) metric on $\mathcal{U}$, as well as the uniform measure, with the result of \cite{BCK} being expressed in terms of the tightness under rescaling of $\mathcal{U}$ in a certain Gromov-Hausdorff-type topology for metric-measure spaces with an embedding into Euclidean space. The main obstacle to extending the work of \cite{BCK} to a full (i.e.\ non-subsequential) scaling limit was the need to prove the existence of the scaling limit of the two-dimensional LERW as a stochastic process, rather than simply as a compact subset of the plane. This was subsequently established in \cite{LV}, and Holden and Sun \cite{HS} then proved that $\sU$ has a full scaling limit as a metric-measure space.

Let us now describe the setting of \cite{BCK} more precisely. To retain information about $\sU$ in the Euclidean topology, $(\sU,\dU)$ can be considered as a spatial tree (cf.\ \cite{DuL}) -- that is, as a real tree (see \cite[Definition 1.1]{rrt}, for example) obtained from the graph by including unit line segments along edges, embedded into $\mathbb{R}^2$ via a continuous map  $\phi_\sU:\sU\rightarrow \mathbb{R}^2$, which is given by the identity on vertices, with linear interpolation along edges. In addition, suppose the space is rooted at the origin of $\mathbb{Z}^2$, giving a random `measured, rooted spatial tree'  $(\sU,\dU,\mu_\sU,\phi_\sU, 0)$. For this quintuplet, it follows from \cite[Theorem 1.1]{HS} (see also \cite[Theorem 1.1]{BCK}) that
\begin{equation}\label{scaling}
\left(\sU,\delta^{\kappa}\dU,\delta^{2}\mu_\sU,\delta\phi_\sU,0\right)\buildrel{d}\over{\rightarrow}\left(\sT,d_\sT,\mu_\sT,\phi_\sT,\rho_\sT\right)
\end{equation}
as $\delta\rightarrow0$ with respect to the Gromov-Hausdorff-type topology introduced in \cite[Section 3]{BCK}. The random limit space is such that, $\mathbf{P}$-a.s.: $(\mathcal{T},d_\mathcal{T})$ is a complete and locally compact real tree; $\mu_\mathcal{T}$ is a locally finite Borel measure on $(\mathcal{T},d_\mathcal{T})$; $\phi_\mathcal{T}$ is a continuous map from $(\mathcal{T},d_\mathcal{T})$ into $\mathbb{R}^2$; and $\rho_\mathcal{T}$ is a distinguished vertex in $\mathcal{T}$. In the original result of \cite{BCK}, the measure $\mu_\sU$ considered was the uniform measure on the vertices, but it is no problem to replace this with the measure we consider here since, after scaling the uniform measure by a factor of two, the Prohorov distance between the two measures is bounded above by two, and so the discrepancy disappears in the scaling limit. Moreover, it readily follows from \eqref{scaling} that the space $(\sT,d_\sT,\mu_\sT,\phi_\sT,\rho_\sT)$ satisfies the following scale invariance property: for any $\lambda>0$,
\begin{equation}\label{scaleinvar}
\left(\sT,\lambda^{\kappa}d_\sT,\lambda^{2}\mu_\sT,\lambda\phi_\sT,\rho_{\mathcal{T}}\right)\buildrel{d}\over{=}\left(\sT,d_\sT,\mu_\sT,\phi_\sT,\rho_\sT\right).
\end{equation}
Further from the SLE description of the limit in \cite{HS}, one also has rotational invariance, i.e.\ for any $\theta\in[0,2\pi)$,
\begin{equation}\label{rotinvar}
\left(\sT,d_\sT,\mu_\sT,R_\theta\phi_\sT,\rho_{\mathcal{T}}\right)\buildrel{d}\over{=}\left(\sT,d_\sT,\mu_\sT,\phi_\sT,\rho_\sT\right),
\end{equation}
where $R_\theta$ is a rotation of Euclidean space about the origin by the angle $\theta$. In Proposition \ref{reroot} below, we further establish an invariance under a rerooting property for the limit space.

One of the motivations for proving \eqref{scaling} was to show that the random walks on $\sU$ converge to a limiting process. It was shown in \cite{BCK} that the random walks on $\sU$ started from 0 satisfy
\begin{equation}\label{rwconvh}
\left(\delta X^\mathcal{U}_{\delta^{-\kappa d_w}t}\right)_{t\geq 0}\rightarrow \left(\phi_{\sT}\left(X^\mathcal{T}_{t}\right)\right)_{t\geq 0}
\end{equation}
in distribution under the averaged or annealed law. (Cf.\ the more general statements concerning the convergence of random walks on trees of \cite{ALW,Cr}.) Here, for $\mathbf{P}$-a.e.\ realisation of $(\sT,d_\sT,\mu_\sT,\phi_\sT,\rho_\sT)$, $X^\mathcal{T}=(X^\mathcal{T}_t)_{t\geq 0}$ is the canonical diffusion, or Brownian motion, on $(\sT,d_\sT,\mu_\sT)$ started from $\rho_\sT$, and $\phi_\sT(X^\mathcal{T})$ is the corresponding random element of $C(\mathbb{R}_+,\mathbb{R}^2)$. In this article, we connect the heat kernel of the discrete process $X^\sU$ to that of $X^\sT$, for which off-diagonal estimates were given in \cite{BCK}. As our first result in this direction, we show the convergence of the quenched and averaged on-diagonal part of the heat kernel. From \eqref{scaling} and \eqref{rwconvh}, the first claim of the following result, which concerns $(p^\sT_t(x,y))_{x,y\in\sT,t>0}$, the quenched heat kernel on the tree $\sT$ (as defined in \cite{BCK}), is essentially an application of the local limit theorem of \cite{CrHa}. To adapt this to yield the corresponding statement for the averaged heat kernels, we check the uniform integrability of the on-diagonal part of the discrete heat kernel by applying an argument similar to that applied to deduce averaged heat kernel estimates for Galton-Watson trees in \cite[Theorem 1.5]{BK}. We note that the exact form of the on-diagonal part of the limiting averaged heat kernel is a simple consequence of the scale invariance property \eqref{scaling}. Moreover, the result at \eqref{odhkscale} improves part of \cite[Theorem 4.4]{BM11}, where it was shown that $n^{{d_f}/{d_w}}\pE \tilde{p}_{n}^\sU(0,0)$ is bounded above and below by constants.

\begin{thm}\label{mainthm2} It holds that
\[\left(n^{{d_f}/{d_w}} \tilde{p}_{\lfloor tn\rfloor }^\sU(0,0)\right)_{t>0}\buildrel{d}\over{\rightarrow}\left( p_{t}^\sT(\rho_\sT,\rho_\sT)\right)_{t>0}\]
in distribution with respect to the topology of uniform convergence on compact subsets of $(0,\infty)$, and moreover,
\begin{equation}\label{odhkscale}
\left(n^{{d_f}/{d_w}}\pE \tilde{p}_{\lfloor tn\rfloor }^\sU(0,0)\right)_{t>0} \rightarrow\left(\pE p_{t}^\sT(\rho_\sT,\rho_\sT)\right)_{t>0}=\left( Ct^{-d_f/d_w}\right)_{t>0}
\end{equation}
in the same topology, where $C\in(0,\infty)$ is a constant.
\end{thm}
\medskip

We next turn our attention to the off-diagonal part of the heat kernel. Whilst it is natural to ask whether the scaling limit of \eqref{odhkscale} can be extended to include the off-diagonal part, we recall that $\phi_{\mathcal{T}}$ is not a bijection (see \cite[Theorem 1.3]{BCK}), and so one cannot \emph{a priori} assume that the limit of $n^{{d_f}/{d_w}}\pE \tilde{p}_{\lfloor tn\rfloor }^\sU(0,[xn^{1/\kappa d_w}])$ (where we write $[ xn^{1/\kappa d_w}]$ for the closest lattice point to $xn^{1/\kappa d_w}$) can be written as $\pE p_{t}^\sT(\rho_\sT,\phi_\sT^{-1}(x))$, or indeed that this latter expectation is well-defined.
This being the case, the following result is presented in terms of the density of the embedded process $\phi_\sT(X^\sT)$, where $X^\sT$ is the canonical Brownian motion on the limiting space; we note $\phi_\sT(X^\sT)$ is not Markov under the annealed law (or strong Markov under the quenched law, see Remark \ref{zqrem} below).
Nevertheless, as we will show, $\phi_\sT^{-1}$ is well defined except on a set of Lebesgue
measure zero, and so  the averaged density of $\phi_\sT(X^\sT)$
is in fact given by the expression
$\pE p_{t}^\sT(\rho_\sT,\phi_\sT^{-1}(x))$.
 The key additional input to the proof of this result is an equicontinuity property for the discrete heat kernel under scaling (see Proposition \ref{equicont}), which in turn depends on our estimate for the probability of seeing long paths in the uniform spanning tree (see Theorem \ref{T:dxy}).

\begin{thm}\label{denslimit}
For each $t\in(0,\infty)$, $\phi_\mathcal{T}(X^\mathcal{T}_t)$ admits a continuous probability density
${q}_t=({q}_t(x))_{x\in\mathbb{R}^2}$ under the annealed probability law $\pP \cdot P^\sU_0$, so that
\[\pE\left(P^\sU_0( \phi_\mathcal{T}(X^\mathcal{T}_t) \in B )\right) = \int_B q_t(x) dx\]
for all Borel $B \subseteq \bR^2$. The functions $({q}_t)_{t>0}$  satisfy the following.\\
(a) There exists a constant $C \in(0,\infty)$ such that
\[  |{q}_t(x)-{q}_t(y)| \leq C t^{-\df/2\dw} |x-y|^{\kappa/2},\qquad \forall x,y\in \bR^2, t >0.\]
(b) For any $\lambda>0$ and $\theta\in[0,2\pi)$, it holds that
\[\left(\lambda^{d_f/d_w}q_{ t\lambda}\left(\lambda^{\frac{1}{\kappa d_w}}R_\theta x\right)\right)_{x\in\mathbb{R}^2}=\left(q_{ t}(x)\right)_{x\in\mathbb{R}^2}.\]
(c) For each $t\in(0,\infty)$, it holds that
\[\left(n^{{d_f}/{d_w}}\pE \tilde{p}_{\lfloor tn\rfloor }^\sU(0,[ x n^{\frac{1}{\kappa d_w}}])\right)_{x\in\mathbb{R}^2} \rightarrow
\left(q_t(x)\right)_{x\in\mathbb{R}^2}\]
uniformly on compact subsets of $\mathbb{R}^2$.\\
(d) For each $t\in(0,\infty)$ and $x\in\mathbb{R}^2$, it holds that
\[q_t(x)=\pE\left( p_{t}^\sT(\rho_\sT,\phi_\sT^{-1}(x))\right).\]
\end{thm}
\medskip

\begin{remark}
{\rm By (b) there exists a continuous function $f: \bR_+ \to \bR_+$ such that
\[  q_t(x)=   t^{-\df/\dw} f( |x| t^{-1/\kappa \dw} ). \]
}\end{remark}
\medskip

As an almost immediate corollary of Theorems \ref{mainthm3} and \ref{denslimit}, we obtain the following.

\begin{cor}\label{cor2} There exist constants $c_1,c_2,c_3,c_4\in (0,\infty)$ and $\theta_1,\theta_2\in(0,1)$ such that the averaged density of $\phi^\mathcal{T}(X^\mathcal{T}_t)$, as given by Theorem \ref{denslimit}, satisfies: for every $x\in\mathbb{R}^2$ and $t>0$,
\[ c_1 t^{-{d_f}/{d_w}} \exp\left\{-c_2  \left(\frac{|x|^{\kappa {d_w}}}{t}\right)^{\frac{\theta_1}{{d_w}-1}} \right\}\leq
 {q}_t(x) \le   c_3 t^{-{d_f}/{d_w}}\exp\left\{-c_4  \left(\frac{|x|^{\kappa {d_w}}}{t}\right)^{\frac{\theta_2}{{d_w}-1}} \right\}. \]
\end{cor}
\medskip

As a further consequence of Theorem \ref{mainthm3}, together with the
the estimate on the probability of seeing long paths of Theorem \ref{T:dxy}, we
obtain the following upper bounds on
the averaged behaviour of the distance travelled by $X^\sU$ up to a given time,
both in terms of the Euclidean and the intrinsic distances.
Bounds for  $\pE \big( E^\sU_0 \big ( \dU (0,X^\sU_n )^p  \big)$ were considered
in \cite[Theorem 4.6]{BM11}, but the upper bound there has an additional
term $(\log n)^{c p}$.

\begin{cor}\label{cor-dist}  For every $p>0$, it holds that for $n \ge 1$,
\begin{align*}
&  c'_p n^{p/{\kappa d_w}} \le      \pE \Big( E^\sU_0 \big |X^\sU_n  \big |^p   \Big) \le c_p n^{p/{\kappa d_w}} , \\
&   c'_p n^{p/{d_w}} \le    \pE \Big( E^\sU_0 \big ( \dU (0,X^\sU_n ) ^p   \Big)   \le c_p n^{p/{ d_w}} .
\end{align*}
\end{cor}

It follows from the argument used to establish the random walk convergence result of \cite[Theorem 1.4]{BCK} that, under the averaged distribution, not only do we have \eqref{rwconvh}, but also $(n^{-1/{\kappa d_w}}|X^\sU_{\lfloor tn\rfloor}|)_{t\geq0}\buildrel{d}\over{\rightarrow}(|\phi_\sT(X^\sT_{t})|)_{t\geq0}$ and $(n^{-1/{d_w}}d_\mathcal{U}(0,X^\sU_{\lfloor tn\rfloor}))_{t\geq0}\buildrel{d}\over{\rightarrow}(d_\mathcal{T}(\rho_\sT,X^\sT_{t}))_{t\geq0}$.
Combining this with the integrability given by Corollary \ref{cor-dist} we obtain the following convergence
result.

\begin{cor}\label{cor3} (a) For every $p>0$, it holds that
\[\left(n^{-p/{\kappa d_w}} \pE \left( E^\sU_0 \left|X^\sU_{\lfloor tn\rfloor }\right|^p\right)\right)_{t\geq 0}\rightarrow
\left(\pE \left( E^\sT_{\rho_\sT} \left(\left|\phi_\sT(X^\sT_{t})\right|^p\right)\right)\right)_{t\geq 0}=\left(C_pt^{p/\kappa d_w}\right)_{t\geq 0}\]
with respect to the topology of uniform convergence on compact subsets of $[0,\infty)$, where $C_p\in(0,\infty)$ is a constant depending only upon $p$.\\
(b) For every $p>0$, it holds that
\[\left(n^{-p/{d_w}} \pE \left( E^\sU_0 \left(d_\mathcal{U}\left(0,X^\sU_{\lfloor tn\rfloor }\right)^p\right)\right)\right)_{t\geq 0}\rightarrow
\left(\pE \left( E^\sT_{\rho_\sT} \left(d_\mathcal{T}\left(\rho_\sT,X^\sT_{t}\right)^p\right)\right)\right)_{t\geq 0}=\left(C_pt^{p/d_w}\right)_{t\geq 0}\]
with respect to the topology of uniform convergence on compact subsets of $[0,\infty)$, where again $C_p\in(0,\infty)$ is a constant depending only upon $p$.
\end{cor}
\medskip

The remainder of the article is organised as follows. In Section \ref{sec:LERW} we review and refine some previous estimates for LERWs and the two-dimensional UST, proving the upper bound of Theorem \ref{T:dxy}(a) and Theorem \ref{T:dxy}(b) in particular. Section \ref{sec:control} provides an approach to showing that particular anomalous paths occur within the UST. This allows us to check the remaining part of Theorem \ref{T:dxy}, as well as the volume and heat kernel fluctuation results of Theorem \ref{mainthm1} and Corollary \ref{cor1} respectively, which will be done in Section \ref{sec:fluct}. Section \ref{sec:VRest} adapts results of \cite{BCK} concerning the structure of the UST to the case where we condition on a particular path being present, and these preliminary statements are then applied in Section \ref{sec:ann-b} to deduce the heat kernel bounds of Theorem \ref{mainthm3}. Then, in Section \ref{failEHI}, we confirm the failure of the elliptic Harnack inequality, as stated in Corollary \ref{ehicor}. And, in Section \ref{sec:scaling}, we apply the random walk scaling limit result of \cite{BCK} in conjunction with the estimates of this article to deduce Theorems \ref{mainthm2} and \ref{denslimit}, as well as Corollaries \ref{cor2} -- \ref{cor3}. Finally, we postpone to the appendix the proofs of some estimates from Section \ref{sec:LERW} that are relatively close variations on the proofs of the corresponding results in \cite{BM10}. NB.\ We will often use a continuous variable in places where a discrete one is required; in this case we implicitly mean that the floor of the relevant variable should be considered.

\section{Loop erased random walk and the UST} \label{sec:LERW}

This section contains some refinements of previous estimates on the geometry of the UST and the behaviour of the LERW. The key input we need for the averaged heat kernel upper bound (Proposition \ref{P:PF1n}) is a relatively straightforward adaptation of \cite[Proposition 2.10]{BCK}, adding resistance estimates to the volume estimates of the latter result. We also set out some new results, which include the upper bounds of Theorem \ref{T:dxy}.

We begin by introducing some notation for paths and operations on paths. A \emph{path} $\gam$ is a (finite or infinite) sequence of adjacent vertices in $\bZ^2$, i.e.\ $\gam=(\gam_0, \gam_1, \dots )$ with $\gam_{i-1}\sim \gam_i$, where for $x,y\in \mathbb{Z}^2$ we write $x\sim y$ if $|x-y|=1$. Given a set $A \subseteq \bZ^2$, we define $\tau_A = \min\{i \ge 0:\gam_i \notin A \}$, and set $\sE_A (\gam) = (\gam_0, \dots,  \gam_{\tau_A})$. Given a finite path $\gam$, we write $\sL (\gam)$ for the chronological loop erasure of $\gam$, see \cite{La1, Law99}.

We now recall Wilson's algorithm, see \cite{Wil}. For $x \in \bZ^2$ let $S^x$ be a simple random walk (SRW) on $\bZ^2$ started at $x$; we take $(S^x)_{x\in\mathbb{Z}^2}$ to be independent. Write $\bZ^2$ as a sequence $\{ z_0,z_1, z_2 , \dots\}$, and define a sequence of trees as follows:
\begin{align}
 \sU_0 &= \{ z_0\},\nonumber \\
 \sU_i &= \sU_{i-1} \cup \mathcal{L}(\sE_{\sU_{i-1}^c}(S^{z_{i}})), \q i \ge 1, \label{wilsonalg}\\
 \sU &= \cup_i \sU_i.\nonumber
\end{align}
By \cite{Wil}, the random tree $\sU$ has the law of the UST. It follows that the law of $\sU$ does not depend on the particular sequence $(z_i)$. In fact, more is true: the $z_i$ can be chosen adaptively as a function of $\sU_{i-1}$. However, if we use independent $(S^x)$ as above, then the final tree $\sU$ depends as a random variable on both the random walks and the sequence $(z_i)$. To circumvent this, for a finite graph Wilson  \cite{Wil} defined
a family of random variables (called `stacks') that enable one to define (non-independent) random walks $(\wt S^x)_{x\in\mathbb{Z}^2}$, and in this setup the final tree $\sU$ does not depend on $(z_i)$ for $i \ge 1$.
(Note though that the random walk $\wt S^{(z_n)}$ does depend on the sequence
$(z_1, \dots z_{n-1})$.) It is straightforward to check that this also holds with probability 1 for a recurrent graph, and it will sometimes be useful for us to apply this construction.

We write $L(x,\infty)$ for the loop-erased random walk from $x$ to infinity; this is the weak limit as $m \to \infty$ of $\sL(\sE_{B_m(x)}(S^x))$. By Wilson's algorithm $L(x,\infty)$ has the same law as the $\gam(x, \infty)$, the unique injective path from $x$ to infinity in $\sU$. (See \cite[Proposition 14.1]{BLPS}, for example.) We moreover write $\gam_x = \gam(x,\infty)$, $\gam_x[i]$ for the $i$th point on $\gam_x$, define the segment of the path $\gam_x$ between its $i$th and $j$th points by $\gam_x[i,j]= (\gam_x[i], \gam_x[i+1], \dots, \gam_x[j])$, and define $\gam_x[i, \infty)$ in a similar fashion. Furthermore, we let $\tau_{y,r}(\gam_x) = \min\{i\geq 0: \gam_x[i] \not\in B_\infty(y,r)\}$; whenever we use notation such as $\gam_x[\tau_{y,r}]$, the exit time $\tau_{y,r}$ will always be for the path $\gam_x$. For $x,y\in \mathbb{Z}^2$, we introduce the `Schramm distance' on $\sU$ (after \cite{Schramm}) by setting
\[\dSU(x,y):=\mathrm{diam}(\gam(x,y)),   \]
where the right-hand side is the diameter of $\gam(x,y)$ (considered as a subset of $\mathbb{Z}^2$) with respect to $d_{\infty}$.

In the next sequence of results of this section we collect and refine some properties of loop-erased random walks from \cite{BCK, BM10, BM11, Mas09}. In the following result, we write $\pd  B_{\infty}(0,r)$ for the outer boundary of $ B_{\infty}(0,r)$, i.e.\ those vertices of $\mathbb{Z}^2\backslash  B_{\infty}(0,r)$ that have as a neighbour a vertex in $ B_{\infty}(0,r)$.

\begin{lemma} \label{L:abscty} Let $\th>1$, $n\ge 1$, and suppose that $D_1$, $D_2$ are subsets of $\bZ^2$ with $\BE(0,\th n) \subseteq D_1 \cap D_2$. There exists a constant $c_1=c_1(\theta)$ such that if $\gam$ is a self-avoiding path from $0$ to $\pd \BE(0,n)$ then
\[ \pP\left(\sE_{\BE(0,n)}( \sL( \sE_{D_1}(S^0))) = \gam\right) \le c_1  \pP\left( \sE_{\BE(0,n)} (\sL( \sE_{D_2}(S^0))) =\gam\right).  \]
If for $i=1$ or $i=2$ one has $D_i=\bZ^d$, then $\sL( \sE_{D_i}(S^0))$ should be taken to be $L(0,\infty)$.
\end{lemma}
\proof See \cite[Proposition 4.4]{Mas09} for the result when $\th\ge 4$. Checking the proof of the latter result, one finds that the result as stated above holds for any $\th>1$. (In \cite{Mas09} the emphasis was on the fact that one can take $C(\theta) = 1 + c (\log \th)^{-1}$ for large $\th$.)
\qed

\begin{definition}\label{regdef}
{\rm Let $D\subseteq \bZ^2$. Let $\lam>1$, $1 \le r_1 \le r_2$.
We say that $D$ is {\em $(\lam, r_1, r_2)$-regular} if we have for each $x,y\in D$,
\begin{align*}
  \lam^{-1} \dSU(x,y)^\kappa\le \dU(x,y) \le\lam \dSU(x,y)^\kappa,&\qquad\hbox{when }r_1 \le \dSU(x,y) \le r_2,\\
 \dU(x,y) \le \lam r_1^\kappa,&\qquad\hbox{when }\dSU(x,y) \le r_1,\\
      \dU(x,y) \ge  \lam^{-1} r_2^\kappa,&\qquad\hbox{when }\dSU(x,y) \ge r_2.
\end{align*}
}\end{definition}

It is straightforward to check the following.

\begin{lemma}
(a)  If $D$  is $(\lam, r_1, r_2)$-regular and $r_2 \ge \lam^{2/\kappa} r_1$, then
\begin{align*}
  \lam^{-1} \dU(x,y) \le \dSU(x,y)^\kappa \le \lam \dU(x,y)^\kappa,
    &\qquad\hbox{when }   \lam r_1^\kappa\leq  \dU(x,y) \leq \lam^{-1} r_2^\kappa,\\
 \dSU(x,y)^\kappa \le \lam r_1^\kappa,  &\qquad\hbox{when } \dU(x,y) \le  \lam r_1^\kappa, \\
      \dSU(x,y) \ge  \lam^{-1} r_2^\kappa, &\qquad\hbox{when } \dU(x,y) \ge \lam^{-1} r_2^\kappa.
\end{align*}
(b) Let $\sT \subset \bZ^2$ be a tree, $w \not\in \sT$, $\gam(w,\sT)$ be a self-avoiding path from $w$ to $\sT$, and let $\sT'= \sT \cup \gam(w,\sT)$. If $\sT$ and $\gam(w,\sT)$ are  $(\lam, r_1, r_2)$-regular, then $\sT'$ is $( 2^\kappa \lam, 2r_1, r_2)$-regular.
\end{lemma}

The next result is a consequence of \cite[Proposition 2.8]{BCK}.

\begin{lemma} \label{L:regbox}
Let $n\ge1$ and $\lam \ge \lam_0$. Then
\[ \bP\left( B_\infty(0,n) \hbox{ is $(\lam, e^{-c_1 \lam^{1/2} } n,n)$-regular}\right) \ge 1  - c_2 e^{-c_3 \lam^{1/2} }. \]
\end{lemma}

\begin{lemma} \label{L:regpath} Let $\th >1$,
$n\ge 1$, and suppose $D \subseteq \bZ^2$ is such that $\BE(0,\th n) \subseteq D$. It then holds that there exist constants $c_i=c_i(\theta)$ such that for $\lam>\lam_0$, where $\lam_0$ is some large, finite constant,
$$ \pP\big( \sE_{\BE(0,n)} ( \sL ( \sE_D (S^0))) \hbox{ is $(\lam, ne^{-c_1 \lam^{1/2}},n)$-regular} \big)
\ge 1 - c_2 e^{-c_3 \lam^{1/2} }. $$
\end{lemma}
\proof By Lemma \ref{L:abscty} it is enough to prove this when $ \sL ( \sE_D (S^0))$ is replaced by $L(0,\infty)$. The bound then follows from Lemma \ref{L:regbox} and Wilson's algorithm. \qed

\begin{lemma} \label{L:gam0y}
There exists $\lam_0 \ge 1$ and constant $c_1>0$ with the following properties. Let $n\ge 1$, $\lam \ge \lam_0$ and $x \in B_\infty(0, 3n/4)$,  let $\pi$ be a shortest path in $\bZ^2$ between $0$ and $x$, and set $A =\{ y \in \bZ^2: \dE(y,\pi) \le n/8 \}$. Then
$$ \pP\big( \gam(0,x) \subseteq A \hbox{ and is  $(\lam, ne^{-c_1 \lam^{1/2}},n)$-regular}   \big) \ge c_1.  $$
\end{lemma}

\proof Let $G_1$ be the event that $B_\infty(0,n)$ is  $(\lam, ne^{-c_1 \lam^{1/2}},n)$-regular. Choose $k \ge 16$ and let $y$ satisfy $d_\infty(0,y) = n/k$. Let $G_2 = \{ \gam(0,y) \subseteq \BE(0, n/8) \}$; by \cite[Lemma 2.6]{BM11} we have $\pP( G_2^c )\le c_2 k^{-1/3}$; choose $k$ so that $c_2  k^{-1/3} \le \half$. Let $S^x$ be a SRW started at $x$, and $G_4$ be the event that $S^x$ makes a closed loop around 0 which separates $0$ and $y$ before it leaves $A$; we have $\bP(G_4) \ge c_3 >0$. Then $\bP( G_2 \cap G_4) \ge \half c_3$. We now choose $\lam_0$ large enough so that $\bP( G_1^c) \le \frac14 c_3$ and hence writing $G= G_1 \cap G_2 \cap G_4$ we have $\bP( G)  \ge \frac14 c_3$. On the event $G$ the SRW $S^x$ hits $\gam(0,y)$ before it exits $A$, so $\gam(0,x) \subseteq A$. Since $B_\infty(0,n)$ is regular, so is the path $\gam(0,x)$. \qed

\begin{thm} \label{T:LERW-lb}
Let $n\ge 1$, and suppose $D \subseteq \bZ^2$ is such that $\BE(0, n) \subseteq D$. If $D \neq \bZ^2$, set $L_{n,D}= \sE_{\BE(0,n)} (\sL(\sE_{D}(S^0)))$, and set $L_{n,\bZ^2}=  \sE_{\BE(0,n)} (L(0,\infty))$.  It then holds that there exist constants $c_i$ such that, for $\lam\geq 1$,
\be \label{e:LERW-lb}
\pP\left( | L_{n,D}| < \lam^{-1} n^\kappa \right) \le c_1 e^{ -c_2 \lam^{1/(\kappa-1)} }.
\ee
In particular, for any $x,y\in\mathbb{Z}^2$,
\be \label{e:duxy-lb}
  \pP\left( \dU(x,y)  < \lam^{-1} d_\infty(x,y)^\kappa \right) \le c_1 e^{- c_2 \lam^{1/(\kappa-1)} }.
\ee
\end{thm}

\proof A bound with exponent $\lam^{4/5 - \eps}$  is given in \cite[Theorem 6.7]{BM10}, and with some more care one can obtain  \eqref{e:LERW-lb} by essentially the same arguments -- see the Appendix for details. Taking $x=0$, $D= \bZ^2\backslash\{y\}$ and $n =\lfloor  d_\infty(0,y)\rfloor$, we have $  \dU(x,y) = | \sL(\sE_{D}(S^0))| \ge |L_{n,D}|$, which (with translation invariance) gives \eqref{e:duxy-lb}.
\qed

To state our next result, Proposition \ref{P:PF1n}, we need to introduce some more notation and basic definitions. Specifically, we write $\Reff$ for the effective resistance on $\sU$ considered as an electrical network with unit conductances along each edge. (See \cite{Bar17, LP} for background.) We recall from \eqref{dfdef} and \eqref{dwdef} the definition of $\df$ and $\dw$.

\begin{definition}\label{def:defjlam}
{\rm We say a ball $B_\sU(x,r)$ is $\lam$-good if we have the following:\\
(1) $\lam^{-1} r^{d_f} \le |B_\sU(x,r) | \le \lam r^{d_f}$, \\
(2) $\Reff(x, B_\sU(x,r)^c) \ge r/\lam$, \\
(3) $B_\sU(x,r) \subseteq \BE(x, \lam r^{1/\kappa})$.}\end{definition}

\noindent
We moreover define
\begin{equation}\label{f1def}
 F_1(\lam,n) =\{ B_\sU(x,r) \hbox{ is $\lam$-good for all } x \in \BE(0,n), \, e^{-\lam^{1/40}} n^\kappa \le r \le n^\kappa\},
 \end{equation}
and note that on the event $F_1(\lam,n)$ we have $B_\sU(x,n^\kappa) \subseteq \BE(x,\lam n)$ for all $x \in \BE(0,n)$.

\begin{propn}\label{P:PF1n}
There exist constants $c_1,c_2,\lambda_0$ such that
\[\pP\left(F_1(\lambda, n)^c\right) \le c_1\exp (-c_2\lambda^{1/16}),\qquad \forall n\ge 1,\:\lambda\ge \lam_0.\]
\end{propn}

\proof
The proof below is a modification of that of \cite[Proposition 2.10]{BCK}. Let $r=ne^{-\lam^{1/32}}$, and assume first that $n$ is large enough so that $r\ge \lam$. Let $J(x,\lambda)$ be the set of those $r\in [1,\infty) $ such that the three conditions in Definition \ref{def:defjlam} hold.

Set $R_1=n$, $R_2= r e^{\lam^{1/16}}$, and let $D_2$ be as in \cite[Proposition 2.9]{BCK}, with $|D_2 | \leq c\lam^4e^{2\lam^{1/16}}$. Set $m_0:=\inf\{m: m \ge e^{\kappa\lam^{1/32}} \}$, and let $ E(r,\lam):=\cap_{x\in D_2}\cap_{m=1}^{m_0+1} \{ m r^{\kappa}\in J(x,\lam)\}$. A simple union bound allows us to deduce from \cite[Theorem 1.1(a) and Proposition 4.2(a)]{BM11} that
\[\mathbf{P}\left(E(r,\lam)^c\right) \le |D_2| e^{\kappa\lam^{1/32}} c e^{-c'\lam^{1/9}}\le C  e^{-c'' \lam^{1/9}}.\]
Let $A_5(r,\lam)$ be the event given in the statement of \cite[Proposition 2.9]{BCK}; we have $\mathbf{P}( E(r,\lam)^c \cup A_5(r,\lam)^c) \le C \exp(-c \lam^{1/16} )$. Moreover, if $E(r,\lam)\cap A_5(r,\lam)$ holds, then, by \cite[(2.14) and the last display in the proof of Proposition 2.9]{BCK}, for each $x \in \BE(0,n)$, there exists $y=y_x \in D_2$ with $d_\mathcal{U}(x,y)\leq 4r^{\kappa}/\lam^{1/4}$ and $d_\infty(x,y)\le 2r/\lam$. Choosing $\lam_0$ large enough, we have $d_\mathcal{U}(x,y_x)\leq r^{\kappa}$ and $d_{\infty}(x,y_x) \le r$.

Now, suppose  $E(r,\lam)\cap A_5(r,\lam)$  holds, and let $x\in \BE(0,n)$, and $s \in [ 4 \lam  r^\kappa , n^\kappa]$. We will prove that $s \in J(x,2\lam)$ by verifying the conditions (1)--(3) in Definition  \ref{def:defjlam}. Choose $m \in[4\lam,m_0+1]$ so that $(m-1) r^\kappa \le s \le mr^\kappa$. It then holds that
\[|B_\sU \left(x,s \right)| \le |B_\sU(y_x, (m+1)r^{\kappa})| \le \lam ((m+1)/(m-1))^\df s^\df \le 2\lam s^\df.  \]
Similarly, $|B_\sU(x,s^{\kappa})|\geq (2 \lam)^{-1} s^\df$, so that the volume bound (1) holds. Next, applying the triangle inequality for resistances (and the fact that $\Reff(x,y_x)=d_\sU(x,y_x)$),
\begin{align*}
R_{\rm{eff}} \left(x,B_\sU(x,s)^c \right)
&\geq R_{\rm{eff}}\left(y_x ,B_\sU(x,s)^c \right) - d_\sU(x,y_x) \\
&\geq  R_{\rm{eff}} (y_x , B_\sU(y_x, s - r^\kappa)^c ) - r^\kappa\\
&\ge \lam^{-1} (m-2) r^\kappa  - r^\kappa  \geq (2 \lam)^{-1} s,
\end{align*}
which gives (2). Finally for (3) we have
\[B_\sU(x,s) \subseteq B_\sU\left(y_x, (m+1) r^{\kappa}\right) \subseteq \BE(y_x, {\lambda (m+1)^{1/\kappa} r} ) \subseteq \BE(x,{2\lam s}),\]
where the last inclusion holds since $r ( 1+ \lam (m+1)^{1/\kappa} ) \le 2 \lam (m-1)^{1/\kappa} r \le 2\lam s$. Thus, for $\lam\geq\lam_0$ with $\lam_0$ suitably large, $E(r,\lam)\cap A_5(r,\lam) \subseteq F_1(2 \lambda,n)$, and this completes the proof of the proposition in the case when $r \ge \lam$.

Finally, suppose $r =ne^{-\lam^{1/32}}<\lam$, i.e.\ $n\leq \lam e^{\lam^{1/32}}$. A union bound then gives
\[\pP\left(F_1(\lambda, n)^c\right)\leq \sum_{ x \in B_n(0)}\sum_{s=1}^{\lceil n^\kappa\rceil}\pP\left(s\notin J(x,\lambda)\right)\leq cn^{2+\kappa}e^{-c'\lambda^{1/9}}\leq Ce^{-c''\lambda^{1/9}},\]
where the last inequality is again an application of \cite[Theorem 1.1(a) and Proposition 4.2(a)]{BM11}. This is enough to complete the proof.
\qed

A further observation that will be useful in the proof of the  averaged heat kernel upper bound is the following.

\begin{lemma} \label{L:F1dist}
Suppose that $F_1(\lam,n)$ occurs, and let $x,y \in B_\infty(0,n)$. It then holds that
$d_{\infty}(x,y)\in[\lam e^{- \lam^{1/40}/\kappa} n,n]$ implies $\dU(x,y)  \ge \lam^{-\kappa} d_{\infty}(x,y)^\kappa$.
\end{lemma}
\proof Let $r <  \lam^{-\kappa} d_{\infty}(x,y)^\kappa$, so that $y \not\in B_\infty(x,{\lam r^{1/\kappa}})$. The condition on $d_{\infty}(x,y)$ implies that we can choose $r$ so that $r \in [e^{-\lam^{1/40}} n^\kappa, n^\kappa]$, and thus property (3) in the definition of a good ball implies that $y \not\in B_\sU(x,r)$, and so $\dU(x,y) >r$. \qed

The next few results will lead to the proof of Theorem \ref{T:dxy}(b), beginning with the case when $x$ and $y$ are neighbours in $\mathbb{Z}^2$.

\begin{prop} \label{nnp} There exist constants $c_i,q$ such that
\[c_1s^{-(2-\kappa)/\kappa}\leq \pP\left(\dU(0,e_1) \ge s \right) \le c_2 ( \log s )^qs^{-(2-\kappa)/\kappa},\]
for all $s\geq 2$, where $e_1=(1,0)$.
\end{prop}

We start with a proof of the lower bound. For this, it will be convenient to introduce $\sU'$, the dual of $\sU$. This is the graph with vertex set $(\mathbb{Z}+\frac12)^2$ whose edges are precisely those nearest neighbour edges that do not cross an edge of $\sU$. It is known that $\sU'$ has the same law as $\sU$ (see \cite{BLPS}). We set $0'=(1/2,1/2)$ for the root of the dual graph.
\bigskip

\noindent
\emph{Proof of the lower bound of Proposition \ref{nnp}.} Applying \cite[Proposition 2.8]{BCK} and Proposition \ref{P:PF1n}, for $r\geq 1$, $\lam\geq \lam_0$, we can find an event $G_0(\lam,r)$ with $\pP( G_0(\lam,r)^c ) \le c_1 e^{-c_2 \lam^{1/40}}$ such that if this event holds then we have $B_{\infty}(0,R)$ is $(c_3\lambda^{1/20},r,R)$-regular and also $F_1(\lam,R)$ holds, where $R:=re^{\lam^{1/40}/\kappa}$. Let $G'_0(\lam,r)$ be the corresponding event for the dual graph, and define $G(\lam,r) = G_0(\lam,r) \cap G'_0(\lam,r)$.

On $G(\lam,r)$, if $\gam'$ is the unique injective path from $0'$ to infinity in $\sU'$, then we have that the section of $\gam'$ from its last exit from $B_{\infty}(0',r)$ to $B_{\infty}(0',2r)^c$ has length greater than $c\lambda^{-1/20}r^\kappa$. (We assume that $\lam$ is chosen large enough so that $e^{\lam^{1/40}} \ge 2$.) Denote by $\gam'_r$ this section of $\gam'$, and let $\{x_1', x_2'\}$ be an edge crossed by $\gam'_r$. If $\{x_1,x_2\}$ is the dual edge to $\{x_1',x_2'\}$, then it must be the case that $\dUS(x_1,x_2) \ge r$, and thus $\dU(x_1,x_2) \ge c\lambda^{-1/20}r^\kappa$.

Finally, for $x_1\sim x_2$ (in $\bZ^2$), set $F(x_1,x_2) =\{ \dU(x_1,x_2) \ge c\lambda^{-1/20}r^\kappa\}$. The argument above gives that
$$  \sum_{x_1 \in B_{\infty}(0,2r)} \sum_{x_2 \sim x_1} \mathbf{1}_{ F(x_1,x_2) } \ge \mathbf{1}_{G(\lam,r)} \sum_{x_1 \in B_{\infty}(0,2r)} \sum_{x_2 \sim x_1}\mathbf{1}_{ F(x_1,x_2) } \ge c\lambda^{-1/20}r^\kappa \mathbf{1}_{G(\lam,r)}. $$
Hence taking expectations
\[\left(1-c_1 e^{-c_2 \lam^{1/40}}\right)c\lambda^{-1/20}r^\kappa \leq \sum_{x_1 \in B_{\infty}(0,2r)} \sum_{x_2 \sim x_1} \pP( F(x_1,x_2) ) \le c' r^2 \pP( F(0,e_1) ),\]
and the result follows by a simple reparameterisation.
\qed

A similar idea gives an upper bound. We begin by looking at the size of the finite component rooted at a vertex. In particular, for $x \in \bZ^2$, this is defined to be the set $A_x =\{ y :  x \in \gam(y,\infty) \}$. We also define the depth of $A_x$ by $ \mathrm{dep}(A_x) = \max\{ \dU(x,y) : y \in A_x \}$.

\begin{lemma}
For $\lam \in \bN$, we have that
\begin{align}
 \label{e:Atail-lb1}
 \pP( \mathrm{dep}(A_0)  = \lam) &\le  c_1  \lam^{-2/\kappa},  \\
\label{e:Atail-lb2}
 \pP( \mathrm{dep}(A_0)  \ge \lam ) &\le  c_2  \lambda^{-(2-\kappa) /\kappa}.
\end{align}
\end{lemma}

\proof Suppose that the event $G(\lam,r)$ defined in the proof of the lower bound holds, and again set $R:=re^{\lam^{1/40}/\kappa}$. NB. We suppose that $\lambda$ is large enough so that $(32\lambda)^{1/\kappa} r\leq R$. Let $x \in \BE(0,r)$, and suppose that $\mathrm{dep}(A_x) =s$, where $r^\kappa \le  s\le  (16\lam)^{-1} R^\kappa$. We then claim that $A_x \subseteq B_{\infty}(0,R/2)$. Indeed, suppose $y \in A_x \cap \BE(0,{R/2})^c$, and let $y'$ be the first point on $\gam(y,x)$ which is in $B_\infty(0,R)$. Since $r\leq \frac14 R \le d_\infty(x,y') \le \dUS(x,y')$, we must then have $s=\mathrm{dep}(A_x) \ge \dU(x,y') \ge \min\{\lam^{-1} R^\kappa,\lam^{-1} \dUS(x,y')^\kappa \}\ge \lam^{-1} 4^{-\kappa} R^\kappa$, which is a contradiction. Now, there exists $y_x \in A_x$ such that $\dU(x,y_x) =s$, and it must be the case that $y_x \in B_{\infty}(0,R/2)$. Thus the ball $B_\sU(y_x, s)$ is $\lam$-good, and we obtain
\be  \label{e:lbvolAx}
   | A_x| \ge    | B_\sU(y_x, s)| \ge \lam^{-1} \mathrm{dep}(A_x)^{2/\kappa}.
\ee
Next, let $s \in \bN \cap [r^\kappa, 2r^\kappa]$, and set $\tilde H_s = \{ x \in B_{\infty}(0,r): \mathrm{dep}(A_x) =s  \}$.
By \eqref{e:lbvolAx}, we then have that
\begin{equation}\label{sedf}
\mathbf{1}_{G(\lam,r) } \sum_{x \in \tilde H_s} |A_x| \ge \lam^{-1} s^{2/\kappa} |\tilde H_s|   \mathbf{1}_{G(\lam,r) }  .
\end{equation}
If $x \in \tilde H_s$ and $y \in A_x$, then $y \in B_{\infty}(0,R/2)$ and $\dU(x,y) \le s$. Hence $d_\infty(x,y)^\kappa \le \dUS(x,y)^\kappa \le \max\{\lam r^\kappa,\lam \dU(x,y)\} \le  2\lam r^\kappa$, and so $ y \in B_{\infty}(0,(1+(2\lam)^{1/\kappa})r)$. Since the sets $(A_x)_{x \in \tilde H_s}$ are disjoint, it follows that  $\mathbf{1}_{G(\lam,r) }  \sum_{x \in \tilde H_s} |A_x| \le   c \lam^{2/\kappa} r^2   \mathbf{1}_{G(\lam,r) }$, and combining this with the estimate at \eqref{sedf} yields
\be\label{lamest}
 |\tilde H_s|\mathbf{1}_{G(\lam,r) } \le c\lam^{1+2/\kappa}  \mathbf{1}_{G(\lam,r) }.
 \ee
Finally, let $\Lambda_*=\inf\{\lambda\geq \lambda_0:\:G(\lam,r)\mbox{ holds}\}$, and note that $\pP(\Lambda_*>\lam)\leq \pP(G(\lam,r)^c)\leq c_1e^{-c_2\lam^{1/40}}$. Thus $\Lambda_*$ is almost-surely finite,
and there exist finite constants $c_p$ such that $\bE (\Lam_*^p) \le c_p$; note that these constants
can be chosen not to depend on $r$.
Hence, from \eqref{lamest}, we deduce that, for $s\in \bN \cap [r^\kappa, 2r^\kappa]$,
\[\pP( \mathrm{dep}(A_0)=s)\leq r^{-2}|B_\infty(0,r)|\pP( \mathrm{dep}(A_0)=s)=r^{-2}\mathbf{E}(|\tilde H_s|)\leq cr^{-2}\mathbf{E}(\Lam_*^{1+2/\kappa})\leq cs^{-2/\kappa}.\]
This gives the bound \eqref{e:Atail-lb1}, and the bound \eqref{e:Atail-lb2} readily follows. \qed

\noindent
\emph{Proof of the upper bound of Proposition \ref{nnp}.}
Again, suppose that the event $G(\lam,r)$ defined in the proof of the lower bound holds, and $\lambda$ is chosen large enough so that $2\leq e^{\lam^{1/40}/\kappa}$. Suppose that  $\dU(0,e_1) >2^\kappa c_3 \lambda^{1/20} r^\kappa$ (where $c_3$ is as in the definition of the aforementioned event). Note that the dual vertices enclosed by the path $\gam(0,e_1)$ (combined with the edge from $0$ to $e_1$) are all elements of the finite component of $\sU'$ rooted at $0'$, which we denote $A'_{0'}$. On $G(\lam,r)$, we have $\dUS(0,e_1) \ge  2r$, and so there exists a point $x_1 \in \gam(0,e_1)$ such that $d_\infty(0,x_1) \ge 2r$. Let $x_2$ be a point on $\gam(0,e_1)$ adjacent to $x_1$, and let $\{x'_1, x'_2\}$ be the edge dual to  $\{x_1, x_2\}$. One of the points $x'_1$, $x_2'$ is in  $A'_{0'}$; we call this point $x'$. Since $d_{\sU'}^\sS(0',x') \ge r$ and $G(\lam,r)$ holds, we thus obtain that $\mathrm{dep}(A'_{0'} ) \ge c_3^{-1}\lambda^{-1/20}r^\kappa$. Hence $  G(\lam,r)\cap \{  \dU(0,e_1) \ge 2^\kappa c_3 \lambda^{1/20} r^\kappa \} \subseteq  G(\lam,r) \cap \{ \mathrm{dep}(A'_{0'} ) \ge c_3^{-1}\lambda^{-1/20}r^\kappa \}$, and so, setting $\tilde{r}:=c_3^{-1}\lambda^{-1/20}r^\kappa$ and $\lambda=c'(\log\tilde{r})^{40}$,
\begin{eqnarray*}
 \pP(  \dU(0,e_1) \ge\tilde{r} ) &\le&   \pP(G(\log\tilde{r},r)^c)  + \pP( \mathrm{dep}(A'_{0'} ) \ge 2^{-\kappa}c_3^{-2}(c')^{-1/10}(\log\tilde{r})^{-4}\tilde{r}^\kappa )\\
& \le& c_1 e^{-c_2c'\log\tilde{r}} +  c \tilde{r}^{-(2-\kappa)/\kappa}(\log\tilde{r})^{4(2-\kappa)/\kappa},
 \end{eqnarray*}
which, taking $c'$ suitably large, yields the desired result.
\qed

\noindent
\emph{Proof of Theorem \ref{T:dxy}(b).} Suppose that the event $G(\lam,r)$ defined in the proof of Proposition \ref{nnp} holds,  and write $R:=re^{\lam^{1/40}/\kappa}$. We will assume that we also have a parameter $t$ that satisfies $4^\kappa c_3\lambda^{1/20}\leq t\leq c_3^{-1}\lambda^{-1/20}e^{\lambda^{1/40}}$. Let $d_\infty(0,x)=r$, and $L$ be a shortest path in $\bZ^2$ between $0$ and $x$. If $\dU(0,x) \ge tr^\kappa$, then $d_{\sU}^{\sS}(0,x)\geq c_3^{-1/\kappa}\lambda^{-1/20\kappa}t^{1/\kappa}r$, and so there exists a point $y$ on $\gam(0,x)$ with $d_\infty(0,y) \ge  c_3^{-1/\kappa}\lambda^{-1/20\kappa}t^{1/\kappa}r$. Let $y'$ be a dual point with $d_\infty(y,y')=\frac12$ which is separated from infinity by $\gam(0,x) \cup L$. The path in the dual tree $\gam'(y',\infty)$ must pass through $L$ at a point $z'$, and the length of the section of $\gam'(y',\infty)$ inside $B_\infty(z',r)$ will be of length at least $c_3^{-1}\lambda^{-1/20}r^\kappa$. Moreover, for each vertex $z'$ of $\gam'(y',\infty)$ inside $B_\infty(z',r)$, it must be the case that $A'_{z'}$ has $d_\infty$-diameter greater than $ d_\infty(0,y')-2r$, and so $\mathrm{dep}(A'_{z'})\geq 2^{-\kappa}c_3^{-1}\lambda^{-1/20}tr^\kappa$. Now, let $H=\{ \dU(0,x) \ge tr^\kappa\}$ and  $F(z') =\{\mathrm{dep}(A'_{z'})\geq 2^{-\kappa}c_3^{-1}\lambda^{-1/20}tr^\kappa\}$. We then have that
$$\mathbf{1}_{H\cap G(\lam,r)} \sum_{z'\in B_{\infty}(0',2r)}\mathbf{1}_{ F(z') }
 \ge \mathbf{1}_{H\cap G(\lam,r)} c_3^{-1}\lambda^{-1/20}r^\kappa. $$
So, by \eqref{e:Atail-lb2},
\begin{align*}
 \pP( H \cap G(\lam,r))
 &\le c_3 r^{-\kappa}\lambda^{1/20} \bE \left( \mathbf{1}_{H\cap G(\lam,r)} \sum_{z'\in B_{\infty}(0',2r)}\mathbf{1}_{ F(z') } \right) \\
 &\le c_3 r^{-\kappa} \lambda^{1/20}  \sum_{z'\in B_{\infty}(0',2r)}  \pP(F(z') )  \\
 &\le c r^{2-\kappa}  \lambda^{1/20}   (2^{-\kappa}c_3^{-1}\lambda^{-1/20}tr^\kappa)^{-(2-\kappa)/\kappa}\\
 &= c  \lambda^{1/10\kappa}  t^{-(2-\kappa)/\kappa}.
\end{align*}
Hence, taking $\lambda=(\log{t})^{41}$, the result follows similarly to the end of the previous proof. \qed

\section{Controlling paths} \label{sec:control}

In this section, we provide a general technique for estimating from below the probability of seeing a particular path configuration in the UST. This will enable us to estimate from below the probability of seeing especially short paths between given points, and so prove the lower bound in Theorem \ref{T:dxy}(a). These estimates will also be a key ingredient in establishing volume and heat kernel fluctuations for the UST, as we do in the subsequent section.

Let $x\in \mathbb{Z}^2$, and $m\geq 1$. A \emph{scale $m$ path} from 0 to $x$, $\pi$ say, is a sequence of distinct vertices $0=x_0,x_1,\dots,x_N=x$ such that $x_i\in (m\mathbb{Z})^2$ and the Euclidean (i.e.\ $\ell^2$) distance between $x_{i-1}$ and $x_i$ is equal to $m$ for each $i=1,\dots,N-1$, and also $x_N\in B_m(x_{N-1})$, where we define
\[B_r(x):=B_\infty(x,r/2).\]
We write $|\pi| = N$ for the length of the path. Now, fix a path $\pi$ of length $N$. The rest of this section is devoted to defining an event ${F}_m(x,\pi)$ with $\pP( {F}_m(x,\pi) ) \ge e^{-c N}$ such that on this event the path from 0 to $x$ in the UST is contained in $\cup_{i=0}^NB_m(x_i)$ and (up to constants) has length $Nm^{\kappa}$.

To complete the program described in the previous paragraph, we will again appeal to Wilson's algorithm. As at the start of Section \ref{sec:LERW}, let $(S^x)_{x\in\mathbb{Z}^2}$ be a collection of independent simple random walks on $\mathbb{Z}^2$, where $S^x$ is started from $x$. Slightly modifying the algorithm at \eqref{wilsonalg}, we use these to construct the part of the UST containing both $0$ and $x$ via an iterative procedure. In particular, let $k$ be some integer that will be fixed later. We begin our construction by taking $\mathcal{U}_0=\gamma_0$ to be the loop-erasure of $S^{(0,\lfloor m/k\rfloor)}$ run until it first hits the origin $x_0=0$. We then continue as at \eqref{wilsonalg}: for $i\geq 1$, let $\mathcal{U}_i=\mathcal{U}_{i-1}\cup \gamma _i$, where $\gamma_i$ is the loop-erasure of $S^{x_i}$ run until it first hits $\mathcal{U}_{i-1}$. We will later use the notation $x_i'$ to represent the unique point in $\mathcal{U}_{i-1}\cap \gamma _i$. From Wilson's algorithm, we obtain that the path from $0$ to $x$ in the graph tree $\mathcal{U}_N$ is distributed identically to the path from 0 to $x$ in $\mathcal{U}$. For convenience, we will henceforth assume that $\mathcal{U}$ has been constructed by continuing with Wilson's algorithm from $\mathcal{U}_N$, and so this equality is almost-sure.

\begin{figure}[t]
\begin{center}
\input{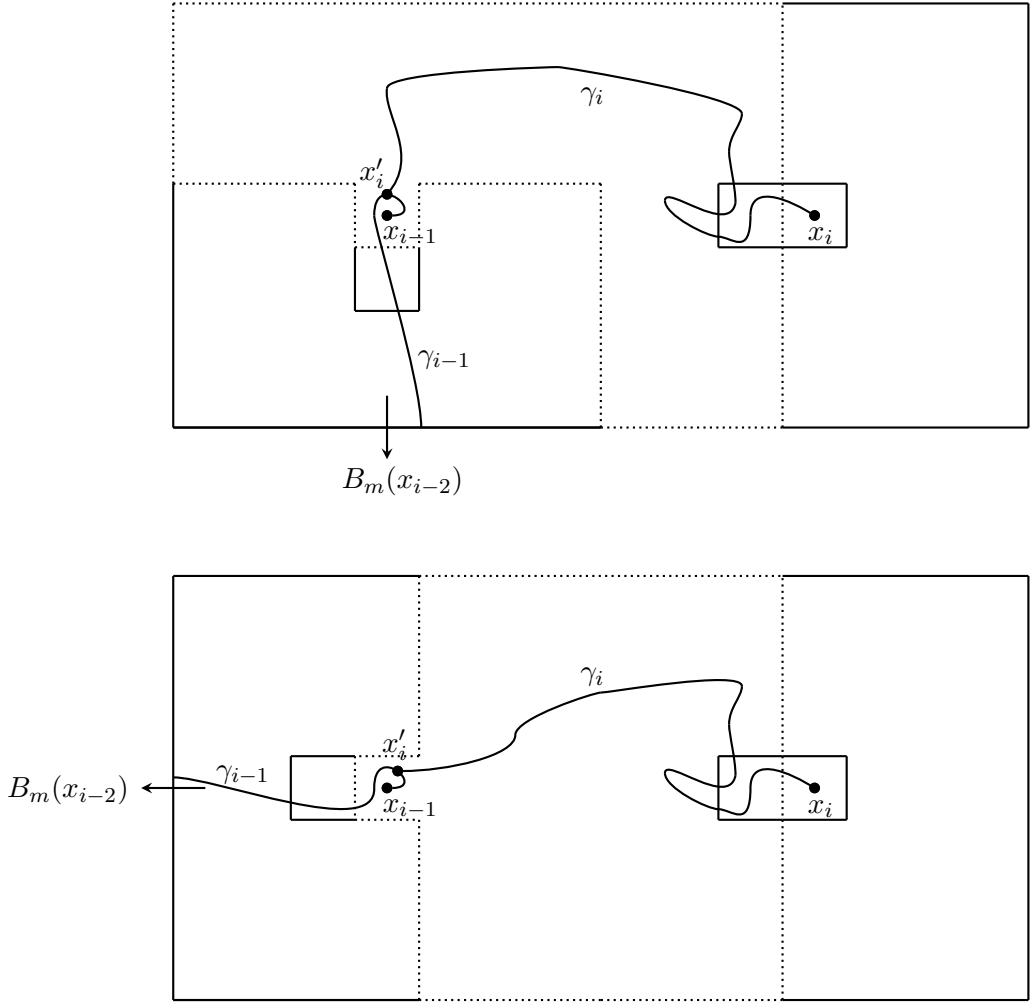}
\rput(-3,10.15){$x_i$}
\rput(-6,12){$\gamma_{i}$}
\rput(-7.95,8.5){$\gamma_{i-1}$}
\rput(-8.9,10.98){$x_i'$}
\rput(-8.4,10.15){$x_{i-1}$}
\rput(-8.5,6.9){$B_m(x_{i-2})$}
\rput(-3,2.55){$x_i$}
\rput(-6,4.3){$\gamma_{i}$}
\rput(-10.6,3.0){$\gamma_{i-1}$}
\rput(-8.6,3.4){$x_i'$}
\rput(-8.4,2.55){$x_{i-1}$}
\rput(-12.9,2.8){$B_m(x_{i-2})$}
\end{center}
\caption{Configuration within $B_m(x_{i-1})\cup B_m(x_i)$ on the event $G_0\cap\dots\cap G_i$.}\label{config}
\end{figure}

We next define a sequence of `good' events $G_i$. Given $\lambda\geq 2$, which will also be chosen later, set $G_0:=\{\gamma_0\subseteq B_m(0)\}\cap\{|\gamma_0|\leq \lambda m^{\kappa}\}$, where $|\gamma_0|$ is the number of elements of the path $\gamma_0$. To define $G_i$ for $i=1,\dots,N-1$, first let $R_i$ be the $m\lambda^{-2}\times 2m\lambda^{-2}$ rectangle consisting of $B_{m\lambda^{-2}}(x_i)$ and the adjacent square of side-length ${m\lambda^{-2}}$ that is closest to $x_{i-1}$; this is the rectangle about $x_i$ with a solid border shown in Figure \ref{config}. Moreover, let $Q_i$ be the union of the $m\times \frac{m(1-\lambda^{-2})}{2}$ rectangle contained in $B_m(x_i)$ that is closest to $x_{i-1}$, the $m\times \frac{m(1-\lambda^{-2})}{2}$ rectangle contained in $B_m(x_{i-1})$ that is furthest from $x_{i-2}$ (if $i=1$, take this to be the rectangle closest to $x_i$), and $B_{{m\lambda^{-2}}}(x_{i-1})$; this is the dotted region shown in Figure \ref{config}. Note in particular $Q_i$ has essentially two forms, depending on whether $x_{i-2},x_{i-1},x_i$ are co-linear or not; these are the two configurations are shown in Figure \ref{config}. For $i=1,\dots N-1$, we then set $G_i=\cap_{j=1}^3G_i^j$, where:
\begin{itemize}
  \item $G_i^1$ is the event that $S^{x_i}$ exits $R_i$ on the side closest to $x_{i-1}$ -- call this exit time $\tau_i$, and also $S^{x_i}_{\tau_i+\cdot}$ hits $\gamma_{i-1}$ before exiting $Q_i$;
  \item $G_i^2$ is the event that $|\gamma_i|\in [\lambda^{-1}m^\kappa,\lambda m^\kappa]$;
  \item $G_i^3$ is the event that $|\gamma_i\cap B_{3m\lambda^{-2}}(x_i)|\leq \lambda(3m\lambda^{-2})^{\kappa} = 3^{\kappa}\lambda^{-3/2}m^{\kappa}$.
\end{itemize}
Finally, we take
\begin{align}\nn
G_N &:=\left\{\gamma_N\subseteq B_m(x_{N-1})\cup B_m(x_{N})\right\}\cap\left\{|\gamma_N|
   \leq \lambda m^{\kappa}\right\}, \\
   \nn
    F_m(x,\pi) &:= \cap_{i=0}^N G_i.
\end{align}

To highlight the relevance of ${F}_m(x,\pi)=\cap_{i=0}^NG_i$ to controlling path lengths, we note that on this event we have that
\begin{equation}\label{distupper}
d_\mathcal{U}(0,x)\leq \sum_{i=1}^Nd_\mathcal{U}(x_{i-1},x_i)\leq \sum_{i=0}^N|\gamma_i|\leq \lambda (N+1)m^\kappa\leq 2\lambda Nm^\kappa.
\end{equation}
Moreover, from the construction it is possible to deduce that, on ${F}_m(x,\pi)$,
\begin{eqnarray}
d_\mathcal{U}(0,x)&\geq&\sum_{i=1}^{N-1}d_\mathcal{U}(x_{i}',x_{i+1}')\nonumber\\
&\geq & \sum_{i=1}^{N-1}\left(d_\mathcal{U}(x_{i}',x_i)-d_\mathcal{U}(x_i,x_{i+1}')\right)\nonumber\\
&\geq &\sum_{i=1}^{N-1}\left(|\gamma_i|-|\gamma_i\cap  B_{3m\lambda^{-2}}(x_i)|\right)\nonumber\\
&\geq & (N-1) \left(\lambda^{-1}-3^\kappa \lambda^{-3/2}\right)m^\kappa,
\label{distlower}
\end{eqnarray}
where the bound $d_\mathcal{U}(x_i,x_{i+1}')\leq |\gamma_i\cap  B_{3m\lambda^{-2}}(x_i)|$ is a consequence of the fact that the random walk $S^{x_i}$ does not backtrack to $B_{m\lambda^{-2}}(x_i)$ once it exits $R_i\subseteq B_{3m\lambda^{-2}}(x_i)$, which means that neither does $\gamma_i$, and so the path in $\mathcal{U}$ from $x_i$ to $x_{i+1}'$ must be contained within $B_{3m\lambda^{-2}}(x_i)$. In particular, for $N\geq 2$ and $\lambda$ suitably large, this implies $d_\mathcal{U}(0,x)\geq \frac{1}{4\lambda}Nm^\kappa$, and so we have the desired control over the lengths of paths.

The following is the key estimate of this section.

{\propn \label{pathprop} If $k$ and $\lambda$ are chosen large enough, then there exists a constant $c\in(0,\infty)$ and $m_0\in\mathbb{N}$ such that
\begin{equation}\label{pathpropest}
\mathbf{P}\left({F}_m(x,\pi)\right)\geq e^{-c |\pi| }
\end{equation}
whenever $m\geq m_0$, $x\in\mathbb{Z}^2$, and $\pi$ is a scale $m$ path $\pi$ from $0$ to $x$.}

\proof It is enough  to show that if $k$ and $\lambda$ are chosen large enough, then there exists a constant $c>0$ such that
\begin{equation}\label{aoest}
\mathbf{P}(G_0)\geq c
\end{equation}
and also, for $i=1,\dots, N$,
\begin{equation}\label{aiest}
\mathbf{P}\left(G_i\:\vline\:G_0,\dots,G_{i-1}\right)\geq c
\end{equation}
uniformly over $m$, $x$ and $\pi$.

To establish \eqref{aoest}, we first note that \cite[Lemma 2.6 and Proposition 2.7]{BM11} imply
\[\mathbf{P}(G_0)\geq \mathbf{P}\left(\gamma_0\subseteq B_m(0)\right)-\mathbf{P}\left(|\gamma_0|>\lambda m^\kappa\right)\geq 1-ck^{-1/3}-c(\lambda k^\kappa)^{-1/5},\]
whenever $m\geq k\geq 1$. Taking any $\lambda\geq 1$ and $k$ large yields a bound of the desired form.

For \eqref{aiest} when $i=1,\dots,N-1$, we start by bounding $\mathbf{P}\left(G_i\:\vline\:G_0,\dots,G_{i-1}\right)$ below by
\[\mathbf{P}\left(G_i^1\:\vline\:G_0,\dots,G_{i-1}\right)
-\mathbf{P}\left(G_i^1\cap (G_i^2)^c\:\vline\:G_0,\dots,G_{i-1}\right)
-\mathbf{P}\left(G_i^1\cap (G_i^3)^c\:\vline\:G_0,\dots,G_{i-1}\right).\]
Now, elementary random walk estimates yield that the first term here is bounded below by $c \lambda^{-4}$ whenever $\lambda^{2}\geq k$ (this latter inequality is required for the case $i=1$). Next, let $\tilde{\gamma}_i$ be the loop-erased random walk from $x_i$ to the boundary of $(B_m(x_i)\cup Q_i)\backslash\gamma_{i-1}$. On $G_i^1$, we have that $\gamma_i=\tilde{\gamma}_i$. Hence, by \cite[Theorem 5.8 and 6.7]{BM10},
\[\mathbf{P}\left(G_i^1\cap (G_i^2)^c\:\vline\:G_0,\dots,G_{i-1}\right)
\leq \mathbf{P}\left(|\tilde{\gamma}_{i}|\not\in [\lambda^{-1}m^k,\lambda m^k]\:\vline\:G_0,\dots,G_{i-1}\right)\leq c_1e^{-c_2\lambda^{3/5}}.\]
Similarly, let $\tilde{\tilde{\gamma}}_i$ be the loop-erased random walk from $x_i$ to the boundary of $(R_i\cup Q_i)\backslash\gamma_{i-1}$. On $G_i^1$, we have that $\gamma_i=\tilde{\tilde{\gamma}}_i$. So, by applying \cite[Theorem 5.8]{BM10} again,
\begin{eqnarray*}
\mathbf{P}\left(G_i^1\cap (G_i^3)^c\:\vline\:G_0,\dots,G_{i-1}\right)
&\leq& \mathbf{P}\left(|\tilde{\tilde{\gamma}}_{i}\cap B_{3m\lambda^{-2}}(x_i)|>\lambda (3m\lambda^{-2})^\kappa\:\vline\:G_0,\dots,G_{i-1}\right)\\
&\leq& c_1e^{-c_2\lambda}.
\end{eqnarray*}
Combining these estimates we obtain
\[\mathbf{P}\left(G_i\:\vline\:G_0,\dots,G_{i-1}\right)\geq c \lambda^{-4}-c_1e^{-c_2\lambda^{3/5}},\]
which is greater than $\frac{c}2\lambda^{-4}>0$ for large $\lambda$.

The estimate (\ref{aiest}) for $i=N$ is obtained similarly.
\qed

We can now conclude the proof of Theorem \ref{T:dxy} -- recall that (b) was proved at the end of Section \ref{sec:LERW}, and the upper bound in (a) follows from Theorem \ref{T:LERW-lb}.
\bigskip

\sm{\em Proof of the lower bound in Theorem \ref{T:dxy}(a)}.
Without loss of generality, we may assume $y=0$. Moreover, in view of Proposition \ref{pathprop}, we can choose constants $c_1,c_2,c_3,c_4,c_5,c_6$ such that if $k=c_1$, $\lambda=c_2$, then the estimate (\ref{pathpropest}) holds with $c=c_3$ and $m_0=c_4$, and also on ${F}_m(x,\pi)$ we have that
\[c_5Nm^\kappa \leq d_\mathcal{U}(0,x)\leq c_6Nm^\kappa.\]
Now, for any $m\geq1$ and $x\in\mathbb{Z}^2$ with $d_{\infty}(0,x)\geq m$, one can choose a scale $m$ path $\pi$ from $0$ to $x$ such that $N\leq c_7 d_{\infty}(0,x) /m$. On ${F}_m(x,\pi)$, we therefore have that
\[d_\mathcal{U}(0,x)\leq c_6c_7d_{\infty}(0,x)m^{\kappa-1}.\]
It readily follows that if $d_{\infty}(0,x)\geq c_4(c_6c_7 \lambda)^4$ and we choose $m={d_{\infty}(0,x)}/{(c_6c_7 \lambda)^4}$, which implies $m\geq c_4$ and $d_{\infty}(0,x)\geq m$, then on ${F}_m(x,\pi)$  it holds that $d_\mathcal{U}(0,x)\leq \lambda^{-1}d_{\infty}(0,x)^\kappa$. So we conclude that
\[\mathbf{P}\left(d_\mathcal{U}(0,x)\leq \lambda^{-1} d_{\infty}(0,x)^\kappa\right)\geq\mathbf{P}\left({F}_m(x,\pi)\right)\geq e^{-c_3 N}\geq e^{-c_3c_6^4c_7^5\lambda^4}.\]
\qed

{\remark \rm Whilst it would be straightforward to apply our approach to construct a corresponding exponential estimate from below for the probability of seeing exceptionally long paths in the UST, a stronger polynomial bound for such an event is already known. Indeed, by considering that with polynomially large probability the loop-erased random walk from $x$ to $y$ exits $B_\infty(x,\lambda d_{\infty}(x,y))$, it was established in \cite[Proposition 2.7]{BM11} that
\[\mathbf{P}\left(d_\mathcal{U}(x,y)\geq \lambda d_{\infty}(x,y)^\kappa\right)\geq c\lambda^{-4/5-\varepsilon}\]
for $x,y\in\mathbb{Z}^2$, $\lambda\geq \lambda_0$ (cf. Proposition \ref{nnp}). The point of our approach is that it also gives control of the macroscopic shape of the long path.}

\section{Volume fluctuations}  \label{sec:fluct}

In this section, we prove Theorem \ref{mainthm1}. The main ingredient in the proofs of these results is the following lemma, which provides tail bounds for the volume of balls in the UST.

{\lem
There exist constants $c_1, c_2$ such that, for all $r,\lambda\geq 1$,
\begin{equation}\label{bigvol}
\mathbf{P}\left(   |B_\mathcal{U}(0, r) |  \geq \lambda r^{2/\kappa} \right) \geq c_1 e^{-c_2\lambda^{\kappa/(\kappa-1)}}=c_1 e^{-c_2\lambda^{5}},
\end{equation}
and also
\begin{equation}\label{smallvol}
\mathbf{P}\left(   |B_\mathcal{U}(0, r) |  \leq \lambda^{-1} r^{2/\kappa} \right) \geq c_1 e^{-c_2\lambda^{\kappa/(2-\kappa)}}
=c_1 e^{-c_2\lambda^{5/3}}.
\end{equation}}

{\remark \rm See \cite[Theorem 1.2]{BM11} for upper bounds of $\exp( -c \lam^{1/3})$ for the probability in \eqref{bigvol} and of $\exp( -c \lam^{1/9})$ for the probability in \eqref{smallvol}.}

\begin{figure}[t]
\begin{center}
\scalebox{0.3}{\includegraphics{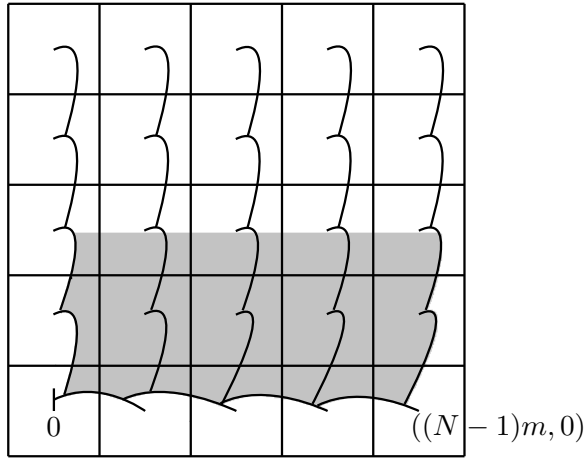}}
\rput(-5.55,.4){$0$}
\rput(0.3,.4){$((N-1)m,0)$}
\end{center}
\caption{The tree $\mathcal{U}_{N,m}$ on the event ${F}(N,m)$.}\label{grid}
\end{figure}

\proof Consider a square of $N\times N$ boxes, each of size $m\times m$, with the bottom left box centred on the origin. Let $\pi$ be the scale $m$ horizontal path from 0 to the point $x=((N-1)m,0)$, and suppose that the part of the UST containing $0$ and $x$ is constructed as in the event ${F}_m(x,\pi)$ of the previous section. Then, for each string of vertical boxes, assume that one has a similar construction, where at the bottom level we assume that the algorithm attaches to the horizontal part of the construction. If both such stages of this construction occur, we say that the event ${F}(N,m)$ holds.  (See Figure \ref{grid}.) Similarly to the proof of Proposition \ref{pathprop}, we have that
\begin{equation}\label{fbound}
\mathbf{P}\left({F}(N,m)\right)\geq e^{-cN^2}
\end{equation}
for all $N\geq 1$, $m\geq m_0$. Moreover, similarly to \eqref{distupper}, we deduce that, on ${F}(N,m)$,
\begin{equation}\label{b1}
d_\mathcal{U}(0,x)\leq c Nm^\kappa
\end{equation}
for every $x\in \mathcal{U}_{N,m}$, where $\mathcal{U}_{N,m}$ is the subset of $\mathcal{U}$ built in defining the event ${F}(N,m)$. Now, on ${F}(N,m)$, we have that every vertex in the $Nm\times Nm$ region of boxes is within a $d_\infty$-distance of $m$ from a vertex in $\mathcal{U}_{N,m}$. Thus, conditioning on $\mathcal{U}_{N,m}$ and continuing to construct the remainder of $\mathcal{U}$ from this tree as the root, by a minor adaptation of the `filling-in' argument of  \cite[Proposition 3.2]{BM11}, it is possible deduce that on an event of (conditional) probability greater than $1-c_1e^{-c_2 N^{{1/3}}}$ every vertex $x$ contained in the bottom $Nm \times \frac{Nm}{2}$ squares that is inside the outer paths (i.e.\ the shaded region of Figure \ref{grid}) satisfies
\begin{equation}\label{b2}
d_{\mathcal{U}}(x,\mathcal{U}_{N,m})\leq (N^{1/2}m)^\kappa\leq Nm^\kappa.
\end{equation}
In particular, putting the bounds at \eqref{b1} and \eqref{b2} together, we deduce that
\[\mathbf{P}\left(|B_\mathcal{U}(0,cNm^\kappa)|\geq c(Nm)^2\right)\geq
\left(1-c_1e^{-c_2 N^{{1/3}}}\right) e^{-cN^2}\geq c_3e^{-c_4N^2}.\]
Setting $r=cNm^\kappa$ and $\lambda = cN^{2(\kappa-1)/\kappa}/c^{2/\kappa}$ yields the result at \eqref{bigvol}.

\begin{figure}[t]
\begin{center}
\scalebox{0.3}{\includegraphics{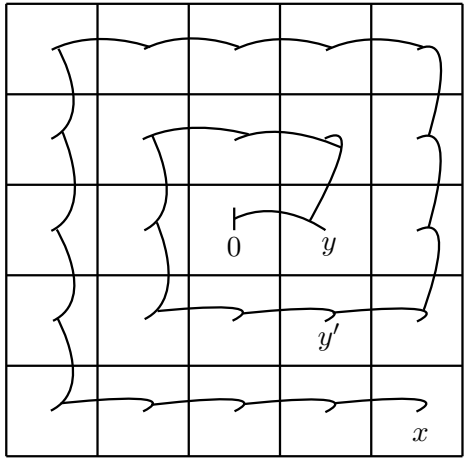}}
\rput(-3.15,2.8){$0$}
\rput(-1.9,2.8){$y$}
\rput(-1.9,1.6){$y'$}
\rput(-0.7,0.3){$x$}
\end{center}
\caption{The tree $\tilde{\mathcal{U}}_{N,m}$ on the event ${G}(N,m)$.}\label{spiral}
\end{figure}

For the result at (\ref{smallvol}), we argue similarly, though with a different initial tree configuration. We again consider a square of $N\times N$ boxes, each of size $m\times m$, centred on the origin. Let $\pi$ be the scale $m$ path that starts at 0 and spirals outwards around the boxes. Denoting the centre of the final box by $x$, we write ${G}(N,m)={F}_m(x,\pi)$. (See Figure \ref{spiral}.) From Proposition \ref{pathprop}, we have that
\[\mathbf{P}\left({G}(N,m)\right)\geq e^{-cN^2}.\]
Furthermore, let $y$ and $y'$ be centres of two adjacent boxes at a Euclidean distance approximately $Nm/3$ from the origin, but with $y$ one circuit closer to the origin than $y'$. (See Figure \ref{spiral}.) By arguing as at \eqref{distlower}, we have on ${G}(N,m)$ that $d_\mathcal{U}(0,y)\geq cN^2m^\kappa$.

Next, denote by $\tilde{\mathcal{U}}_{N,m}$ the tree constructed in the definition on ${G}(N,m)$, and note that every vertex in $B_{Nm}(0)$ is within $d_\infty$-distance $m$ of this set. Thus, similarly to the first part of the proof, we can again apply the `filling-in' argument of \cite[Proposition 3.2]{BM11} to deduce that on an event of (conditional on $\tilde{\mathcal{U}}_{N,m}$) probability greater than $1-c_1e^{-c_2 N^{{1/3}}}$ every vertex $x$ contained in $B_{Nm/2}(0)$ is within a $d_{\mathcal{U}}$-distance of $Nm^\kappa$. If this is the case and ${G}(N,m)$ occurs, it further holds that every point $z$ on the straight line between $y$ and $y'$ satisfies
\[d_\mathcal{U}(0,z)\geq d_\mathcal{U}(0,y)-d_\mathcal{U}(y,z)\geq cN^2m^\kappa -Nm^\kappa\geq cN^2m^\kappa,\]
where we have applied the lower bound on $d_\mathcal{U}(0,y)$ from the previous paragraph. In particular, since by construction any path in $\mathcal{U}$ from $B_{Nm}(0)^c$ to $0$ must pass through the line between $y$ and $y'$, it follows that $B_{Nm}(0)^c\subset B_{\mathcal{U}}(0, cN^2m^\kappa)^c$, which implies in turn that
\[\mathbf{P}\left(|B_\mathcal{U}(0,cN^2m^\kappa)|\leq (Nm)^2\right)\geq
\left(1-c_1e^{-c_2 N^{{1/3}}}\right) e^{-cN^2}\geq c_3e^{-c_4N^2}.\]
Setting $r=cN^2m^\kappa$ and $\lambda = cN^{\frac{2(2-\kappa)}{\kappa}} $ yields \eqref{smallvol}.
\qed

\noindent
{\em Proof of Theorem \ref{mainthm1}.} We start by showing large volumes occur almost-surely, i.e.\ \eqref{bigvolas1}. To this end, we define a sequence of scales:
\[D_i=e^{i^2},\qquad m_i =e^{i^2}/\varepsilon(\log i)^{1/2}.\]
We now run Wilson's algorithm, using the family of independent SRW $(S^x, x \in \bZ^2)$. At stage $i$ we use all the vertices in $B_{2 D_i}(0)$ which have not already been explored, in an order described in more detail below; write $\sU_i$ for the tree obtained. Let $\sF_i$ be the $\sigma$-field generated by the construction at the end of stage $i$.

By \cite[Theorem 1.1]{BM11}, we have that
\begin{equation}\label{inclusion}
B_{2D_i}(0)\subseteq B_\mathcal{U}(0,\lambda^4 D_i^\kappa)\subseteq B_{\lambda^5 D_i}(0)
\end{equation}
with probability greater than $1-c\lambda^{-17/16}$.

Hence, if we run Wilson's algorithm from the vertices contained inside $B_{2D_i}(0)$ (in any order), taking $0$ as the root, then the probability of seeing the part of the tree we generate, $\mathcal{U}_i$ say, leaving $B_{\lambda^5 D_i}(0)$ is less than  $c\lambda^{-17/16}$. By applying a Borel-Cantelli argument we thus obtain that
\begin{equation}\label{inclusion2}
\mathcal{U}_i\subseteq B_{i^5D_i}(0)\subseteq B_{m_{i+1}}(0)
\end{equation}
for large $i$, almost-surely. Moreover, from \eqref{inclusion}, we see that we may also assume that the $d_\mathcal{U}$-diameter of $\mathcal{U}_i$ is bounded above by $i^4D_i^\kappa\leq m_{i+1}^\kappa$ for large $i$, almost-surely. Define the event ${F}(i)$ to be the event that \eqref{inclusion2} and the diameter estimate for $\mathcal{U}_i$ both hold.

Next, we define an event ${G}(i+1)$ as follows. In particular, we first suppose that it incorporates the event ${F}(i)$ holding. We then mimic the definition of the event
\[{F}(D_{i+1}/m_{i+1}, m_{i+1})={F}(c(\log {(i+1)})^{1/2}, m_{i+1})\]
from the proof of (\ref{bigvol}). However, we run the first random walk in the box $B_{m_{i+1}}(0)$ until it hits the root $\mathcal{U}_i\subseteq  B_{m_{i+1}}(0)$, rather than the root $0$. Let $\mathcal{U}'(i+1)$ be the part of the UST that is thus constructed. Next, extend $\mathcal{U}'(i+1)$ to a tree $\mathcal{U}(i+1)$ by running loop-erased random walks from each of the vertices contained in the bottom $\frac{D_{i+1}}{m_{i+1}}\times \frac{D_{i+1}}{2m_{i+1}}$ squares that is inside the outer paths until they hit the part of the tree already constructed as in Wilson's algorithm (again, we refer to the shaded region of Figure \ref{grid}). We then complete the definition of ${G}(i+1)$ by supposing on this event that the $d_\mathcal{U}$-diameter of $\mathcal{U}(i+1)$ is bounded above by $c D_{i+1}m_{i+1}^{\kappa-1}$. (Since on ${F}(i)$ we also have an estimate for the $d_\mathcal{U}$-diameter of $\mathcal{U}(i)$ of $m_{i+1}^\kappa$, we can control the lengths of paths in the appropriate way.) Similarly to \eqref{fbound}, this construction yields that, for large $i$,
\[\mathbf{P}\left({G}(i+1)\:\vline\:\mathcal{F}_i\right)=\mathbf{P}\left({G}(i+1)\:\vline\:\mathcal{F}_i\right)\mathbf{1}_{{F}(i)}\geq e^{-c (D_{i+1}/m_{i+1})^2}\geq i^{-c}\]
for some $c<1$. Since it is clear that ${G}(i)$ is $\mathcal{F}_i$-measurable, then it follows from the conditional Borel-Cantelli lemma that ${G}(i)$ occurs infinitely often, almost-surely. Finally, we note that on ${G}(i)$ we have that $|B_\mathcal{U}(0,c D_{i+1}m_{i+1}^{\kappa-1})|\geq cD_{i+1}^2$. From this, the reparameterisation $r_i=cD_im_i^{\kappa-1}$ yields the result.

To prove (\ref{smallvolas1}), we proceed in essentially in the same way. In particular, define an event ${H}(i+1)$ similarly to ${G}(i+1)$, but based on the event ${G}(D_{i+1}/m_{i+1},m_{i+1})$ from the proof of (\ref{smallvol}) (i.e.\ using the spiral path of Figure \ref{spiral}, rather than the finger-like structure of Figure \ref{grid}), and then `filling-in' from all vertices in $B_{D_{i+1}/2}(0)$. Arguing as in the proof of (\ref{smallvol}), we deduce that $\mathbf{P}({H}(i+1)\:\vline\:\mathcal{F}_i)\geq i^{-c}$ for some $c<1$, ${H}(i)$ is $\mathcal{F}_i$-measurable, and moreover, on ${H}(i)$ we have that $|B_\mathcal{U}(0,cD_im_i^{\kappa-2})|\leq D_{i}^2$.
\qed

\section{Volume and resistance estimates on the UST}  \label{sec:VRest}

The aim of this section is to derive estimates for `good events', on which we have control on the volume and resistance of the two-dimensional UST. These will be applied in the subsequent sections to deduce the heat kernel estimates and other results stated in the introduction. Much of what we do here will build on previous work from \cite{BCK, BM11}. As already noted, the main input for the averaged heat kernel upper bound was the adaptation of \cite[Proposition 2.10]{BCK} that was established in Proposition \ref{P:PF1n}. A key difference in deriving the averaged heat kernel lower bound is that we will be need to understand the structure of the UST conditional on the presence of a given path, and deriving the relevant estimates requires substantial effort; our main result is Theorem \ref{thm:thm5-8}.

Our first two lemmas relate to the following situation. Suppose we have begun the construction of the UST using Wilson's algorithm, and have constructed a tree $\sU_0$. Write $\pP_\sT$ for the law of  $\sU$ conditional on the event $\{ \sU_0 = \sT\}$. We wish to adapt the unconditioned results of \cite{BCK, BM11} to the law $\pP_\sT$. We begin with the following `filling-in' lemma, based on  \cite[Lemma 2.3]{BCK} and  \cite[Proposition 3.2]{BM11}. If $\sT$ is a tree contained in $\sU$, then for each $x \in \bZ^2$ there exists a unique self-avoiding path in $\sU$ connecting $x$ and $\sT$; we denote this by $\gam(x,\sT)$. We write $\dU(x,\sT)  = \min\{ \dU(x,y), y \in \sT\}$ for the length of this path.

\begin{lem} \label{L:fillin} Let $r\ge 1$, and $\sT$ be a finite connected tree. There exist constants $c_i\in(0,\infty)$ not depending on $r$ and $\sT$ such that for each $\delta \le \frac14 $ the following holds. Let $A_1 \subseteq A_2$ be subsets of $\bZ^2$, with the property that any path in $\bZ^2 \setminus \sT$ between $A_1$ and $A_2^c$ is of length greater than $r$. Suppose that $\dE(x,\sT) \le \delta r$ for all $x \in A_2$.
Then there exists an event $G$ such that $$\pP_\sT(G^c)  \le c_1 r^{-2} |A_2|  \exp(-c_2 \delta^{-1/2}),$$ and on $G$ we have that, for all $x \in A_1$:
\[\dU(x, \sT) \le  ( \delta^{1/2} r)^\kappa;\qquad\dSU(x, \sT) \le  \delta^{1/2} r; \qquad \gam(x,\sT) \subseteq A_2. \]
\end{lem}

\proof If $A_1$ is a Euclidean ball of radius $r/2$, and $A_2$ is a Euclidean ball of radius $r$ centred on the same point, then this is immediate from the proof of  \cite[Lemma 2.3]{BCK}. (Checking the proof in \cite{BM11} one sees that one can take the power of $\delta$ to be $\delta^{-1/2}$ rather than $\delta^{-1/3}$.) The proof for more general sets $A_i$ is similar.
\qed

{\lem[{Cf.\ \cite[Lemma 2.5]{BCK}}]  \label{L:basicLP}
Let $x\in\mathbb{Z}^2$, $r\ge 1$, $k\geq 2$, and $D_0 \subseteq \bZ^2$ satisfy $B_{9r/8}(x) \subseteq \cup_{y \in D_0} B_{r/18k}(y)$. Let $\sT$ be a finite connected tree such that $\mathcal{T}\subseteq B_{2r}(x)^c$, and write $\gam = \gam(x,\sT)$. There exists an event $F_1=F_1(x,r,k)$ which satisfies
\[ \pP_\sT( F_1^c ) \le e^{-c_1 k^{1/8}},\]
and on $F_1(x,r,k)$ there exists $T \le \tau_{x,r}(\gam)$ such that, writing $W_{x}=\gam[T]$:\\
(a) $k^{-1/4} r^\kappa \le T \le  k^{1/4} r^\kappa$;\\
(b) $ a^{-2}r \le d_{\infty}(x,W_{x}) \le r$;\\
(c) there exists $Y_{x} \in D_0$ such that $d_{\infty}( Y_{x},  W_{x})\le {r}/{3k}$, $\dSU(Y_{x}, W_{x} ) \le {2r}/{3k}$ and also $\dU(Y_{x},W_{x})\le c_1 (r/k)^\kappa$.}

\sm \proof This follows as in \cite{BCK}. The most delicate part of the argument is to verify that \cite[Lemma 2.4]{BCK} holds in this context. For this, we need to show that if $\gam = \gam(x, \sT)$, then $\gam$ does not make too many close returns  to the segment $\gam[x, \tau_{x,r}]$ after time $\tau_{x,(1 + k^{-1/8})r}$. The argument in \cite[Lemma 2.4]{BCK} is for $\gam(x,\infty)$, and the proof for $\gam(x, \sT)$ is very similar.
\qed

\begin{figure}[t]
\begin{center}
\scalebox{0.4}{\includegraphics{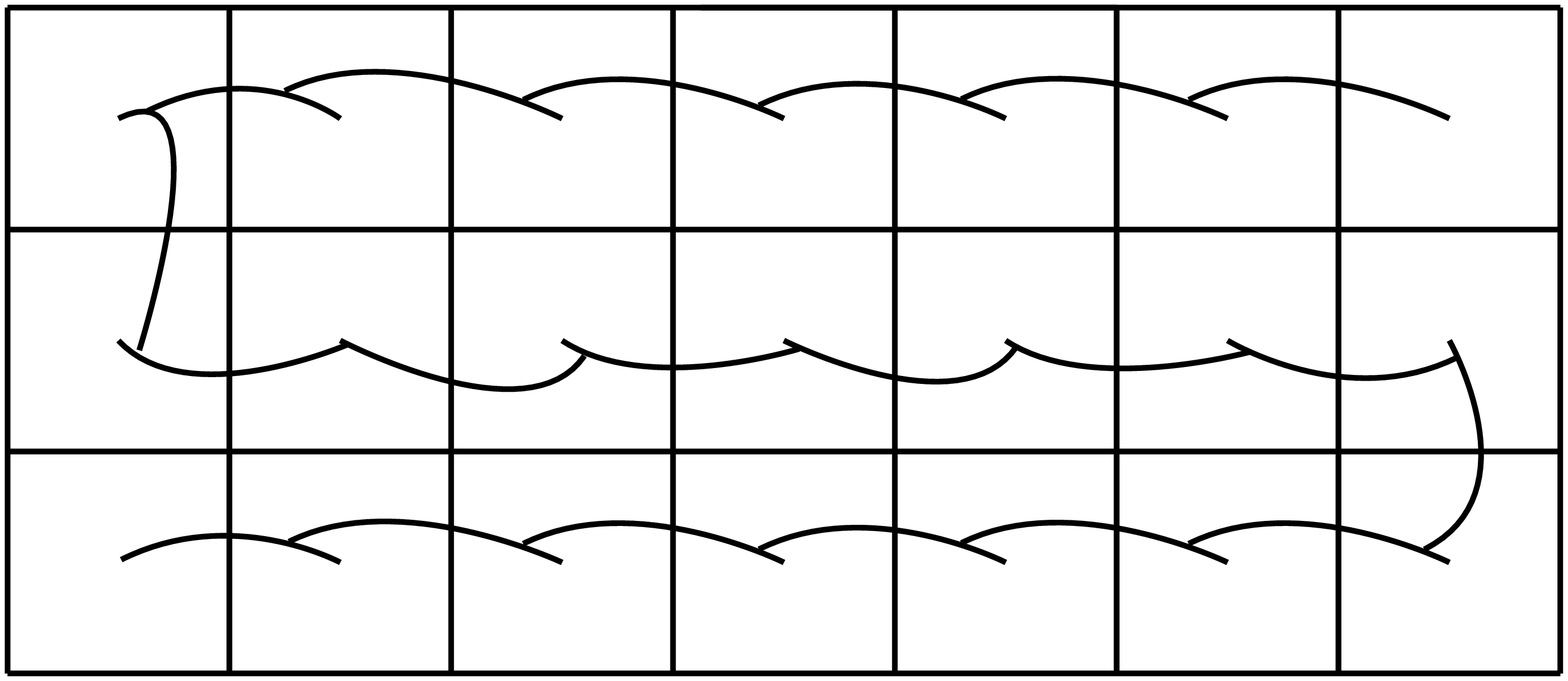}}
\end{center}
\caption{The part of $\sU$ constructed on the event $H_{N_1}$ (with
$N_1=3$).}\label{sec4fig}
\end{figure}

Towards stating Theorem \ref{thm:thm5-8}, we now introduce an event similar to those of the kind considered in Sections  \ref{sec:control} and \ref{sec:fluct}, but incorporating the regularity of Definition \ref{regdef}. We now choose $N \in \bN$ to be suitably large; in particular $N \ge 128$.  Let $(\tilde{x}_i)_{i=0}^{N_1}$ be the path with an `S shape' given by
\begin{align*}
&((-N,-1),(-N+1,-1),\dots,(N,-1),\\
&\qquad(N,0),(N-1,0),\dots,(-N,0),\\
&\qquad\qquad(-N,1),(-N+1,1),\dots,(N,1));
\end{align*}
note that $N_1=3(2N+1)-1$. Let $m\in\mathbb{N}$ with $m\geq 256$, we then let $({x}_i)_{i=0}^{N_1}$ be the corresponding scale $m$ path given by setting $x_i=m\tilde{x}_i$. (Ultimately we will only be interested in the situation when both $N$ and $m$ are very large.) Let $\Gam(0)$ be the path in $\bR^2$ which is the union of the line segments $[x_{i-1},x_i]$ for $i=1, \dots, N_1$. For $t >0$  we write
\[\Gam(t) = \{ x\in \bR^2: \dE(x, \Gam(0)) \le mt \}.\]
We now use Wilson's algorithm to construct $\sU$, and begin the construction using the points $x_i$, $0\le i \le N_1$. We wish the tree constructed to be inside $\Gam(1/8)$ -- see Figure \ref{sec4fig}, and also to have some additional regularity properties given below. As above, let $(S_n^x)_{n\geq 0}$, $x\in\mathbb{Z}^2$, be independent SRW in $\bZ^2$, where $S^x$ is started from $x$. For $i=1, \dots, N_1$ let $R'_i$ be the rectangle, with sides $m/4$ and $5m/4$ which contains both $\BE(x_{i-1}, m/8)$ and $\BE(x_i, m/8)$. Let $R_i \subset R'_i$ be the rectangle with sides $m/4$ and $m/2$ which contains $\BE(x_i, m/8)$, and let $Q_i = R'_i \setminus \BE(x_i, m/8)$. Let $\sU_0= \gam_0= \{0\}$, and for $i \ge 1$ let
$$ \gam_i = \sL( \sE_{\sU_{i-1}} (S^{x_i})), \qquad \sU_i = \sU_{i-1}  \cup \gam_i,\qquad \hbox{for } i=1, \dots, N_1. $$
Define events $G^j_i$, $i=1, \dots, N_1$, $j=1,2$, as follows.
\begin{itemize}
  \item $G_i^1$ is the event that $S^{x_i}$ first exits $R_i$ on the side closest to $x_{i-1}$, within a $d_\infty$-distance $m/16$ of the line segment between $x_{i-1}$ and $x_i$ -- call this exit time $\tau_i$, and also $S^{x_i}_{\tau_i+\cdot}$ hits $\gamma_{i-1}$ before exiting $Q_i$;
  \item $G_i^2$ is the event that $\gam_i$ is $(\lam, (3/2)m e^{-c_1 \lam^{1/2}}, 3m/2)$-regular.
\end{itemize}
Let $G_i = G^1_i \cap G^2_i$, and
\[H_i = \bigcap_{j=1}^i G_j, \qquad i=1, \dots, N_1.\]
Similarly to Proposition \ref{pathprop}, we then have the following.

\begin{proposition} \label{P:HN1}
There exist constants $c_1 \in(0,\infty)$ and $\lam_0,m_0\ge 2$ such that if $\lam\ge \lam_0$ and $m\geq m_0$, then
\be \label{e:PGNlb}
\pP( H_{N_1} ) \ge e^{- c_1 N_1}.
\ee
\end{proposition}

\proof Set $H_0 = \Omega$. Arguing as in the proof of  Lemma \ref{L:gam0y}, we have that $\pP(H_1)= \pP( G_1) \ge c_2>0$ for $m$ sufficiently large. From this, the result at \eqref{e:PGNlb} will follow if we can prove that for some $c_3$ that
\be \label{e:indh}
 \pP( H_i | H_{i-1} ) \ge e^{-c_3}, \qquad i=2, \dots, N_1.
\ee
We use induction. Suppose that \eqref{e:indh} holds for $i-1$ for some $i\geq 2$. Since $\sU_{i-1}$ contains a path from $x_{i-1}$ to the boundary of $\BE(x_{i-1},m/8)$, standard properties of the simple random walk $S^{x_i}$ give us that
\be\label{aaa}
\pP( G_{i}^1 |H_{i-1}) \ge 2e^{-c_4}.
\ee
Now, let $\wt \gam_i$ be the path in $\sU_i$ from $x_i$ to $\sU_{i-2}$. Note that, if we condition on $\sU_{i-2}$, then $\wt \gam_i$ is equal in law to $\sL( \sE_{\sU^c_{i-2}} (S^{x_i}) )$. Thus, if we set
$$ \wt G_{i}^2  = \{  \sE_{B_{3m}(x_i)} (\wt \gam_i)  \hbox{ is $(\lam, (3/2)me^{-c_1 \lam^{1/2}}, 3m/2)$-regular}\},$$
then Lemma \ref{L:regpath} yields that
\be\label{bbb}
\mathbf{P}\left( (\wt G_{i}^2)^c  | H_{i-2}\right) \le c_5e^{-c_6 \lambda^{1/2}}.
\ee
It is also straightforward to verify that $G_i^{1} \cap \wt G_{i}^2 \cap H_{i-1} \subseteq G_{i}^1 \cap G_{i}^2\cap H_{i-1}$, and it therefore follows that
\[ \pP( H_i | H_{i-1} )=\pP( G_i^1\cap{G}_i^2 | H_{i-1} )\geq \pP( G_i^1\cap \tilde{G}_i^2 | H_{i-1} )\geq \pP( G_{i}^1 |H_{i-1})-\mathbf{P}\left(\left(\wt G_{i}^2\right)^c  | H_{i-1}\right).\]
Using \eqref{bbb} and the inductive hypothesis we have that
\be  \label{e:pf3c}
 \mathbf{P}\left(\left(\wt G_{i}^2\right)^c  | H_{i-1}\right)=\frac{\mathbf{P}\left(\left(\wt G_{i}^2\right)^c\cap H_{i-1}  | H_{i-2}\right)}{\mathbf{P}\left(H_{i-1}  | H_{i-2}\right)}
\le e^{c_3}  \times c_5 e^{-c_6\lambda^{1/2}} \leq e^{-c_4},
\ee
where the final inequality holds by taking $\lambda_0$ suitably large. Combining \eqref{aaa} and \eqref{e:pf3c}, we thus obtain that \eqref{e:indh} holds for $i$.
\qed

We now fix $\lambda \ge \lam_0$, where $\lambda_0\ge 2$ is as in the previous proposition, and consider the uniform spanning tree obtained by conditioning on the event $H_{N_1}$. Let $\sT$ be a fixed tree such that $\pP( \sU_{N_1} = \sT \:|\: H_{N_1})>0$, and, as above, write $\pP_\sT( \cdot) = \pP( \cdot | \sU_{N_1} = \sT )$. We will derive volume and resistance bounds for balls $\BU(x,r)$ where $x$ is close to the middle section of $\sT$ and $r$ is of order $m^\kappa$ that will hold with $\pP_\sT$ probability close to 1. To this end, we introduce some more notation. The tree $\sT$ contains a path from $\tilde x_0$ to $\tilde x_{N_1}$; denote this by $\sT_{trunk}$. For $a>0$, let
\begin{align*}
    \sT_{0} &:= \sT \cap \left( [- (N-1) m,  (N-1) m]\times [-m/8,m/8]\right), \\
    \sT_{0,+} &:= \sT \cap \left( [-(N-1) m, (N-1) m]\times [7m/8,9m/8]\right).
\end{align*}
The following lemma, relating to distances on $\sT$, follows easily from the definition of the event $H_{N_1}$.

\begin{lemma} \label{L:dUonT}
Let $z_1, z_2 \in \sT_0$ with $d_\infty(z_1,z_2 ) \ge 3m$. Then
\[  c_1 \lam^{-1} m^{\kappa -1} d_\infty(z_1,z_2 ) \le \dU(z_1,z_2) \le   c_2 \lam m^{\kappa -1} d_\infty(z_1,z_2 ).\]
\end{lemma}

Now, for $a\in(0,1)$, let
\begin{equation}\label{ddef}
D(a):= [-aNm,aNm] \times [-9m/8,9m/8].
\end{equation}
We wish to define the region in $D(a)$ which lies `between' $\sT_0\cap \sT_{trunk}$ and $\sT_{0,+}\cap \sT_{trunk}$. To do this precisely, write $\sT_{trunk}^\bR$ for the continuous piecewise  linear self-avoiding  path  in $\bR^2$ obtained by connecting neighbouring points in $\sT_{trunk}$ by a line segment. Let $D^+_\bR(a)$ be the closure of the connected component of $D(a) \setminus  \sT_{trunk}^\bR$ which contains the point $(0, m/2)$, and define $D^+(a) = D^+_\bR(a) \cap \bZ^2$. To simplify our notation we will concentrate on the regions $D^+(a)$; exactly the same arguments apply to the corresponding region $D^-(a)$ lying `below' $\sT_0\cap \sT_{trunk}$.

Let $k \in \bN$ satisfy $k \ge \lam^4$ and $k \le m$. We now choose a grid $\Lam_1 \subset D^+(7/8)$ of points with separation of order $\frac14 k^{-1/4} m$ such that
\begin{align*}
   &D^+(3/4) \subseteq \bigcup_{ z \in \Lam_1 } \BE\left(z, \fract14 k^{-1/4} m\right).
\end{align*}
Since $| D(3/4)| \le c N m^2$, we can choose this set so that $ |\Lam_1| \le c N k^2$. Let $z \in  \Lam_1$ and set
\begin{align*}
 G_{11}(z) &:= \left\{ S^z \hbox{ hits $\sT$ before it leaves $\BE(z, k^{1/7} m )$} \right\}, \\
  G_{12}(z) &:=    \left\{ | \gam(z,\sT)| \le k^{1/7} m^\kappa \right\}, \\
   G_{13}(z) &:=    \left\{ | \gam(z,\sT)| \ge k^{-1/7} d_\infty(z, \sT)^\kappa \right\}.
\end{align*}

\begin{lemma} \label{L:connect2}
If $\lambda^4\leq k\leq m \wedge (N/8)^7$ and $z \in  \Lam_1$ then
\[  \pP_\sT (G_{1j}(z)^c )  \le c e^{-c k^{1/10} } \hbox{ for } j=1,2,3.\]
\end{lemma}

\proof Note that $S^z$ can only leave $\BE(z, k^{1/7} m )$ without hitting $\sT$ if it leaves horizontally at a distance of order $k^{1/7}m$ from $z$. Since every point in $D^+(3/4)$ is within a $\dE$ distance $5m/8$ of $\sT$ we obtain the bound on $\pP(G_{11}(z)^c)$. The bound for $G_{12}$ follows by part 1 of \cite[Theorem 2.2]{BM11} (with $D=D'=D(1)$ and $n=m$), and the bound for  $G_{13}$ follows by part 2 of the same theorem.
\qed

Now let
$$ F_2(k) = \bigcap_{ z \in \Lam_1} \left(G_{11}(z) \cap  G_{12}(z) \cap  G_{13}(z)\right); $$
by Lemma \ref{L:connect2} we have
\be \label{e:F2prob}
   \pP(  F_2(k)^c ) \le c N k^{2}   e^{-c k^{1/10} } \le c N e^{- c' k^{1/10} }.
\ee

\begin{proposition}\label{P:volub}
There exists $\delta_1>0$ such that the following holds.
Suppose that $\lam^7 \leq k \le m \wedge (\delta_1 N)^7$.
There exists an event $F_3=F_3(k)$ with
\be \label{e:PF5c}
 \pP_\sT( F_3^c) \le c N e^{-c_3 k^{1/10}}
\ee
such that on $F_3(k)$ the following properties hold.  \\
(a) If $y \in D^+(5/8)$, then there exists $x \in \sT$ with
$$ \dUS(y,x) \le 5 k^{1/7}m, \q \dU(y,x) \le  2 k^{1/7}m^\kappa. $$
(b) If $y \in  D^+(5/8)$ and $1 \le s \le \delta_1 \lam^{-1}  N $, then
\be \label{e:ballcontain}
   \BU(y, s m^\kappa)  \subseteq B_\infty(y, c (\lam s + k^{1/7}) m ) \cap D(3/4),
\ee
and thus
\be \label{e:volub}
 |  \BU(y, s m^\kappa)   | \le c'  (\lam s + k^{1/7})^2 m^2.
\ee
\end{proposition}

\proof
We continue the construction of the UST from $\sU_{N_1}$ by adding in the points in the grid $\Lam_1 \cap D^+(3/4)$; write $\sU_1^*$ for the tree thus obtained. We then complete the uniform spanning tree inside  $D^+(5/8)$. We use the filling-in of Lemma \ref{L:fillin} with $\delta = k^{-1/2}$, $r=m/8$, $A_1 = D^+(5/8)$, $A_2=D^+(3/4)$, and write $\tilde F_{3}(k)$ for the `good event' given by Lemma \ref{L:fillin}. Then
$$ \pP( \tilde F_{3}(k)^c ) \le c N e^{- k^{1/4} }. $$
Now let $ F_3(k) = F_{2}(k)  \cap \tilde F_{3}(k)$; the bound \eqref{e:PF5c} follows from \eqref{e:F2prob} and the bound on $\pP(\tilde F_{3}^c(k))$ given above.

In the remaining part of the proof, we assume $F_3(k)$ holds. Let $y_1 \in D^+(5/8)$. Then the event $\tilde F_{3}$ implies that there exists $w_1 \in \sU_1^*$ with $\dSU(y_1,w_1)\le k^{-1/4} m$ and $\dU(y_1,w_1) \le (  k^{-1/4} m )^\kappa$. By the construction of $\sU_1^*$ there exists a point $z_1 \in \Lam_1$ such that $w_1 \in \gam(z_1, \sT)$, and $\gam(z_1, \sT) \subset B_\infty(z_1, k^{1/7}m)$. Let $x_1$ be the point where $\gam(z_1, \sT)$ meets $\sT$. The events $G_{1i}(z_1)$ imply that $\dUS(w_1, x_1) \le 4  k^{1/7} m$, and $\dU(x_1, w_1) \le  k^{1/7} m^\kappa$, and the bounds in (a) follow immediately.

For part (b), let $y_1, y_2 \in D^+(5/8)$ with $\dU(y_1, y_2) \le  s m^\kappa$. Let $x_2$ be the point where $\gam(y_2,\sT)$ meets $\sT$ -- we may have $x_1=x_2$. As $\sU$ is a tree, we have $\dU(x_1,x_2) \le  \dU(y_1,y_2) \le s m^\kappa$. Using Lemma \ref{L:dUonT}, and taking $\lam_0$ large enough so that $c_1^{-1} \lam s \ge 3$, we obtain
$$ d_\infty(x_1 , x_2 ) \le c_1^{-1} \lam m s . $$
Since $\dE(y_j , x_j) \le 5 k^{1/7} m$, it follows that $\dE(y_1,y_2) \le  c \lam m s + 10 m  k^{1/7}$. This proves \eqref{e:ballcontain} and the volume upper bound \eqref{e:volub} is then immediate. \qed

We now consider resistance bounds.

\begin{proposition}
There exist $\delta_2>0$ and $c_1$ such that the following holds. Suppose that $\lam^7 \le k \le m \wedge (\delta_2 N)^7$, and $c_1 k^{1/7} \le s \le \delta_2 N$. Let $F_3(k)$ be as in the previous proposition. On the event $F_3(k)$, we have
\be \label{e:rebound}
 s m^\kappa \ge   \Reff(x,  \BU(x, s m^\kappa)^c ) \ge  \fract14 s m^\kappa
\hbox{ for } y \in D^+(9/16).
\ee
\end{proposition}

\proof
The upper bound is immediate. For the lower bound let $y \in D^+(5/8)$, and write $B_1=  \BU(y, \half s  m^\kappa)$, $B_2 = \BU(y, s m^\kappa)$. It is sufficient to prove that there are exactly two points in $\pd B_1$ which are connected to $\pd B_2$ by a path outside $B_1$; a cut set argument then gives the bound \eqref{e:rebound}.

Note first that by the construction of $\sT$ there exists $c$
such that each component $\sT'$ of $\sT \setminus \sT_{trunk}$  satisfies
$$ \dU(z,z') \le c \lam m^\kappa, \q  \dUS(z,z') \le c  m, \hbox{ for } z,z' \in \sT'. $$
By Proposition \ref{P:volub} and the  observation above there is a point $x \in \sT_{trunk}$ with $\dU(x,y) \le c_1 k^{1/7} m^\kappa$ and $\dUS(x,y) \le c_1  k^{1/7} m$; we used here the fact that $k \ge \lam^7$. Note that $x \in B_1$.

Let $w_1, w_2$ be the two points in  $\sT_{trunk}\cap \pd B_1$; it is clear that each of these is connected to $\pd B_2$ by a path outside $B_1$. Now let $z \in \pd B_2$ and suppose there exists $w \in \pd B_1$ with $w \neq w_1, w_2$ such that $\gam(w,z)$ is disjoint from $B_1$. By Proposition \ref{P:volub} we have $\dE(y,z) \le c (\lam s + k^{1/7}) m \le c' s m$, so choosing $\delta_2$ small enough we have that $z \in D^+(5/8)$. Let $z'$ be the closest point in $\sT$ to $z$. By  Proposition \ref{P:volub} we have $\dU(z,z') \le 2 k^{1/7} m^\kappa$, and it follows that the path $\gam(z,y)$ must intersect $\sT$. Let $z''$ be the closest point in $\sT_{trunk}$ to $z'$. The definition of $z$ implies that $z''$ must lie on $\sT_{trunk}$ between $w_1$ and $w_2$, and hence
$\dU(y,z'') \le \half s m^\kappa$. Thus
$$ \dU(y,z) \le \dU(y,z'') + \dU(z'',z') + \dU(z',z)  \le \half s m^\kappa +c \lam m^\kappa
+ 2 k^{1/7} m^\kappa, $$
which contradicts the fact that $z\in \partial B_2$ if $c_1$ is large enough.
 \qed

\begin{proposition}
There exist $\delta_3>0$ and $c_1$ such that the following holds. Suppose that $\lam \le k^{1/7} \le \delta_3 N$ and $k \le m^{1/2}$. There is an event $F_4(k)$ with $\pP_\sT( F_4^c) \le c N e^{-c' k^{1/24} }$ such that on $F_4(k)$, if $c_1 k^{1/4} \le s \le \delta_3 N$, then
\[  | \BU(x,s) | \ge c \lambda^{-1} k^{-5/2} sm^2  \q \hbox { for all  $x \in D^+(\half)$}.   \]

\end{proposition}
\proof
We follow the general lines of  \cite[Theorem 3.4]{BM11}, but note that the event $H_{N_1}$ means that the path $\sT$ cannot loop back on itself too much. This means that  the hardest part of the proof in \cite{BM11}, which uses \cite[Lemma 3.7]{BM11} is not needed.

We choose points $z_i \in \sT$, $1\le i \le N_2$, such that $\BE(z_i, m/2)$ are disjoint and $\sT \subseteq \cup_i \BE(z_i, m)$. We have $c N \le N_2 \le c' N$. Write  $m_1 = k^{-1/4} m$. For each $i$ choose points $w_{ij} \in \sT \cap B_\infty(z_i, m/2)$ with $1\le j\le N_3$ such that $\BE(w_{ij}, m_1)$ are disjoint and $N_3 \ge c k^{1/4}$.

The event $H_{N_1}$ implies that if $y_1, y_2 \in \mathcal{T}\cap\BE(z_i, m_1)$, then $\dU(y_1,y_2) \le c \lam m_1^\kappa$. Choose also $a, b >0$. Let $Q_1(w_{ij}) = \BE(w_{ij}, k^{-a} m_1)$ and $Q_2(w_{ij}) =  \BE(w_{ij}, 2 k^{-a} m_1)$. We cover $Q_2(w_{ij})$ by a grid $\Lam_3$ of points with separation $k^{-b} k^{-a} m_1$, so that $|\Lam_3| =4 k^{2b}$. We run Wilson's algorithm for the points in $\Lam_3$, and declare this stage of the construction a success for $w_{ij}$ if for all $y \in \Lam_3$, the random walk $S^y$ hits $\sT$ before it  leaves $\BE(w_{ij}, m_1)$. By the discrete Beurling estimate, \cite{LL}, the probability of failure $p_1$ satisfies
\[ p_1 \le |\Lam_3| c k^{-a/4} \le c k^{ 2b - a/4}. \]
We choose $a=1,b=1/12$ and $k$ large enough so that $p_1 < \half$.

Using the `stacks' construction in Wilson's algorithm, we can successively explore the UST in each box $\BE(w_{ij}, m_1)$, for $j=1, \dots, N_3$, and continue until we get a success. (See the argument in Theorem 5.4 of \cite{BJ} for more details.) Conditional on previous failures, the probability of a failure at stage $j$ is still bounded by $p_1$, and so since we have $N_3$ tries, the overall probability of failure is less than $c e^{- c k^{1/4}}$.

If this stage is a success for some $j$, we write $j_i =j$ and $w_i'= w_{i {j_i}}$. We then fill in the UST inside $Q_2(w'_i)$. By the filling-in result of Lemma \ref{L:fillin} (with $A_j = Q_j(w'_i)$, $j=1,2$ and $r = k^{-1} m_1$, $\delta = k^{-1/12}$) the `good' event $G=G(w'_i)$ given there satisfies $ \pP_\sT (G^c) \le c  e^{-c k^{1/24}}$. Moreover, if $G$ holds, then every path from a point in $\BE(w'_i, k^{-1} m_1)$ to $\sT$ is contained in $\BE(w'_i, m_1)$. Write $G'_i$ for the event that both stages of the construction are successful.

Now set
$$ \tilde F_{4} = \bigcap_{i=1}^{N_2} G'_i, $$
so that
$$ \pP_\sT ( \tilde F_{4}^c ) \le c N e^{ -c k^{1/24} } e^{- c k^{1/4}}\leq c N e^{ -c' k^{1/24} } , $$
and define $F_4 := \tilde F_{4} \cap F_3$,  where $F_3$ is the event defined in Proposition \ref{P:volub}.

Suppose now that $F_4$ holds. Let $y \in \BE(w'_i, k^{-1} m_1)$. The event $F_3$ implies that $\dU(y, \sT) \le c k^{1/4} m^\kappa$, and the event $G'_i$ implies that $\gam(y, \sT) \subseteq \BE(w'_i,m_1)$. It follows that $\dU(y, z_i) \le c_2 ( k^{1/4} + \lam)  m^\kappa$, and thus we obtain that
$$ | \BU(z_i, c_2 ( k^{1/4} + \lam)  m^\kappa )| \ge | \BE(w'_i, k^{-1} m_1) | \ge m^2 k^{-5/2}. $$
Next let $x \in D^+(1/2)$, and set $J = \{ z_i:   \dU(z_i, x) \le \half s m^\kappa \}$. By Lemma \ref{L:dUonT} we have $|J| \ge c \lam^{-1} s$. Taking $c_1$ large enough we have $\BU(z_i, c_2 ( k^{1/4} + \lam)  m^\kappa ) \subset \BU(x,sm^\kappa)$, and so
$$  \BU(x,sm^\kappa) \ge |J| m^2 k^{-5/2} = c \lam^{-1} s m^2 k^{-5/2}. $$
\qed

We summarize the estimates of this section in the following theorem.

\begin{thm}\label{thm:thm5-8}
There exist constants $c_i$ such that the following holds. Suppose that $\lam^7 \le k \le m^{1/2} \wedge c_1 N^4$. There exists an event $F_*=F_*(k)$ with $\pP_\sT( F_*^c) \le c_2 N e^{-c_3 k^{1/24}}$, such that on $F_*(k)$ if $x\in D^+(1/2)$, and $c_4 k^{1/4} \le s \le c_5 N$ then
\begin{equation}\label{vest}
c\lambda^{-1} k^{-5/2} m^2 s \le | \BU(x, s  m^\kappa)| \le c\lam m^2 s,
\end{equation}
\begin{equation}\label{rest}
\fract14 s m^\kappa  \le   \Reff(x,  \BU(x, s m^\kappa )^c ) \le s m^\kappa .
\end{equation}
\end{thm}

Finally, we give a local version of the previous result. For $x\in \mathcal{T}$, and $s \in [c_4 k^{1/4}, c_5 N]$ define $F_*(x,k,s)$ to be the event that the estimates at \eqref{vest} and \eqref{rest} hold. By only considering boxes of size $m$ within a distance $c sm$ of  $x$, rather than the order $N$ boxes considered in the previous result, one readily obtains the following.

\begin{cor}\label{corest}
Let $\lam^7 \le k \le m^{1/2} \wedge c_1 N^4$ and $c_4 k^{1/4} \le s \le c_5 N$.
Then if $x \in \sT_0 \cap D^+(1/2)$,
\[\pP_\sT( F_*(x,k,s)^c) \le c\lambda s e^{-c' k^{1/24}}.\]
\end{cor}

\section{Heat kernel bounds}  \label{sec:ann-b}

In this section, we will obtain heat kernel bounds using the estimates given in the previous section, starting with the quenched fluctuations.

\bigskip

\noindent
{\em Proof of Corollary \ref{cor1}.} By \cite[Theorem 4.1]{BK}, for any realisation of $\mathcal{U}$ we have that
\begin{equation}\label{eq:th41bk}
p^\mathcal{U}_{2r|B_{\mathcal{U}}(0,r)|}(0,0)\leq \frac{2}{|B_{\mathcal{U}}(0,r)|}.\end{equation}
Let $a_n=|B_{\mathcal{U}}(0,n)|n^{-d_f}$, and $t_n=n|B_{\mathcal{U}}(0,n)|$.
Plugging these into \eqref{eq:th41bk}, we have
\[
p^\mathcal{U}_{2t_n}(0,0)\leq \frac{2}{|B_{\mathcal{U}}(0,n)|}=2t_n^{-d_f/d_w}a_n^{-(1+2/\kappa)}
=2t_n^{-d_w/d_f}a_n^{-5/13}.\]
By \eqref{bigvolas1}, $a_n>(\log \log n)^{1/5}$ infinitely often, almost-surely, giving the
liminf statement.

We next prove the limsup statement. We use the construction of $\sU$ given after the proof of Lemma \ref{L:basicLP}  with a $3(2N+1)$ array of boxes of side $m$. As in Proposition \ref{P:HN1}, this has probability of success of at least $e^{-c_1 N}$.
Given $m$, we define
$R_*=m^\kappa$,
$N = \frac 1{2c_1} \log \log R_*$, $R=\half N m^\kappa$, and define
\begin{align*}
  v(t) &=   \begin{cases}t^\df,  &~~\hbox{ if }~ t \le  R_*,\\
   t R_*^{\df-1},    &~~\hbox{ if }~ t \ge  R_*.\end{cases}
 \end{align*}
Let $\sT$ be the tree given right after Proposition \ref{P:HN1}. Then by Theorem \ref{thm:thm5-8}, there is an event $F_*(k)$ with  $\pP_\sT( F_*(k)^c) \le c N e^{-c' k^{1/24}}$ such that the following hold on $F_*(k)$:
 \begin{align*}
  c_{1,k}v(R_0)\le |   \BU(0, R_0)  | \le c_{2,k}v(R_0),\\
  c_{3,k}R_0\le   \Reff(0,  \BU(0, R_0)^c ) \le R_0,~~~~~~
 \end{align*}
 where $R_0:=ck^\theta m^\kappa$ and $c_{i,k}:=k^{q_i}$, where $q_i\in\mathbb{R}$. We now follow \cite{KM}. Let $r(t)=t$ and let $\sI(t)$ be the inverse function of $r(t)\cdot v(t)$. After some calculations we get, noting $\dw=1+\df$,
\begin{align*}
  \sI(t) =
  \begin{cases}
   t^{1/\dw},& ~\hbox{ if }~~ t \le  R_*^{d_w}, \\
   t^{1/2} R_*^{-(\df-1)/2},& ~\hbox{ if }~~ t \ge   R_*^{\dw},
  \end{cases}
 \end{align*}
and the function $\tilde k(t) = v(\sI(t))$ is
\begin{align*}
  \tilde k(t) = \begin{cases}
  t^{\df/(1+\df)}=t^{\df/\dw},  &\hbox{ if }~~ t \le  R_*^{\dw}, \\
  t^{1/2} R_*^{(\df-1)/2} ,   &\hbox{ if }~~ t \ge   R_*^{\dw}.
 \end{cases}
  \end{align*}
We can rewrite the final line as
$$ \tilde k(t) = t^{\df/\dw} \left( \frac{R_*^{\dw}}{t} \right)^\alpha, ~~~\mbox{where}~~~
\alpha = \frac{ \df-1}{2\dw} >0.$$
One then finds from \cite[Proposition 3.3]{KM} that
$$ p^\sU_{2n}(0,0) \ge \frac{ c_{1,k}} { \tilde k(n) } \, \hbox { for }  \, \frac {c_{2,k}}2 R^{\dw} \le n \le c_{2,k}R^{\dw},$$
for some $c_{i,k}=c_ik^{-q_i}$, $i=1,2$, where $c_i,q_i>0$. (Note that \cite[Proposition 3.3]{KM} holds for $R\ge R_*$ if the assumption of the proposition holds for $R\ge R_*$.) So taking $T = c_{2,\lambda} R^{\dw}/2= \frac{c_{2,\lambda}}{2(4c_1)^{\dw}}(R_* \log \log R_*)^{\dw}$, it holds that, given $\sT$, with probability greater than $1-c N e^{-c' k^{1/24}}$, we have
$$  T^{\df/\dw} \tilde{p}^\sU_T(0,0) \ge c_{3,k}(\log \log R_*)^{\alpha \dw} \ge c_{3,k}' (\log \log T)^{\alpha \dw}. $$
In order to have $1-c N e^{-c' k^{1/24}}\ge 1/2$, it is enough to take $k\asymp (\log N)^{24}$ which is comparable to $(\log\log\log R_*)^{24}$. (Note that this choice of $k$ enjoys $k\leq \sqrt{m}\wedge (N/8c\lambda)^{1/\theta}$ that is required in Theorem \ref{thm:thm5-8}.) Hence we have
$$  T^{\df/\dw} \tilde{p}^\sU_T(0,0) \ge c_6 (\log \log T)^{\alpha \dw-\varepsilon}, $$
for some $c_6>0$ and $\varepsilon>0$ which is small.

The rest of the argument goes through similarly to the proof of Theorem \ref{mainthm1} in Section \ref{sec:fluct}.
We choose $m(i) = e^{i^2/\kappa}$, so $R_*(i) = e^{i^2}$,
 $N(i) = (2c_1)^{-1} \log\log R_*(i)$, and
 $\sum_i e^{-c_1 N(i)} = \sum_i i^{-1}=\infty$.
 Similarly to \eqref{inclusion2}, we have good separation of scales. Using the conditional Borel-Cantelli lemma, we obtain the desired lower bound.
\qed

\subsection{Averaged heat kernel upper bound}

To establish the upper bound of Theorem \ref{mainthm3}, we start by deducing upper estimates for the transition density that hold on the event  $F_1(\lam,n)$, which was defined at \eqref{f1def}. In this subsection, we fix $\eps_0= 1/40$. Moreover, throughout this and the next subsection, we will write
\[\Phi(t,r)    = \Big( \frac{r^\dw}{t} \Big)^{1/(\dw-1)},\]
where $d_w=13/5$ was introduced at \eqref{dwdef}. We also set
\be\label{sigxr}
\sigma_{x,r} = \inf\left\{n \ge 0: \dU(x,X^\sU_n)=r\right\}, \qquad T_x= \inf\left\{ n\ge 0: X^\sU_n =x \right\}.
\ee

\begin{lemma}[{\cite[Proposition 3.3]{KM}}]\label{L:KM34} There exist constants $c_i$ and $q_i$ such that if  $B_\sU(x,r)$ and $B_\sU(x,c_1\lambda^{-q_1}r)$ are $\lam$-good, then
\begin{align*}
  &p^\sU_t(x,x) \le c_2 \lam^{q_2} t^{-{d_f}/{d_w}},\qquad\hbox{if }\half r^{d_w} \le t\le r^{d_w}~\hbox{and }  t\in \bN,\\
 &P_x^\sU\left( \sigma_{x,r} > c_3 \lam^{-q_3} r^{d_w}\right)  \ge c_4 \lam^{-q_4}.
\end{align*}
\end{lemma}

\begin{lemma} \label{L:F1easyub} There exists $\lam_0$ such that if $\lam\ge \lam_0$ then
on the event $F_1(\lam,n)$ it holds that
\begin{align}
  &p^\sU_{t}(x,x) \le   c_2\lam^{q_2} t^{-{d_f}/{d_w}}, \q \hbox{ for all  $x \in B_\infty(0,n)$, $n^{\dw \kappa} e^{- \lam^{\eps_0}/2} \le t \le n^{\dw \kappa}$
  and $t\in \bN$}, \nonumber\\
  \label{e:smallexit}
  &P^\sU_x( \sigma_{x,r} > c_3 \lam^{-q_3} r^{d_w})  \ge c_4 \lam^{-q_4},    \q \hbox{ for all  $x \in B_\infty(0,n)$
  and    $n^{\kappa} e^{- \lam^{\eps_0}/2 } \le r \le n^{\kappa}$.    }
  \end{align}
\end{lemma}
\proof This follows immediately from Lemma \ref{L:KM34} and the definition of $F_1(\lam,n)$.
\qed

\begin{lemma} \label{L:tailhit}
There exist $\lam_0$, $q>0$ such that the follows holds. Suppose $F_1(\lam,n)$ holds with $\lam \ge \lam_0$. If $n/8 \le d_\infty(0,x) \le 7n/8$, then
\be \label{e:med-hit}
 P^\sU_0( T_x \le t ) \le 2\exp( - c \lam^{-q} \Phi(t,n^{\kappa}) )
\q \hbox { for } \;   t \ge n^{\kappa \dw} e^{-\lam^{\eps_0}/2}.
\ee
\end{lemma}

\proof Let $y$ be the first point on the path $\gam(0,x)$ with $d_\infty(0,y)\ge n/9$. By Lemma \ref{L:F1dist} we have $\dU(0,y) \ge c \lam^{-\kappa} n^{\kappa}$, and it is clear that $T_y \le T_x$. Let $m \ge 1$, and set
$$ R = c \lam^{-\kappa} n^{\kappa}, \q r = \frac{R}{m}, \q s = \frac{t}{m}. $$
Moreover, let $x_i$ be points on $\gam(0,x)$ such that $x_0=0$ and $d_\sU(x_{i-1},x_i)=r$, and $\xi_i$ be the duration from $T_{x_i}$ until $X^\sU$ leaves $B_\sU(x_i, r)$. Using this notation, we have that
$$ T_y \ge \sum_{i=0}^{m-1} \xi_i \ge s  \sum_{i=0}^{m-1} \mathbf{1}_{\{\xi_i>s\}}. $$
We next choose $m \in \bN$ so that $m \in [m_1,m_1+1]$ where
$$ m_1^{\dw-1} = \lam^{-b} \frac{ n^{\kappa \dw}}{t}, $$
and $b >0$ will be chosen later. We set $q = (\dw q_4 + b)/(\dw-1)$. If $t \ge \lam^{-b} n^{\kappa \dw}$, then $\Phi(t,n^\kappa)\leq \lambda^{b/(d_w-1)}$, and so our choice of $q$ ensures that  the probability bound in \eqref{e:med-hit} holds.  Thus we will assume that $m_1\ge 1$, so that $m_1 \le m \le 2m_1$ and $r = {R}/{m} \ge {R}/{2m_1}$. Now, the condition on $r$ in \eqref{e:smallexit} holds if ${R}/{2m_1} = {c \lam^{-\kappa} n^\kappa }/{2m_1} \ge e^{-\lam^{\eps_0}/2} n^\kappa$. This is equivalent to
$$ 2^{\dw-1} \lam^{-b } \frac{ n^{\kappa \dw}}{t} \le (c\lam^{-\kappa}  e^{\lam^{\eps_0}/2})^{d_w-1}, $$
i.e.\ $ t \ge c \lam^{\kappa(d_w-1)-b} e^{-\lam^{\eps_0}(d_w-1)/2 } n^{\kappa \dw}$, and we observe that this holds if $t \ge n^{\kappa \dw} e^{-3\lam^{\eps_0}/5}$ and $\lambda$ is large enough.
To apply \eqref{e:smallexit}, we will also need that $s = t/m$ satisfies $s \le  c_3\lam^{-q_3} r^\dw$. After some algebra we find this requires $\lam^{-b} \le  c_3 \lam^{-q_3-\kappa d_w}$.  So, taking $b > q_3+ \kappa d_w$ and $\lam_0$ large enough, this condition is also satisfied.

With the choice of $m$ in the previous paragraph, we can apply the bounds in \eqref{e:smallexit} to deduce that $\sum_{i=0}^{m-1} \mathbf{1}_{\{\xi_i>s\}}$ stochastically dominates a binomial random variable with parameters $m$ and $p= c_4 \lam^{-q_4}$. Applying the following general bound for a binomial random variable $\eta$,
$$ \mathbf{P}\left( |\eta- \mathbf{E}\eta | > t \right) \le 2 \exp( -t^2/(2 \mathbf{E}\eta + 2t/3) ),$$
we thus deduce that
\begin{align*}
P^\sU_0(  T_y <  \half c_4\lambda^{-q_4}t )=P^\sU_0(  T_y <  \half s mp  ) \le 2 e^{ - c m p} =  2 \exp( -c \lam^{-q_4} m ),
\end{align*}
and the result follows by a reparameterisation of $t$.
\qed

The following result improves upon the corresponding bound in \cite[Proposition 4.15]{BM11} by obtaining an upper bound for  $p^\sU_T(x,y)$ on a set for which the probability has a uniform (in $n$) lower bound. Whilst the estimate holds for a more limited range of times, it is enough for our purposes. We take $\eps_1 < \eps_0=1/40$.

\begin{proposition} \label{P:G-ub}
Suppose $F_1(\lam,n)$ holds with $\lam \ge \lam_0$, then
\[ p^\sU_t(0,x) \le \lam^{q_5} t^{-\df/\dw} \exp\left( -\lam^{-q_5} \Phi(t,n^{\kappa})  \right)\]
whenever $x \in B_{3n/4}(0)\setminus B_{n/2}(0)$, $ e^{-\lam^{\eps_1}} n^{\kappa \dw}\leq t\leq n^{\kappa d_w}$ and $t\in \bN$.
\end{proposition}

\proof  Let $z_1$ be the first point on the path $\gam(0,x)$ with $d_\infty(0,z_1)\geq n/8$, and $z_2$ be the first point on the path $\gam(x,0)$ with $d_\infty(x,z_2)\geq n/8$. Let $A_0$ be the set of points $y$ in $\bZ^2$ such that the path $\gam(0,y)$ does not contain $z_1$, and $A_x= \bZ^2 \setminus A_0$. Then, as in the proof of \cite[Theorem 4.9]{BK},
\begin{align}
 P^\sU_0( X^\sU_t =x) &= P^\sU_0( X^\sU_t =x, X^\sU_{[t/2]} \in A_0) + P^\sU_0( X^\sU_t =x, X^\sU_{[t/2]} \in A_x) \nonumber\\
&\le 4 P^\sU_x( X^\sU_t =0, X^\sU_{[t/2]} \in A_0)+ P^\sU_0( X^\sU_t =x, X^\sU_{[t/2]} \in A_x) \nonumber\\
&\le 4 P^\sU_x( X^\sU_t =0, T_{z_2} < t/2 ) + P^\sU_0( X^\sU_t =x,  T_{z_1} < t/2).\label{eeee}
\end{align}
For the second term above,
\begin{align*}
  P^\sU_0( X^\sU_t =x,  T_{z_1} < t/2)
  &=  E^\sU_0 \left( \mathbf{1}_{\{T_{z_1} < t/2\}} P^\sU_{z_1}( X^\sU_{t-T_{z_1}}=x)  \right)\\
&\le P^\sU_0(T_{z_1} < t/2 ) \sup_{t/2\le s \le t, s\in \bN} P^\sU_{z_1}( X^\sU_s=x) \\
&\le 4P^\sU_0(T_{z_1} < t/2 ) \sup_{t/2\le s \le t, s\in \bN} \sqrt{p^\sU_s(z_1,z_1)p^\sU_s(x,x)} \\
&\le c\lam^{q }   \exp\left( -c'\lam^{-q} \Phi(t,n^{\kappa})  \right) \lam^q  t^{-\df/\dw}.
\end{align*}
Here we used the Cauchy-Schwarz for the penultimate bound, and Lemmas \ref{L:F1easyub} and \ref{L:tailhit} to obtain the final one. The first term of \eqref{eeee} is bounded in the same way. \qed

We now have all the pieces in place, and the one remaining lemma we give provides the means to put these together.

\begin{lemma} \label{L:Epsum} Let $G_k$, $k \ge 1$, be a
sequence of sets with $\pP(G_k) \rightarrow 1$
and let $T\in \bN$.
If we have $p^\sU_T(0,x) \le a_k$ on $G_k$ for each $k$, then
\be \label{e:aksum}
 \pE p^\sU_T(0,x) \le  a_1  + \sum_{k=2}^\infty  a_{k} \pP( G_{k-1}^c).
 \ee
 \end{lemma}
\proof Set $A_1=G_1$ and $A_k = G_k\backslash G_{k-1}$ for $k\ge 2$. Since
$\pP(\cup_k G_k)=1$, we have $\pP(\cup_k A_k)=1$, and thus
$$ \pE p^\sU_T(0,x) = \sum_{k=1}^\infty  \pE( p^\sU_T(0,x)\mathbf{1}_{A_k})
\le \sum_{k=1}^\infty a_k  \pP(A_k)
\le a_1 \pP(A_1) + \sum_{k=2}^\infty  a_k \pP( G_{k-1}^c). $$
\qed

\noindent
\emph{Proof of the upper bound of Theorem \ref{mainthm3}.} By \cite[Theorem 4.4]{BM11}, we have that $ \pE p^\sU_{2T} (0,0) \le c T^{-\df/\dw}$.  Hence, applying the Cauchy-Schwarz as in the proof of Proposition \ref{P:G-ub}, we further have that, for all $x \in \bZ^2$,
 \be\label{e:aub-near}
 \pE p^\sU_{2T} (0,x)\le \pE[p^\sU_{2n}(0,0)^{1/2}  p^\sU_{2T}(x,x)^{1/2}]\le  \pE(p^\sU_{2T}(0,0))^{1/2} \pE(p^\sU_{2T}(x,x))^{1/2} \le c T^{-\df/\dw}.
 \ee
Hence if $d_\infty(0,x)\leq 16$, then the result follows.

Now let $d_\infty(0,x) \ge 16$.
Choose $n$ such that $x \in B_{3n/4}(0)\backslash B_{n/2}(0)$, and set $\Phi=\Phi(T,n^\kappa)$.
Set for $k \ge 1$,
$$ \lam_k =  k^{1/\eps_0} \Phi^{1/(\eps_0+q_5)}.$$
Choose $c_1$ so that $\Phi \ge c_1$ implies that
$\lambda_k^{-(d_w-1)(\varepsilon_0+q_5)}\ge e^{-\lambda_k^{\varepsilon_1}}$.
If $\Phi \le c_1$, then the estimate again follows from \eqref{e:aub-near}, so we
assume that $\Phi > c_1$.
We now use Lemma \ref{L:Epsum} with $G_k = F_1(\lam_k,n)$.
The definition of $\lam_k$ gives that
$T \geq n^{\kappa d_w}\lambda_k^{-(d_w-1)(\varepsilon_0+q_5)}\geq n^{\kappa d_w}e^{-\lambda_k^{\varepsilon_1}}$,
so Proposition \ref{P:G-ub} allows us to take
$$ a_k = T^{-\df/\dw} \lam_k^{q_5}   \exp\big( -\lam_k^{-q_5} \Phi \big). $$
Thus the first term in the sum \eqref{e:aksum} is given by
$$ a_1= T^{-\df/\dw} \Phi^{q_5/(q_5+\eps_0)} \exp( - \Phi^{\eps_0/(q_5+\eps_0)} ). $$
Appealing to Proposition \ref{P:PF1n}, and for convenience replacing $\exp(- c\lam^{1/16})$ with the weaker
bound $\exp(-c \lam^{\eps_0})$,
we see the $k$th term for $k \ge 2$ is bounded above by
$$ T^{-\df/\dw} k^{q_5/\eps_0} \Phi^{q_5/(q_5+\eps_0)}
\exp\big( - \Phi^{\eps_0/(q_5 + \eps_0)} ( k^{-q_5/\eps_0} + (k-1) ) \big). $$
Summing this series, the bound follows with $\theta_2 = \eps_0/(q_5 + \eps_0)$. \qed

\subsection{Averaged heat kernel lower  bound}

In this subsection, we will use Theorem \ref{thm:thm5-8} to establish the averaged heat kernel lower bound from Theorem \ref{mainthm3}. The ideas of the following arguments are from \cite[Section 4]{BK}. We first obtain deterministic diagonal and near-diagonal lower bounds that hold on realisations of $\sU$ that occur with suitably high probability. We recall the notation $D^+(a)$ from $4$ lines
below \eqref{ddef}, define $D^-(a)$
analogously for the corresponding part of the UST below $\mathcal{T}_{0}\cap\mathcal{T}_{trunk}$, and set $D^\pm(a):=D^+(a) \cup D^-(a)$.

\begin{lemma}\label{on-dialbqq}
Let $\lam\geq \lam_0$, $m\geq m_0$ and $\lam^7\leq k\leq m^{1/2}\wedge c_1N^4$. Moreover, let  $\tilde{F}_*(k)$ be an event with the properties described in the statement of Theorem \ref{thm:thm5-8} for both ${D}^+(1/2)$ and ${D}^-(1/2)$, and in particular satisfies $\pP_\sT( \tilde{F}_*(k)^c) \le c N e^{-c' k^{1/24}}$. Then there exist constants $c_i,q_i$ such that on $\tilde{F}_*(k)$, if $x\in D^\pm(3/8)$ and $c_2k^{1/4}\leq s\leq c_3 N$, then
\be\label{eq:niebi2-4}
\tilde{p}^{\mathcal{U}}_{n}(x,x) \ge c_4 \lambda^{-3}k^{-5/2}s^{-1}m^{-2}\geq  c_5 \lambda^{-q_1}k^{-q_2}n^{-d_f/d_w}
\ee
for $c_6 \lambda^{-1}k^{-5/2}s^2 m^{\kappa d_w}\le  n \le c_7\lambda^{-1}k^{-5/2}s^2 m^{\kappa d_w}$.
\end{lemma}
\proof This can be obtained by modifying standard arguments. By a line-by-line modification of the proof of \cite[Proposition 4.4.1, 4.4.3]{Kum}, for example, we have on $\tilde{F}_*(k)$ that
\[ c'\lambda^{-1}k^{-5/2}s^2m^{\kappa d_w}\leq E^\sU_x\left(\sigma_{x,sm^\kappa}\right)\le c\lambda s^2m^{\kappa d_w}\]
for all $x\in D^\pm(3/8)$ and $s$ in the given range. The above estimates and the Markov property (see
\cite[Proposition 4.4.3]{Kum}) imply that the following holds on $\tilde{F}_*(k)$,
\begin{eqnarray*}
P^\sU_x\left(\sigma_{x,s_k}>n\right)&\ge& \frac{c'\lambda^{-1}k^{-5/2}s^2m^{\kappa d_w}-n}
{c\lambda s^2m^{\kappa d_w}},
\end{eqnarray*}
for all $x\in D^\pm(3/8)$ and $n\geq 0$. Given this and the upper volume estimate that holds on $\tilde{F}_*(k)$, \eqref{eq:niebi2-4} can be proved as in \cite[Proposition 4.4.4]{Kum}.   \qed

\begin{lemma}\label{near-diag-lower}
(a) Let $\lam\geq \lam_0$, $m\geq m_0$, $\lam^7\leq k\leq m^{1/2}\wedge c_1N^4$ and $\alpha\in(0,1)$.
Moreover, let  $\tilde{F}_*(k)$ be an event as in Lemma \ref{on-dialbqq} that satisfies
$\pP_\sT( \tilde{F}_*(k)^c) \le c N e^{-c' k^{1/24}}$.
Then there exist constants $c'_i,q_i$ such that on $\tilde{F}_*(k)$, we have for  $x\in D^\pm(3/8)$ and $y\in\sU$ satisfying $d_{\sU}(x,y)\le 2s^{\alpha}m^\kappa$ for some $c'_1 \lambda^{4/(1-\alpha)} k^{5/(1-\alpha)}\leq s\leq c'_2 N$,
\be \label{e:ndlb67}
\tilde{p}^\sU_n(x,y) \ge c'_3 \lambda^{-q_1}k^{-q_2}n^{-d_f/d_w}
\hbox{ for $c'_4 \lambda^{-1}k^{-5/2}s^2 m^{\kappa d_w}\le  n \le c'_5\lambda^{-1}k^{-5/2}s^2 m^{\kappa d_w}$}.
\ee
(b) If $x_0\in D(3/8)\cap\mathcal{T}_0$ and $x,y\in B_\sU(x_0,s^{1/2}m^\kappa)$, then the same lower bound holds on an event $\tilde{F}_*(x_0,k,s)$ that satisfies $\pP_\sT( \tilde{F}(x_0,k,s)_*^c) \le c \lambda s e^{-c' k^{1/24}}$.
\end{lemma}
\proof By the discrete-time adaptation of \cite[Lemma 4.3]{BK} (which can be obtained by applying estimates in \cite[Section 4]{llt}) and Lemma \ref{on-dialbqq}, we have
\begin{eqnarray*}
\left|\frac{\tilde{p}^\sU_n(x,y)}{\tilde{p}^\sU_n(x,x)}-1\right|^2\le \frac{c d_\sU(x,y)}{n\tilde{p}^\sU_n(x,x)}
\le \frac{c'\lambda^4 k^5 }{ s^{1-\alpha}}\leq \frac{1}{4}.
\end{eqnarray*}
Hence $|\tilde{p}^\sU_n(x,y)-\tilde{p}^\sU_n(x,x)|\le \tilde{p}^\sU_n(x,x)/2$, so we obtain
\begin{eqnarray*}
\tilde{p}^\sU_n(x,y)\ge \tilde{p}^\sU_n(x,x)-|\tilde{p}^\sU_n(x,y)-\tilde{p}^\sU_n(x,x)|\ge \tilde{p}^\sU_n(x,x)/2
\ge  c \lambda^{-q_1}k^{-q_2}n^{-d_f/d_w},
\end{eqnarray*}
where we used Lemma \ref{on-dialbqq} in the last inequality. This establishes part (a), and part (b) is obtained in the same way, but using Corollary \ref{corest} in place of Theorem \ref{thm:thm5-8}.
\qed

\begin{defn}
{\rm
Let $M,N\in {\mathbb N}$, $\alpha=\frac12$, $\lambda\geq \lambda_0$, $m\geq m_0$, $\mathcal{T}$ be a fixed tree as described after Proposition \ref{P:HN1}, and $x\in \sT_{0}$ with $x\ne 0$. Set $r=d_\sU(0,x)/N$ and let $z_0=0,z_1,\cdots,z_N=x$ be points on the path between $0$ and $x$ with $|d_{\sU}(z_{i-1},z_i)-r|\le 1$ that are chosen in some fixed way. For $i=1,\dots,N$, let $\xi_i$ be the smallest integer $k$ such that $\tilde{F}_*(z_i,k,k^{12})$, $\{|B_\sU(z_{i-1},\frac{1}{4}k^{6}m^\kappa)|\geq m^2\}$ and $\{|B_\sU(z_i,\frac14k^{6}m^\kappa)|\geq m^2\}$ hold. (Set $\xi_i=\infty$ if the requirements are not satisfied.) We then say that $G(q,x,N,M)$ holds if $\sum_{i=1}^N\xi_i^q\le MN$.
}\end{defn}

\begin{propn} It holds that
\[\pP_\sT(G(q,x,N,M))\ge 1-\frac {c_q}M-  cN\lambda \left(m^{6}\wedge N\right)e^{-c' m^{1/48}\wedge (cN^{1/288})},\]
where $c_q$ is a constant that depends on $q$.
\end{propn}
\proof By Corollary \ref{corest} and a simple union bound, it holds that
\[\mathbf{P}_{\mathcal{T}}\left(\xi_i> m^{1/2}\wedge (cN^{1/12})\mbox{ for some }i\right)\leq  N\times c\lambda \left(m^{6}\wedge N\right)e^{-c' m^{1/48}\wedge (cN^{1/288})}.\]
Moreover, Corollary \ref{corest} and the Markov inequality yield that
\[\mathbf{P}_{\mathcal{T}}\left(\sum_{i=1}^N\xi_i^q\mathbf{1}_{\{\xi_i\leq m^{1/2}\wedge (cN^{1/12})\}}> MN\right)\leq \frac{c_q}{M}.\]
Putting these estimates together completes the proof.
\qed

\medskip

\noindent
\emph{Proof of the lower bound of Theorem \ref{mainthm3}.}
For simplicity, we only consider the case when $x=(R,0)$, where $R\in\mathbb{Z}_+$; see Remark \ref{R:genx}
for the modifications necessary for the general case. Let $T \geq R$, and set
$y = R^{\kappa \dw}/T$. We need to consider several cases. These will depend on constants $b,b'\geq 8$, which will be chosen below.

\sm {\em Case 1: $R \le T \le b R$.} Let $F$ be the event that the UST $\sU$ contains the
straight path along the $x_1$-axis between $0$ and $x$. By considering the construction of $\sU$
which starts by running $S^x$ until it hits $0$, we have $\bP(F) \ge 4^{-R}$.
Let $z$ be the point adjacent to $x$ on the path $\gam(0,x)$.
If the event $F$ holds, then $P^\sU_0( X_T \in \{x,z\}) \ge 4^{-T}$, and it follows immediately
that $\bE \tilde p^\sU_T(0,x) \ge 4^{-R+T} \ge e^{-c T}$, which yields the desired lower bound for
these values of $R$ and $T$.

\sm {\em Case 2: $T\geq b'R^{\kappa d_w}$}. To begin with, suppose that $R\geq 1$.
We use Lemma \ref{near-diag-lower}, and take $\lam = \lam_0$.
Given $k$ we set $c_1 N^4= k$, and choose $k=k_0\geq\lam_0$ large enough so that
$\bP_\sT( \tilde F_*(k_0)^c) \le \half$.
We take $s = c'_1 \lam_0^8 k_0^{10}$. As the constants
$\lam,k,s$ do not depend on $R$ or $T$, we can absorb them into constants $c_i$.
The bound \eqref{e:ndlb67} holds for $T \in [c_1 m^{\kappa \dw}, c_2   m^{\kappa \dw}]$, so we
choose  $m$ so that $T$ is in this range; this gives that $m \ge c R$. (NB. By increasing the value of $b'$ if necessary, we can further ensure that $m\geq m_0$ and $m^{1/2}\geq k_0$.) The construction of $ \tilde F_*(k_0)$ in Section \ref{sec:VRest} implies that
on this event $\dU(0,x) \le c s^{1/2} m^\kappa$, and thus we have the lower bound
$$  \tilde p^\sU_T(0,x) \ge c T^{-\df/\dw}. $$
Since $\bP( H_{N_1}) \ge \exp(-c N) \ge \exp( -c' k_0)$ and
$\bP_\sT( \tilde F_*(k_0)) \ge \half$, the averaged lower bound $  \bE \tilde p^\sU_T(0,x) \ge c T^{-\df/\dw}$
follows. If $R=0$, then one can use the same event as for $R=1$ to deduce the result, since one also has that $ \tilde p^\sU_T(0,0) \ge c T^{-\df/\dw}$ on that event.

\sm  {\em Case 3: $bR < T < b' R^{\kappa \dw}$}.
Choose $N,m \in\mathbb{N}$ to satisfy
\[N\le  (b')^{\frac{\kappa d_w}{\kappa d_w-1}}\left(\frac{|x|^{\kappa d_w}}T\right)^{1/(\kappa d_w-1)}< N+1, \q m \le R/N< m+1. \]
Note that $N+1 \ge (b')^{\frac{\kappa d_w}{\kappa d_w-1}}y^{1/(\kappa \dw-1)} \ge 8$, and $R/N \ge b^{1/(\kappa \dw-1)}$.
Hence choosing $b$ large enough we can ensure that $m \ge m_0$ and also that
if $\sT$ is a tree selected in the way described after Proposition \ref{P:HN1}, then
$\bP_\sT(G(43/2,x, N,M))\geq \frac{1}{2}$.
The reason we take $q=43/2$ is that this is the power of $\xi_i$ that arises in the time range
for the estimate \eqref{e:ndlb67}. More precisely
on $G(q,x, N,M)$, for each $i=1,\dots,N$, it holds that
\be \label{e:chain-lb}
\tilde{p}^\sU_{n_i}(z,y) \ge c_1 \lambda^{-q_1} \xi_i^{-q_2}n_i^{-\df/\dw}
\ee
for $c_2\lambda^{-1} \xi_i^{-5/2}(\xi_i^{12})^2m^{\kappa \dw}\le  n_i \le c_3\lambda^{-1} \xi_i^{-5/2}(\xi_i^{12})^2 m^{\kappa \dw}$ and $z,y\in B_\sU(z_i,\xi_i^{6}m^\kappa)$.
Since $24-5/2=43/2$, in
the argument below, we will  need to sum over the quantities $\xi^{43/2}_i$; restricting to the event $G(43/2,x, N,M)$ ensures that we can control this sum.

Now, since it holds that $d_\sU(z_{i-1},z_{i})\leq c \lambda m^\kappa$, the estimate
\eqref{e:chain-lb} includes the case when $z\in B_\sU(z_{i-1},\frac14\xi_i^{6}m^\kappa)$ and $y\in B_\sU(z_{i},\frac14\xi_i^{6}m^\kappa)$. Setting $\tilde{T}:=\sum_{i=1}^Nn_i$, where $n_i$, $i=1,\dots, N$,
satisfy the previous constraints, we then have that
\begin{equation}\label{eq:8-11do20}
c_2 \lambda^{-1}\sum_i\xi_i^{43/2} m^{\kappa \dw}\le
\tilde{T} \le c_3\lambda^{-1} \sum_i\xi_i^{43/2} m^{\kappa \dw}.
\end{equation}
Moreover, writing $B_i:=B_\sU(z_i,\frac14(\xi_{i}^{6}\wedge \xi_{i+1}^{6})m^\kappa)$, we have that
\begin{eqnarray*}
\lefteqn{\tilde{p}^\sU_{\tilde{T}}(0,x)}\\
&\ge &\frac12\sum_{y_1\in B_1}\cdots \sum_{y_{N-1}\in B_{N-1}}\tilde{p}^\sU_{n_1}(0,y_1)
\cdots \tilde{p}^\sU_{n_N}(y_{N-1},x)\mathbf{1}_{\{n_i-d_\sU(y_{i-1},y_i)\mbox{ is even for each }i\}}\\
&\ge & \prod_{j=1}^{N-1}\Big(|B_j|\, c_1 \lambda^{-q_1} \xi_j^{-q_2}n_j^{-\df/\dw}\Big)
 c_1 \lambda^{-q_1} \xi_N^{-q_2}n_N^{-\df/\dw}\\
&\ge &c_\lambda\tilde{T}^{-d_f/d_w}\exp\left(-c'_\lambda N-(q_2+(43/2)d_f/d_w)\sum_{i=1}^N \log \xi_i-\frac{d_f}{d_w}\log\left(\sum_{i=1}^N \xi_i^{43/2}\right)\right)\\
&\ge& c_\lambda\tilde{T}^{-d_f/d_w}\exp(-c_\lambda'N),
\end{eqnarray*}
where in the last inequality we used $\sum_{i=1}^N \log \xi_j\le \sum_{i=1}^N\xi^{43/2} \le MN$.
Note that in \eqref{eq:8-11do20}, we may take $c_2>0$ as small as we like. (This is because $\tilde{p}^\sU_n(x,x)$ is monotone decreasing, which means we can take $c_6$ in Lemma \ref{on-dialbqq} as small as desired. Moreover, we can take $c_4'$ in Lemma \ref{near-diag-lower} to match this). In particular, taking $c_2\leq 8^{-\kappa d_w}M^{-1}$, we obtain
\[c_2 \lambda^{-1}\sum_i\xi_i^{43/2} m^{\kappa \dw}\leq 8^{-\kappa d_w}Nm^{\kappa \dw}\le T.\]
Hence we may take $\tilde{T}\le T$. If $T\le c_3 \lambda^{-1} \sum_i\xi_i^{43/2} m^{\kappa \dw}$, then we can choose $n_j$ so that $T=\tilde{T}$. If not, let $T'=T-\tilde{T}\le  Nr^{d_w}$. Let $j_0$ be such that $\xi_{j_0}$ is minimal, and add $N'$ extra steps between $B_{j_0-1}$ and $B_{j_0}$ in the chaining argument above, each with time length satisfying the constraint of $n_{j_0}$ and the total time of the additional steps is equal to $T'$. The latter constraints readily imply that $N'\leq cN$, and we further observe that each extra step contributes a factor of $c_\lambda \xi_{j_0}^{-(q_2+(43/2)d_f/d_w)}$ to the lower bound. Thus the total contribution is no less than $e^{-c_\lambda N}$. Taking the average over $G(43/2,x,N,M)$ and $\mathcal{T}$, we obtain the result.
\qed

\begin{remark} \label{R:genx}
{\rm For a general $z = (z_1, z_2) \in \bZ^2$ we need to replace the tree $\sT$ defined in Section \ref{sec:VRest}
by a tree which connects $0$ and $z$. We replace the `S-shaped' path $(\tilde x_i)_{i=0}^{N_1}$
defined just after Lemma \ref{L:basicLP} with a path for the which the central section has `L' shape
which connects $0$ and $(z_1/m, z_2/m)$, and the rest of the path shields the central section from
the remainder of $\bZ^2$.
  The estimates of Section \ref{sec:VRest} and \ref{sec:ann-b} all work for this path, and the proof of the
lower bound on $\bE \tilde p^\sU_n(0,z)$ then follows.
} \end{remark}

\section{Failure of the elliptic Harnack inequalities}\label{failEHI}

The aim of this section is to make precise and prove Corollary \ref{ehicor}. We start by giving the definition of the elliptic Harnack inequality that we consider, as well as a related metric doubling property.

\begin{defn}\label{ehidef}
{\rm Let $(X_\omega,d_\omega, \mu_\omega)$ be a weighted random graph.\\
(i)  We say that the \emph{large scale elliptic Harnack inequalities (LS-EHI)} hold (for the random walk associated with $(X_\omega,d_\omega, \mu_\omega)$) if there exists a deterministic constant $C>1$ and, for each $x_0\in X_\omega$, there exists an $R_{1,x_0}(\omega)>0$ such that the following inequality is satisfied
\begin{equation*}
\sup_{B_{d_\omega} (x_{0},R) }{u}\leq C\inf_{B_{d_\omega} (x_{0},R) }{u}.
\end{equation*}
for any $x_0\in X_\omega$, $R\ge R_{1,x_0}(\omega)$ and any non-negative bounded harmonic function $u$ on  $B_{d_\omega}( x_{0},2R)$.\\
(ii) We say that the \emph{large scale metric doubling property (LS-MD)} holds if there exists a deterministic constant $M\in {\mathbb N}$ and, for each $x_0\in X_\omega$, there exists $R_{2,x_0}(\omega)>0$ such that, for any $x_0\in X_\omega$ and $R\ge R_{2,x_0}(\omega)$, $B_{d_\omega}(x_0,R)$ can be covered by $M$ balls of radius $R/2$.}
\end{defn}

The  main result in this section is the following.

\begin{thm}\label{thm:MD-LSthm}
(LS-EHI) does not hold for the random walk on $\sU$.
\end{thm}

For the proof, we use the following proposition.

\begin{prop}\label{thm:MD-LSpr}
(LS-EHI) implies (LS-MD).
\end{prop}

\proof The proof is a line-by-line modification of
\cite[Theorem 3.11]{BM}. Hence we omit it.
\qed

The following lemma will be used to check that (LS-MD) is violated for $\sU$.

\begin{lemma} \label{L:noMD}
There exists a constant $\delta>0$ such that, $\mathbf{P}$-a.s., one can find a divergent sequence $(R_n)_{n\geq 1}$ for which there exist at least $n$ disjoint $d_\sU$-balls of radius $\delta R_n$ contained in $B_\sU(0,R_n)$.
\end{lemma}

\proof Let $(G(i))_{i\geq 1}$ be the events described in the proof of Theorem \ref{mainthm1},
where it was shown that $G(i)$ holds infinitely often, $\mathbf{P}$-a.s. Now, let
$(z_j)_{j=1}^{\varepsilon(\log i)^{1/2}}$ be the vertices at the centres of the top
row of boxes in the configuration shown in Figure \ref{grid} for $N=D_i/m_i$ and $m=m_i$.
On $G(i)$, we have that
\[d_\sU(0,z_j)\leq C\left(\frac{D_i}{m_i}\right)m_i^{\kappa},\qquad \forall j=1,2,\dots,\varepsilon(\log i)^{1/2},\]
and also
\[d_\sU(z_j,L_i)\geq c\left(\frac{D_i}{m_i}-2\right)m_i^{\kappa},\qquad \forall j=1,2,\dots,\varepsilon(\log i)^{1/2},\]
where $L_i$ is the bottom row of boxes in the configuration shown in Figure \ref{grid}.
It readily follows that there exist at least $\varepsilon(\log i)^{1/2}$ disjoint $d_\sU$-balls of
radius $\frac{c}{2}\left(\frac{D_i}{m_i}-2\right)m_i^{\kappa}$ contained in
$B_\sU(0,C\left(\frac{D_i}{m_i}\right)m_i^{\kappa})$. Hence taking
\[R_n=C\left(\frac{D_i}{m_i}\right)m_i^{\kappa},\qquad \delta=\frac{c}{4C},\]
where $i=e^{(n/\varepsilon)^2}$ yields the result.
\qed

\noindent
\emph{Proof of Theorem \ref{thm:MD-LSthm}.} Let $\delta>0$, and suppose that LS-MD holds. Then there exists a constant $M' = M'(M,\delta)$ such that: for each $x \in X_\omega$, there exists $R'_{x,\omega} < \infty$, such that if $R \ge R'_{x,\omega}$, then
the ball $B_{\sU}(x,R)$
can be covered by $M'$ balls of radius $\delta R$. However, Lemma \ref{L:noMD} shows that this fails for $\sU$. Hence Proposition \ref{thm:MD-LSpr} yields the result.
\qed

\section{Scaling limits}\label{sec:scaling}

In this section, we prove the results stated in the introduction concerning scaling limits of the random walk, namely Theorems \ref{mainthm2} and \ref{denslimit}, and Corollaries \ref{cor2} and \ref{cor3}.
\bigskip

\noindent
{\em Proof of Theorem \ref{mainthm2}.} By the separability of the Gromov-Hausdorff-vague topology (see, for example, \cite[Proposition 5.12]{ALWGap}), it is possible to suppose that we have a sequence $(\sU_n)_{n\geq 1}$ of copies of $\sU$, all built on the same probability space, so that
\[\left(\sU_n,n^{-\kappa}d_{\sU_n},n^{-2}\mu_{\sU_n},0\right){\rightarrow}\left(\sT,d_\sT,\mu_\sT,\rho_\sT\right)\]
holds $\mathbf{P}$-a.s. (Note that for this part of the article, we do not need the spatial embeddings into $\mathbb{R}^2$.) It follows from \cite[Proposition 5.9]{ALWGap} that, $\mathbf{P}$-a.s., there exists a metric space $(M,d_M)$ so that the spaces $(\sU_n,n^{-\kappa}d_{\sU_n})$, $n\geq 1$, and $(\sT,d_\sT)$ can be isometrically embedded into $(M,d_M)$ in such a way that: $0$ and $\rho_{\sT}$ are mapped to a common point, $0_M$ say; the embedded measures $n^{-2}\mu_{\sU_n}$ converge vaguely to the embedded version of $\mu_\sT$; and, for all but countably many $r$, the sets $\sU_n\cap \bar{B}_M(0_M,r)$, where $\bar{B}_M(0_M,r)$ is the closed ball in $M$ of radius $r$ centred at $0_M$, converge to $\sT\cap \bar{B}_M(0_M,r)$ with respect to the Hausdorff distance between compact subsets of $(M,d_M)$. As a consequence (see, \cite[Theorem 7.1]{Cr}), we moreover have that the laws of the random walks $(X^{\sU_n}_{tn^{\kappa+2}})_{t\geq0}$ converge weakly to the law of $(X^\sT_t)_{t\geq0}$, when these are considered as measures on $D(\mathbb{R}_+,M)$. Consequently, we have that the Assumptions 1 and 5 of \cite{CrHa} are satisfied (actually Assumption 1 requires the convergence of measures of balls under the various laws, but this condition is readily relaxed to the requirement that the balls in question are continuity sets for the limiting measure), and hence we can apply \cite[Theorem 1 and Proposition 14]{CrHa} to deduce that the associated transition densities satisfy, $\mathbf{P}$-a.s.,
\begin{equation}\label{hkllt}
\left(n^{2} \tilde{p}_{\lfloor tn^{2+\kappa}\rfloor }^{\sU_n}(0,0)\right)_{t>0}{\rightarrow}\left( p_{t}^\sT(\rho_\sT,\rho_\sT)\right)_{t>0}.
\end{equation}
Reparameterising this, the first part of the theorem follows.

In view of the distributional limit we have just proved, to prove the scaling limit at \eqref{odhkscale} it will suffice to check the following integrability condition: for any $p\geq 1$, there exists a constant $C\in(0,\infty)$ such that
\begin{equation}\label{ofmoments}
\sup_{n\geq 1} n^{d_f/d_w}\left\|\tilde{p}_{n}^\sU\left(0,0 \right)\right\|_{p}\leq C,
\end{equation}
where $\|\cdot\|_p$ is the $L_p$ norm with respect to $\mathbf{P}$. Now, by Lemma \ref{L:F1easyub}, on the event $F_1(\lambda, n^{1/d_w\kappa})$, it holds that $\tilde{p}_{n}^\sU(0,0 )\leq c\lambda^q n^{-d_f/d_w}$. Hence, if $\Lambda_n :=\inf\{\lambda\geq 1:\:F_1(\lambda, n^{1/d_w\kappa})\mbox{ holds}\}$, then
\begin{equation}\label{lpbound}
n^{d_f/d_w}\left\|\tilde{p}_{n}^\sU\left(0,0 \right)\right\|_{p}\leq c\left\|\Lambda_n^q\right\|_{p}.
\end{equation}
Since Proposition \ref{P:PF1n} yields that the right-hand side above is uniformly bounded in $n$, this completes the proof. \qed

In preparation for the proof of Theorem \ref{denslimit}, we verify the equicontinuity of the averaged heat kernel under scaling.

\begin{propn}\label{equicont} There exists a constant $C\in(0,\infty)$ such that
\[\sup_{n\geq 1}n^{{d_f}/{d_w}}\left|\pE \tilde{p}_{\lfloor tn\rfloor }^\sU\left(0,[x n^{\frac{1}{\kappa d_w}}]\right)
-\pE \tilde{p}_{\lfloor tn\rfloor }^\sU\left(0,[y n^{\frac{1}{\kappa d_w}}]\right)\right|\leq Ct^{-d_f/2d_w}|x-y|^{\kappa/2},\]
for all $x,y\in\mathbb{R}^2$, $t>0$.
\end{propn}
\proof From \cite[Lemmas 9 and 10]{CrHa}, we have for every $x,y\in\mathbb{Z}^2$ and $n\geq 1$ that
\begin{equation}\label{ffffff}
\left( \tilde{p}_{n}^\sU\left(0,x \right)
- \tilde{p}_{n}^\sU\left(0,y \right)\right)^2\leq \frac{2d_\mathcal{U}(x,y)\tilde{p}_{2\lceil n/2\rceil}^\sU\left(0,0 \right)}{n}.
\end{equation}
Hence Jensen's and H\"{o}lder's inequalities yield that, for any $\varepsilon>0$,
\[\left| \pE\tilde{p}_{n}^\sU\left(0,x \right)
-\pE \tilde{p}_{n}^\sU\left(0,y \right)\right|\leq \sqrt{\frac{2}{n}}\left\|d_\mathcal{U}(x,y)^{1/2}\right\|_{1+\varepsilon}
\left\|\tilde{p}_{2\lceil n/2\rceil}^\sU\left(0,0 \right)^{1/2}\right\|_{\frac{1+\varepsilon}{\varepsilon}},\]
where we again write
$\|\cdot\|_p$ for the $L_p$ norm with respect to $\mathbf{P}$.
Now, by Theorem \ref{T:dxy},
it holds that, for suitably small $\varepsilon$,
\begin{equation}\label{dest1}
\left\|d_\mathcal{U}(x,y)^{1/2}\right\|_{1+\varepsilon}\leq C|x-y|^{\kappa/2}.
\end{equation}
Moreover, from \eqref{lpbound} (and Proposition \ref{P:PF1n}), we have that
\[\left\|\tilde{p}_{2\lceil n/2\rceil}^\sU\left(0,0 \right)^{1/2}\right\|_{\frac{1+\varepsilon}{\varepsilon}}\leq C n^{-d_f/2d_w}.\]
Since $n^{{d_f}/{d_w}}\times n^{-1/2} \times (n^{\frac{1}{\kappa d_w}})^{\kappa/2}\times n^{-d_f/2d_w}=1$, combining these estimates readily yields the result.
\qed

We moreover note the following rerooting invariance property of the limiting tree.

\begin{propn}\label{reroot} (a) For any $x\in\mathbb{R}^2$,
\[\mathbf{P}\left(\left|\phi_{\sT}^{-1}(x)\right|>1\right)=0.\]
(b) For any $x\in\mathbb{R}^2$,
\begin{equation}\label{gggggg}
\left(\sT,d_\sT,\mu_\sT,\phi_\sT-x,\phi_{\sT}^{-1}(x)\right)\buildrel{d}\over{=}\left(\sT,d_\sT,\mu_\sT,\phi_\sT,\rho_\sT\right).
\end{equation}
\end{propn}
\proof We first prove the result of part (a) for $x\in\mathbb{R}^2\backslash\{0\}$. In particular, by the scale and rotational invariance properties of \eqref{scaleinvar} and \eqref{rotinvar}, respectively, we have that $\mathbf{P}(|\phi_{\sT}^{-1}(x)|>1)$ is a constant for $x\in\mathbb{R}^2\backslash\{0\}$. Moreover, as was noted in the proof of \cite[Theorem 1.3]{BCK}, we know that the Lebesgue measure of $\{x:\:|\phi_{\sT}^{-1}(x)|>1\}$ is zero, $\mathbf{P}$-a.s. Hence, it follows from Fubini's theorem that $\mathbf{P}(|\phi_{\sT}^{-1}(x)|>1)=0$ for all $x\in\mathbb{R}^2\backslash\{0\}$.

We next prove part (b) for  $x\in\mathbb{R}^2\backslash\{0\}$. To begin with, we note from part (a) that the left-hand side of \eqref{gggggg} is a well-defined measured, rooted spatial tree, $\mathbf{P}$-a.s. Moreover, by the separability of the Gromov-Hausdorff-type topology that we are considering (see \cite[Proposition 3.4]{BCK}), it is possible to suppose that we have realisations of the relevant random objects built on a common probability space so that
\[\left(\sU,\delta_n^{\kappa}\dU,\delta^{2}_n\mu_\sU,\delta_n\phi_\sU,0\right){\rightarrow}\left(\sT,d_\sT,\mu_\sT,\phi_\sT,\rho_\sT\right),\]
almost-surely as $n\rightarrow\infty$ (cf.\ the proof of Theorem \ref{mainthm2}). It follows that it is almost-surely possible to choose a (random) $x_n^R\in\delta_n\mathbb{Z}^2$ such that
\[\left(\sU,\delta_n^{\kappa}\dU,\delta^{2}_n\mu_\sU,\delta_n\phi_\sU-x_n^R,x_n^R\right){\rightarrow}\left(\sT,d_\sT,\mu_\sT,\phi_\sT-x,\phi_\sT^{-1}(x)\right).\]
In particular, this implies that $|x_n^R-x|\rightarrow 0$, almost-surely. Moreover, let $x_n\in\delta_n\mathbb{Z}^2$ be a deterministic sequence such that $|x_n-x|\rightarrow 0$. One can then deduce from \cite[Theorem 1.1]{BM11} (and the Borel-Cantelli lemma) that there exists a deterministic subsequence $n_i$ along which $\delta_{n_i}^\kappa d_\sU(x_{n_i},x_{n_i}^R)\rightarrow 0$, almost-surely. Hence we find that, almost-surely,
\[\left(\sU,\delta_{n_i}^{\kappa}\dU,\delta^{2}_{n_i}\mu_\sU,\delta_{n_i}\phi_\sU-x_{n_i},x_{n_i}\right){\rightarrow}\left(\sT,d_\sT,\mu_\sT,\phi_\sT-x,\phi_\sT^{-1}(x)\right).\]
By the translation invariance of $\sU$ (see \cite[Theorem 2.3]{Pem91}), the left-hand side here has the same distribution as $(\sU,\delta_{n_i}^{\kappa}\dU,\delta^{2}_{n_i}\mu_\sU,\delta_{n_i}\phi_\sU,0)$, which we know converges in distribution to $(\sT,d_\sT,\mu_\sT,\phi_\sT,\rho_\sT)$, and so the result follows.

Finally, let $x\in\mathbb{R}^2\backslash\{0\}$. Then, from part (b) (for such $x$), we know that $|\phi^{-1}_{\sT}(0)|$ is equal in distribution to $|\phi^{-1}_{\sT}(x)|$. And, from part (a) (again, for such $x$), we know the latter is $\mathbf{P}$-a.s.\ equal to 1. In particular, we find that $\phi^{-1}_{\sT}(0)=\{\rho_\sT\}$, $\mathbf{P}$-a.s. Hence both part (a) and part (b) are readily extended to include the point $x=0$.
\qed

\noindent
{\em Proof of Theorem \ref{denslimit}.} From \cite[Theorem 1.4]{BCK}, we know that, under the averaged law $\int P^\sU_0(\cdot)\mathbf{P}(d\sU)$,
\begin{equation}\label{rwconv}
\left(n^{-1/d_w\kappa}X^{\sU}_{tn}\right)_{t\geq 0}\buildrel{d}\over{\rightarrow} \left(\phi_{\sT}\left(X^{\sT}_t\right)\right)_{t\geq 0}.
\end{equation}
Applying this in conjunction with Theorem \ref{mainthm2} (specifically \eqref{odhkscale}) and Proposition \ref{equicont}, elementary analysis arguments yield that, for each fixed $t\in(0,\infty)$, $\phi_{\sT}(X^{\sT}_t)$ admits a density $q_t(x)\in C(\mathbb{R}^2,\mathbb{R})$ satisfying the convergence result of part (c). From this, part (a) of the theorem is a simple consequence of Proposition \ref{equicont}. Moreover, given the continuity of the density $q_t$ in the spatial variable, part (b) follows from the scale and rotational invariance properties at \eqref{scaleinvar} \eqref{rotinvar}.

For part (d), we again recall from the proof of \cite[Theorem 1.3]{BCK} that the Lebesgue measure of $\{x:\:|\phi_{\sT}^{-1}(x)|>1\}$ is zero, and also from the latter result that $\mu_{\sT}=\mathcal{L}\circ\phi_\sT$, where $\mathcal{L}$ is two-dimensional Lebesgue measure. Putting these observations together yields
\begin{eqnarray*}
\int_B q_t(x) dx&=&\pE\left(P^\sU_0( X^\mathcal{T}_t \in\phi_\mathcal{T}^{-1}( B ))\right)\\
& =&\pE\left( \int_{\phi_\sT^{-1}(B)}p_{t}^\sT(\rho_\sT,x)\mu_\sT(dx)\right)\\
& =&\pE\left( \int_{B}p_{t}^\sT(\rho_\sT,\phi^{-1}_\sT(x))\mathbf{1}_{\{|\phi^{-1}_\sT(x)|=1\}}dx\right)\\
&=&\int_{B}\pE\left( p_{t}^\sT(\rho_\sT,\phi^{-1}_\sT(x))\right)dx,
\end{eqnarray*}
for all Borel $B \subseteq \bR^2$, where we have applied Fubini's theorem and Proposition 8.2(a) to obtain the final equality. It follows that the desired equality holds for Lebesgue almost-every $x$, and so to complete the proof, it will suffice to show that, for each fixed $t$, $p_t(x):=\pE( p_{t}^\sT(\rho_\sT,\phi^{-1}_\sT(x)))$ is continuous in $x$. Now, from the rotational invariance of \eqref{rotinvar}, we have that $p_t(x)$ is constant on circles centred at the origin. And thus, to check continuity at $x\neq 0$, it will suffice to show that $p_t(\lambda x)\rightarrow p_t(x)$ as $\lambda\rightarrow 1$. Moreover, by the scale invariance property \eqref{scaleinvar}, this is equivalent to checking that $p_{\lambda t}(x)\rightarrow p_t(x)$ as $\lambda\rightarrow1$, and doing this is our next aim.  Arguing as in the proof of \cite[Theorem 10.4]{Kigq}, for example, and applying the monotonicity of the on-diagonal part of the heat kernel, one can deduce that, for $s,t>r$,
\[\left|p_{s}^\sT(\rho_\sT,\phi^{-1}_\sT(x))-p_{t}^\sT(\rho_\sT,\phi^{-1}_\sT(x))\right|\leq 2r^{-1}|t-s|\sqrt{p_{r/2}^\sT(\rho_\sT,\rho_\sT)p_{r/2}^\sT(\phi^{-1}_\sT(x),\phi^{-1}_\sT(x))}.\]
From this, the Cauchy-Schwarz inequality, the rerooting invariance of Proposition \ref{reroot}(b), and \eqref{odhkscale}, we see that
\[\left|p_s(x)-p_t(x)\right|\leq 2r^{-1}|t-s|p_{r/2}(0)=Cr^{-1-d_f/d_w}|t-s|,\]
which implies that $p_{\lambda t}(x)\rightarrow p_t(x)$ as $\lambda\rightarrow1$, as desired. To deal with the case $x=0$, we again argue as in the proof of \cite[Theorem 10.4]{Kigq}, for example (cf.\ \eqref{ffffff}), to deduce that
\[\left|p_{t}^\sT(\rho_\sT,\rho_\sT)-p_{t}^\sT(\rho_\sT,\phi^{-1}_\sT(x))\right|\leq t^{-1}\sqrt{d_\mathcal{T}(\rho_\sT,\phi_{\sT}^{-1}(x))p_{t}^\sT(\rho_\sT,\rho_\sT)}.\]
This implies
\[\left|p_{t}(0)-p_{t}(x)\right|\leq t^{-1}\left\|d_\mathcal{T}(\rho_\sT,\phi_{\sT}^{-1}(x))^{1/2}\right\|_{1+\varepsilon}\left\|p_{t}^\sT(\rho_\sT,\rho_\sT)^{1/2}\right\|_{\frac{1+\varepsilon}{\varepsilon}}.\]
From \eqref{hkllt} and \eqref{ofmoments}, we have that the term $\|p_{t}^\sT(\rho_\sT,\rho_\sT)^{1/2}\|_{\frac{1+\varepsilon}{\varepsilon}}$ is finite for any $\varepsilon>0$. Moreover, arguing as in the proof of Proposition \ref{reroot}, for each $x$, one has that there exists a sequence $(x_{n})$ such that $|x_{n}-x|\rightarrow 0$ and, along a subsequence $(n_i)$,
\[n_i^{-\kappa}d_{\sU}(0,x_{n_i})\buildrel{d}\over{\rightarrow}d_\mathcal{T}(\rho_\sT,\phi^{-1}_{\sT}(x)).\]
Hence from \eqref{dest1} we obtain that $\|d_\mathcal{T}(\rho_\sT,\phi^{-1}(x))^{1/2}\|_{1+\varepsilon}\leq C|x|^{\kappa/2}$, where the constant does not depend on $x$. In particular, these estimates imply that $p_t(x)\rightarrow p_t(0)$ as $|x|\rightarrow 0$, and so the proof is complete.
\qed

\noindent
{\em Proof of Corollary \ref{cor2}.} This is an easy application of Theorems \ref{mainthm3} and \ref{denslimit}(c).
\qed

\noindent
{\em Proof of Corollary \ref{cor-dist}.}
We begin with the bounds for $|X^\sU_n|$.
Integrating the upper bound of Theorem \ref{mainthm3}, we find that
\[ {n^{-p/d_w\kappa}} \bE ( E_0^\sU |X^{\sU}_{n} |^p )
\le n^{-p/d_w\kappa}\sum_{x\in\mathbb{Z}^2}c_1n^{-{d_f}/{d_w}}|x|^p\exp\left\{-c_2  \left(\frac{|x|^{\kappa {d_w}}}{n}\right)^{\frac{\theta_2}{{d_w}-1}}\right\}
\leq c_3,\]
as required. The lower bound follows in a similar fashion.

For the upper bound on $\dU(0, X^\sU_n)$
set $R_k =  \lceil e^k n ^{1/d_w}\rceil$, $B_k = B_\sU(0,  R_k)$ for $k \ge 0$, $B_{-1} =\emptyset$,
and $D_k = B_k \setminus B_{k-1}$ for $k \ge 0$.
Let $k_0 = ((d_w-1)/d_w) \log n$. Note that if $k > k_0$ then $R_k>n$, so that $P_0^\sU ( X^\sU_n \in D_k) =0$. (Recall we are looking at the discrete time walk.) Write $\sigma_k = \sigma_{0,R_k+1}$, where $\sigma_{x,r}$ was defined at \eqref{sigxr}. We then have
\begin{align}
  \bE E_0^\sU d_\sU(0, X^\sU_n)^p
     &\le \bE \sum_{k=0}^{k_0} 2^pe^{ pk} n^{p/d_w} P_0^\sU ( X^\sU _n \in D_k) \nonumber\\
     &\le   2^pn^{p/d_w} \sum_{k=0}^{k_0}  e^{pk} \bE P_0^\sU ( \sigma_k \le n).\label{sumtobound}
\end{align}
Now, by an almost identical argument to Lemma \ref{L:tailhit}, it is possible to check that on the event $F_1(\lambda_k,\lam_k R_k^{1/\kappa})$ with $\lambda_k:=(4k)^{40}$ we have
\[  P_0^\sU ( \sigma_k \le n) \le C \exp( - c \lam_k^{-q_4} m_k)=C\exp(-c\lam_k^{-q}e^{kd_w/(d_w-1)});\]
here, $m_k:=(c_3\lambda_k^{-q_3})^{1/(d_w-1)}\Phi(R_k^{d_w}/n)$ represents the number of steps into which the stopping time is decomposed, where $c_3,q_3,q_4$ are as in \eqref{e:smallexit}. Hence, by Proposition \ref{P:PF1n},
\[ \bE  P_0^\sU( \sigma_k \le n)\le \pP\left(  F_1(\lambda_k,\lam_k R_k^{1/\kappa})^c\right) + C\exp(-c\lam_k^{-q}e^{kd_w/(d_w-1)}) \leq Ce^{-ck^{40/16}},\]
which implies that the sum in \eqref{sumtobound} is finite, and so establishes
the upper bound. The lower bound is proved by the same argument as is used in \cite[Theorem 4.4]{BM11}.
\qed

\noindent
{\em Proof of Corollary \ref{cor3}.}
 From \eqref{rwconv} we have under the averaged law that
\[\left(n^{-1/d_w\kappa}\left|X^{\sU}_{tn}\right|\right)_{t\geq 0}\buildrel{d}\over{\rightarrow} \left(\left|\phi_{\sT}\left(X^{\sT}_t\right)\right|\right)_{t\geq 0}.\]
Part (a) now follows using the uniform integrability given by Corollary \ref{cor-dist}.

For part (b), we start by noting that the convergence at \eqref{scaling} implies that the same result holds if $\delta\phi_\sU$ is replaced by the map $x\mapsto \delta^{\kappa}d_{\sU}(0,x)$, and $\phi_\sT$ is replaced by the map $x\mapsto d_{\sT}(\rho_\sT,x)$. As a consequence, in place of the random walk convergence result of \eqref{rwconv}, one obtains that
\[\left(n^{-1/d_w}d_{\sU}(0,X^{\sU}_{tn})\right)_{t\geq 0}\buildrel{d}\over{\rightarrow} \left(d_{\sT}(\rho_\sT,X^{\sT}_t)\right)_{t\geq 0}.\]
(Concretely, apply \cite[Theorem 7.2]{Cr}.)
Part (b) then also follows from Corollary \ref{cor-dist}.
\qed

\begin{remark}\label{zqrem}
{\rm
Let $\mathcal{R}:=\{x:\:|\phi^{-1}_\sT(\{x\})|=1\}$.
With $\mathbf{P}$-probability one, we have that $\mathcal{L}(\mathcal{R}^c)=0$, where we again use $\mathcal{L}$ to denote Lebesgue measure on $\mathbb{R}^2$, and moreover $0\in\mathcal{R}$ (see Proposition \ref{reroot} and its proof). Since $\mu_\sT(\phi^{-1}_\sT(\mathcal{R}^c))=\mathcal{L}(\mathcal{R}^c)$ (by \cite[Theorem 1.3]{BCK}), it follows that, $\mathbf{P}$-a.s., for any $x\in\mathcal{R}$ and $t\geq 0$,
\[P^\mathcal{T}_{\phi^{-1}_\sT(x)}\left(\phi_\sT(X^\mathcal{T}_t)\in\mathcal{R}\right)=\int_{\phi_\sT^{-1}(\mathcal{R})}p_t^\mathcal{T}(\phi_\sT^{-1}(x),y)\mu_\sT(dy)=1,\]
where $P^\mathcal{T}_{\phi^{-1}_\sT(x)}$ is the quenched law of $X^\mathcal{T}$ started from $\phi_\sT^{-1}(x)$. It readily follows that, when started from $x\in\mathcal{R}$ (including from $x=0$), $\phi_\sT(X^\mathcal{T})$ is a Markov process, and moreover has transition density that is determined by $(p^\mathcal{T}_t(\phi^{-1}_\sT(y),\phi^{-1}_\sT(z)))_{y,z\in\mathcal{R}}$ (and which is defined arbitrarily elsewhere). On the other hand, if $\tau$ is a stopping time for $\phi_\sT(X^\mathcal{T})$ such that $P^\mathcal{T}_{\phi^{-1}_\sT(x)}(\phi_\sT(X^\mathcal{T}_\tau)\in\mathcal{R}^c)>0$, then it is clear that the quenched law of $(\phi_\sT(X^\mathcal{T}_{\tau+t}))_{t\geq0}$ does not only depend on $\phi_\sT(X^\mathcal{T}_{\tau})$, and so $\phi_\sT(X^\mathcal{T})$ is not strong Markov. Indeed, the situation is somewhat similar to that of reflecting Brownian motion in a planar domain with a slit removed (cf.\ comments in \cite[Section 3]{BZ}), though the slit is replaced in our case by the dense set $\mathcal{R}^c\subseteq\mathbb{R}^2$, which we note coincides with the `dual trunk' studied in \cite[Section 10]{Schramm}.
} \end{remark}

\appendix

\section{Appendix: Short LERW paths}

In this section we improve the estimates in \cite{BM10} to prove Theorem \ref{T:LERW-lb}. We begin by considering the following situation, which is described in terms of parameters $m,n,N\in\mathbb{N}$ satisfying $4\le n \le m \le m+2n \le N$, cf. \cite[Definition 1.4]{BM10}. Let $B_m = \BE(0,m)$, $B_N=\BE(0,N)$,
and $x \in \pd_R B_m$, where for a square $B$ we write $\pd_R B$
for the right-hand side of the interior boundary of $B$.
Moreover, let $x_1= x + (\frac{n}{2},0)$, and define $A_n(x) = \BE(x,n/4)$.
Finally, we also suppose we are given a subset $K \subseteq B_m$ that contains a path in
$B_m$ from $0$ to $x$. Importantly, we note that the latter assumption was not made in \cite{BM10};
it is the key to removing the terms in $\log(N/n)$ in \cite[Lemmas 4.6 and 6.1, and Propositions 6.2 and 6.3]{BM10}.
We also remark that in \cite{BM10} the balls $B_n$ and $B_N$ were in the $\ell_2$ norm on $\bZ^2$ rather
than the $\ell_\infty$ norm, but this makes no essential difference to the arguments.

The first result of the section concerns the Green's function $G$ of a simple random walk $S$ on $\mathbb{Z}^2$.
Given a subset $A\subsetneq \mathbb{Z}^2$, we write $G_A(y,z)$ for the expected number of visits that $S$ makes to
$z$ when it starts at $y$ up until it exits $A$. In the proof, we write $\pP_x$ for the law of the random walk started from
$x$, and $\bE_x$ for the corresponding expectation.

\begin{lemma} \label{L:GBNK} There exist constants $c_i$ such that, for $y,z \in A_n(x)$,
\be \label{e:GBNK}
 c_1 \log \left( \frac{n}{1\vee |y-z|}\right) \le G_{B_N \backslash K }(y,z) \le  c_2  \log \left( \frac{n}{1\vee |y-z|}\right).
\ee
\end{lemma}

\proof Set $A_1= B_{5n/16}(x_1)$ and $A_2= B_{3n/8}(x_1)$. We note that
\[ c_1 \log \left( \frac{n}{1\vee |y-z|}\right) \le G_{A_1}(y,z)  \le G_{A_2}(y,z)\le  c_2  \log \left( \frac{n}{1\vee |y-z|}\right).\]
(Cf.\ The applications of results from \cite[Chapter 6]{LLRW} that appear as \cite[Proposition 2.4]{BM10}.)
Hence, since $G_{B_N\backslash K }(y,z) \ge G_{A_1}(y,z)$, the lower bound is immediate. For the upper bound,
writing $T_A$ and $\tau_A$ for the hitting and exit time of a subset $A\subseteq\mathbb{Z}^2$ by the simple random walk $S$, respectively, we have
\begin{align*}
  G_{B_N \backslash K }(y,z) &=   G_{A_2}(y,z) + \bE_y\left( G_{B_N \backslash K }( S_{\tau_{A_2}},z)\right) \\
 &\le c_2  \log \left( \frac{n}{1\vee |y-z|}\right)
 + \max_{w \in \pd A_2} \pP_w( T_{A_1} < T_K) \max_{w' \in \pd A_1} G_{B_N \backslash K }(w',z).
\end{align*}
By the discrete Harnack inequality (see \cite[Theorem 6.3.9]{LLRW}, for example)
and the fact that $K$ contains a path from $x$ to 0, we have that $\pP_w( T_{A_1} < T_K) \le 1-c_3$ for all $w \in \pd A_2$.
Further, for $w' \in \pd A_1$ we have
$$ \pP_{w'}( T_z < \tau_{A_2}) \le  1\wedge\frac{c_4}{\log n} . $$
(Again, cf.\ \cite[Proposition 2.4]{BM10}.)
Combining these estimates gives $G_{B_N \backslash K }(z,z) \le c_2  \log (n)
 +(1-c_3)G_{B_N \backslash K }(z,z)$, and thus $G_{B_N \backslash K }(z,z)\leq \frac{c_2}{c_3}\log(n)$. Hence
\begin{align*}
  G_{B_N \backslash K }(y,z) &\le c_2  \log \left(\frac{n}{1\vee |y-z|}\right)
 +c_5,
\end{align*}
which yields the bound \eqref{e:GBNK}. \qed

Next, let $\tilde{S}$ be a random walk started at $x$ and conditioned to leave $B_N$ before its first return to $K$. We write $\tilde{G}(\cdot, \cdot)$ for the Green's function of $\tilde{S}$.

\begin{lemma}[{Cf. \cite[Lemma 4.6]{BM10}}] \label{L:newL46} There exist constants $c_i$ such that, for $z \in A_n(x)$ we have $c_1 \le \tilde{G}(x,z) \le c_2$.
\end{lemma}

\proof We follow  the proof in \cite{BM10}. Taking $y=z$ in \eqref{e:GBNK} we can improve the upper bound on $G_{B_N \backslash K }(z,z)$ in \cite[(4.10)]{BM10} to $c \log n$. Using Lemma \ref{L:GBNK} again, we can improve the upper bound in the
equation above  \cite[(4.11)]{BM10}, and hence improve the upper bound
in  \cite[(4.11)]{BM10} from $c \log(N/n)/\log N$ to $c /\log n$. With these new bounds the argument of \cite[Lemma 4.6]{BM10} gives that $\tilde{G}(x,z) \le c_2$. \qed

The following two results refine some conditional hitting time estimates from \cite{BM10}.

\begin{lemma}[{cf.\ \cite[(6.1)]{BM10}}] \label{L:neweq6.1} There exists a constant $c_1$ such that if
$D_1 =  \pd_R \BE(x,n/16)$ and $K'=K\backslash\{x\}$, then, for $v \in D_1$,
\[ \pP_v \left( T_x < \tau_{\BE(x,n/8)} \:|\: T_x <T_{K'} \wedge \tau_{B_N} \right) \ge c_1>0.\]
\end{lemma}

\proof Write $B'={B_{n/8}(x)}$. The second displayed equation on \cite[p.\ 2409]{BM10} gives
\be \label{e:fromBM10}
\pP_v \left( T_x < \tau_{B_{n/8}(x)} \:|\: T_x <T_{K'} \wedge \tau_{B_N} \right)
 = \frac{ G_{B'\backslash K}(v,v) }{   G_{B_N\backslash K}(v,v)  }\times
 \frac{  \pP_x ( T_v < \tau_{B'} \wedge T^+_K  ) }{    \pP_x ( T_v < \tau_{B_N} \wedge T^+_K  )},
\ee
where $T^+_K= \min \{ j \ge 1: S_j \in K\}$. As in Lemma \ref{L:GBNK} we have that $G_{B_N\backslash K}(v,v)  \le c \log n$, and so the ratio of Green's functions in \eqref{e:fromBM10} is bounded below by a constant $c>0$. Using the strong Markov property at $\tau_{B'}$ we obtain
\begin{align*}
   \pP_x ( T_v < \tau_{B_N} \wedge T^+_K  )
   &\le \pP_x ( T_v < \tau_{B'} \wedge T^+_K  )
   + \pP_x ( \tau_{B'} \le T^+_K  ) \max_{y \in \pd B'} \pP_y( T_v < \tau_{B_N} \wedge T^+_K ).
  \end{align*}
The argument at the top of \cite[p.\ 2410]{BM10} gives that
$$ \pP_x ( \tau_{B'} \le T^+_K  ) \le c (\log n)  \, \pP_x ( T_v < \tau_{B'} \wedge T^+_K  ) . $$
Moreover, for $y \in \pd B'$,
$$ \pP_y(T_v < \tau_{B_N} \wedge T^+_K )
\le \frac{ G_{\bZ^2 \backslash K}(y,v) }{  G_{\bZ^2\backslash K}(v,v) }, $$
and as in Lemma \ref{L:GBNK} we have $G_{\bZ^2 \backslash K}(y,v) \le c$,
$G_{\bZ^2 \backslash  K}(v,v) \ge c \log n$.
Combining these estimates concludes the proof. \qed

\begin{lemma}[{cf.\ \cite[(6.2)]{BM10}}] \label{L:neweq6.2}
There exists a constant $c>0$ such that if $w \in\pd_R \BE(x,n)$, then
\be \label{e:new6.2}
 \pP_w \left( \tau_{B_N} < T_{\BE(x,7n/8)} \:\vline\:  \tau_{B_N} < T_K \right) \ge  c.\ee
\end{lemma}
\proof As on \cite[p.\ 2410]{BM10}, we let $y_0$ be the point in $B_n(x)$ that maximises $\pP_y(  \tau_{B_N} < T_K )$.
Writing $B_7= \BE(x,7n/8)$, $T_7= T_{B_7}$,  we have
\begin{align*}
 \pP_{y_0} (  \tau_{B_N} < T_K )
 &=   \pP_{y_0}( \tau_{B_N} < T_K \wedge T_7 )
 + \bE^{y_0} ( \mathbf{1}_{\{T_7 < \tau_{B_N} \wedge T_K\} } \pP_{S_{T_7}}( \tau_{B_N} < T_K ) ) \\
&\le  \pP_{y_0}( \tau_{B_N} < T_K \wedge T_7 ) +     \max_{v \in \pd B_7} \pP_v( \tau_{B_N} < T_K).
\end{align*}
Since $K$ contains a path from $0$ to $x$, the discrete Harnack inequality (again, see \cite[Theorem 6.3.9]{LLRW}, for example)
gives us that there exists a constant $p_1>0$ such that
$$  \pP_v( \tau_{\BE(x,n)} < T_K) \le 1-p_1,\qquad\mbox{for all }v \in \pd B_7.$$
Thus
$$   \pP_{y_0} (  \tau_{B_N} < T_K )
\le \pP_{y_0}( \tau_{B_N} < T_K \wedge T_7 )   + (1-p_1)    \pP_{y_0} (  \tau_{B_N} < T_K ), $$
which proves \eqref{e:new6.2} in the case $w=y_0$. We can now use a reflection argument
as on \cite[p.~2410-2411]{BM10} to obtain the general case. \qed

These estimates now lead to an improved lower bound on the length of a LERW. Recall the definition of the conditioned r.w. $\tilde S$, and set  $L_1 = \sL( \sE_{B_N}(\tilde{S}))$, $L_2 = \sE_{B_n(x)}(L_1)$.

\begin{lemma}[{cf.\ \cite[Lemma 6.1]{BM10}}] \label{L:newL61}
There exists a constant $c>0$ such that, for any $z \in A_n(x)$,
\[ \pP( z \in L_2 ) \ge c n^{\kappa -2}.\]
\end{lemma}

\proof Using Lemmas \ref{L:neweq6.2} and \ref{L:neweq6.1} to replace \cite[(6.1),(6.2)]{BM10}, this follows as in \cite{BM10}.
\qed

\begin{proposition} [{cf.\ \cite[Proposition 6.2 and 6.3]{BM10}}]
There exist constants $c_1,c_2$ and $p>0$ such that
\begin{align} \nn
   &c_1 n^\kappa \le \bE M  \le c_2  n^\kappa , \\
   \nn
   &\bE (M^2 ) \le c_2 n^{2\kappa}, \\
    \label{e:Mlowertail}
   &\pP( M \le c_3 n^\kappa ) \le 1-p.
\end{align}
\end{proposition}

\proof Given Lemmas \ref{L:newL61} and \ref{L:newL46}
the bounds on $\bE(M)$ and $\bE(M^2)$  follow as in \cite{BM10}.
The final inequality is then immediate from a second moment bound. \qed

\sm {\em Proof of Theorem \ref{T:LERW-lb}.} We follow the proof of \cite[Proposition 6.6]{BM10}, first proving the result in the case when $D= B_N(0)$, where $N/2 \le nk \le 3N/4$ for some $k \ge 4$. Set $L = \sL( \sE_{B_N(0)}(S^0))$, and, for $j=1, \dots k$, let $\gam_j = \sE_{B_{j n}(0)}(L)$. Let $x_j$ be the point where $L$ first  exits $B_{j n}(0)$, and $B_j = B_{n}(x_j)$. Let $\alpha_j$ be the path $L$ from $x_{j-1}$ to its first exit from $B_{j-1}$, and $V_j$ be the number of hits by $\alpha_j$ on the set $B_{j-1}$. Let $\sF_j$ be the $\sigma$-field generated by $\gam_j$. Using the domain Markov property for the LERW (see \cite{Law91}) and then \eqref{e:Mlowertail}, we have
\be \label{e:Njub}
 \pP\left( V_j \le c_3 n^\kappa | \sF_{j-1} \right) = \pP\left( M^{\gam_{j-1}}_{(j-1)n,n,N,x_{j-1}} \le c_3 n^\kappa\right)\le1-p.
\ee
Let $\eta_j = \mathbf{1}_{\{V_j \le c_3 n^\kappa\}}$. By \eqref{e:Njub}, $\sum_{j=1}^k \eta_j$ stochastically dominates
a binomial random variable with parameters $k$ and $p$, and so there exists a constant $c>0$ such that
$$ \pP\left(\sum_{j=1}^k \eta_j < \half p k \right) \le e^{-c k}. $$
Setting $L' =  \sE_{B_{nk}(0)}(L)$, we have $ |  L' | \ge c_3 n^\kappa \sum_{j=1}^k \eta_j$,
and thus as $N/2 \le nk \le 3N/4$ we obtain
\[ \pP\left( |L' | < c k^{-1/4} N^\kappa \right)  \le e^{-c k};\]
taking $k = c \lam^{1/(\kappa-1)} = c\lam^4$ this gives the result when $D=B_N$.
Note that the proof above actually gives the lower bound for the length of $L'$ rather than $L$,
so we can use Lemma \ref{L:abscty} with $D_1= B_N$, $D_2=D$ to obtain a lower bound of the same form for
$|\sE_{B_N(0)} (\sL(\sE_{D}(S^0)))|$.
\qed

\noindent \textbf{Acknowledgements.} We thank Zhen-Qing Chen for asking us whether the elliptic Harnack inequalities hold for large scale in this model, which motivated us to work on Section \ref{failEHI}, and for the analogy with reflecting Brownian motion in a planar domain with a slit removed given in Remark \ref{zqrem}. This research was supported by the Research Institute for Mathematical Sciences, an International Joint Usage/Research Center located in Kyoto University.

\bibliographystyle{amsplain}
\bibliography{usth}

\normalsize
\vskip 0.1truein

\noindent {\bf Martin T. Barlow}

\noindent Department of Mathematics,
University of British Columbia, Vancouver, B.C., V6T 1Z2, Canada.

\noindent E-mail: {\tt barlow@math.ubc.ca}

\medskip
\noindent {\bf David A. Croydon}

\noindent  Research Institute for Mathematical Sciences,
Kyoto University, Kyoto 606-8502, Japan.

\noindent E-mail: {\tt croydon@kurims.kyoto-u.ac.jp}

\medskip
\noindent {\bf Takashi Kumagai}

\noindent Research Institute for Mathematical Sciences,
Kyoto University, Kyoto 606-8502, Japan.

\noindent E-mail: {\tt kumagai@kurims.kyoto-u.ac.jp}

\end{document}